\input amssym.def
\input amssym.tex

\def\item#1{\vskip1.3pt\hang\textindent {\rm #1}}


\newskip\litemindent
\litemindent=0.7cm  
\def\Litem#1#2{\par\noindent\hangindent#1\litemindent
\hbox to #1\litemindent{\hfill\hbox to \litemindent
{\ninerm #2 \hfill}}\ignorespaces}
\def\litem{\Litem1}

\tolerance=300
\pretolerance=200
\hfuzz=1pt
\vfuzz=1pt

\hoffset=0in
\voffset=0.5in

\hsize=5.8 true in 
\vsize=9.2 true in
\parindent=25pt
\mathsurround=1pt
\parskip=1pt plus .25pt minus .25pt
\normallineskiplimit=.99pt

\countdef\revised=100
\mathchardef\emptyset="001F 
\chardef\ss="19
\def\3{\ss}
\def\anf{$\lower1.2ex\hbox{"}$}
\def\frac#1#2{{#1 \over #2}}
\def\>{>\!\!>}
\def\<{<\!\!<}

\def\into{\hookrightarrow}
\def\onto{\to\mskip-14mu\to} 
\def\ssssarr{\hbox to 15pt{\rightarrowfill}}
\def\sssarr{\hbox to 20pt{\rightarrowfill}}
\def\ssarr{\hbox to 30pt{\rightarrowfill}}
\def\sarr{\hbox to 40pt{\rightarrowfill}}
\def\arr{\hbox to 60pt{\rightarrowfill}}
\def\larr{\hbox to 60pt{\leftarrowfill}}
\def\Arr{\hbox to 80pt{\rightarrowfill}}
\def\mapdown#1{\Big\downarrow\rlap{$\vcenter{\hbox{$\scriptstyle#1$}}$}}

\def\sssmapright#1{\smash{\mathop{\sssarr}\limits^{#1}}}
\def\ssmapright#1{\smash{\mathop{\ssarr}\limits^{#1}}}
\def\smapright#1{\smash{\mathop{\sarr}\limits^{#1}}}

\def\Ad{\mathop{\rm Ad}\nolimits}

\def\Aut{\mathop{\rm Aut}\nolimits}

\def\coker{\mathop{\rm coker}\nolimits}
\def\Comp{\mathop{\rm Comp}\nolimits}

\def\der{\mathop{\rm der}\nolimits}

\def\Diff{\mathop{\rm Diff}\nolimits}

\def\Ext{\mathop{\rm Ext}\nolimits}

\def\Gau{\mathop{\rm Gau}\nolimits}

\def\Hom{\mathop{\rm Hom}\nolimits}%
\def\id{\mathop{\rm id}\nolimits} 
\def\im{\mathop{\rm im}\nolimits}


%

\def\Out{\mathop{\rm Out}\nolimits}

\def\PSL{\mathop{\rm PSL}\nolimits}



\def\SL{\mathop{\rm SL}\nolimits}



\def\trile{\trianglelefteq}

\def\0{{\bf 0}}
\def\1{{\bf 1}}

\def\a{{\frak a}}

\def\b{{\frak b}}

\def\g{{\frak g}}

\def\h{{\frak h}}

\def\m{{\frak m}}

\def\n{{\frak n}}

\def\z{{\frak z}}

\def\L{\mathop{\bf L{}}\nolimits}

\def\C{{{\Bbb C}{\mskip+1mu}}} 
\def\K{{{\Bbb K}{\mskip+2mu}}} 

\def\R{{\Bbb R}} 
\def\Z{{\Bbb Z}} 
\def\N{{\Bbb N}}

\def\K{{\Bbb K}}

\def\:{\colon}  
\def\.{{\cdot}}
\def\|{\Vert}
\def\bsk{\bigskip}

\def\giantskip{\vskip2\bigskipamount}
\def\gsk{\giantskip}
\def \la {\langle}

\def\msk{\medskip}
\def \ra {\rangle}
\def \res {\!\mid\!\!}

\def\bbr{\bigbreak}
\def\giantbreak{\par \ifdim\lastskip<2\bigskipamount \removelastskip
         \penalty-400 \giantskip\fi}

\def\nin{\noindent}
\def\cen{\centerline}
\def\pagebreak{\vskip 0pt plus 0.0001fil\break}
\def\linebreak{\break}

\def\hat{\widehat}

\def\epsilon{\varepsilon}
\def\eset{\emptyset}

\def\nin{\noindent}
\def\oline{\overline}

\def\ubrace{\underbrace}
\def\pder#1,#2,#3 { {\partial #1 \over \partial #2}(#3)}
\def\pde#1,#2 { {\partial #1 \over \partial #2}}
\def\phi{\varphi}


\def\subeq{\subseteq}

\def\Rarrow{\Rightarrow}

\def\tilde{\widetilde}

\font\ninerm=cmr9
\font\eightrm=cmr8

\font\eightbf=cmbx8


\font\smc=cmcsc10
\font\bfone=cmbx10 scaled\magstep1 
\font\bftwo=cmbx10 scaled\magstep2 

\def\qed{{\unskip\nobreak\hfil\penalty50\hskip .001pt \hbox{}\nobreak\hfil
          \vrule height 1.2ex width 1.1ex depth -.1ex
           \parfillskip=0pt\finalhyphendemerits=0\medbreak}\rm}

\def\qeddis{\eqno{\vrule height 1.2ex width 1.1ex depth -.1ex} $$
                   \medbreak\rm}

\def\Lemma #1. {\bigbreak\vskip-\parskip\noindent{\bf Lemma #1.}\quad\it}

\def\Sublemma #1. {\bigbreak\vskip-\parskip\noindent{\bf Sublemma #1.}\quad\it}

\def\Proposition #1. {\bigbreak\vskip-\parskip\noindent{\bf Proposition #1.}
\quad\it}

\def\Corollary #1. {\bigbreak\vskip-\parskip\nin{\bf Corollary #1.}
\quad\it}

\def\Theorem #1. {\bigbreak\vskip-\parskip\noindent{\bf Theorem #1.}
\quad\it}

\def\Definition #1. {\rm\bigbreak\vskip-\parskip\noindent
{\bf Definition #1.}
\quad}

\def\Remark #1. {\rm\bigbreak\vskip-\parskip\noindent{\bf Remark #1.}\quad}

\def\Example #1. {\rm\bigbreak\vskip-\parskip\noindent{\bf Example #1.}\quad}
\def\Examples #1. {\rm\bigbreak\vskip-\parskip\noindent{\bf Examples #1.}\quad}

\def\Problems #1. {\bigbreak\vskip-\parskip\noindent{\bf Problems #1.}\quad}
\def\Problem #1. {\bigbreak\vskip-\parskip\noindent{\bf Problem #1.}\quad}
\def\Exercise #1. {\bigbreak\vskip-\parskip\noindent{\bf Exercise #1.}\quad}

\def\Conjecture #1. {\bigbreak\vskip-\parskip\noindent{\bf Conjecture #1.}\quad}

\def\Proof#1.{\rm\par\ifdim\lastskip<\bigskipamount\removelastskip\fi\smallskip
            \noindent {\bf Proof.}\quad}

\def\Axiom #1. {\bigbreak\vskip-\parskip\noindent{\bf Axiom #1.}\quad\it}

\def\Satz #1. {\bigbreak\vskip-\parskip\noindent{\bf Satz #1.}\quad\it}

\def\Korollar #1. {\bbr\vskip-\parskip\nin{\bf Korollar #1.} \quad\it}

\def\Folgerung #1. {\bbr\vskip-\parskip\nin{\bf Folgerung #1.} \quad\it}

\def\Folgerungen #1. {\bbr\vskip-\parskip\nin{\bf Folgerungen #1.} \quad\it}

\def\Bemerkung #1. {\rm\bigbreak\vskip-\parskip\noindent{\bf Bemerkung #1.}
\quad}

\def\Beispiel #1. {\rm\bigbreak\vskip-\parskip\noindent{\bf Beispiel #1.}\quad}
\def\Beispiele #1. {\rm\bigbreak\vskip-\parskip\noindent{\bf Beispiele #1.}\quad}
\def\Aufgabe #1. {\rm\bigbreak\vskip-\parskip\noindent{\bf Aufgabe #1.}\quad}
\def\Aufgaben #1. {\rm\bigbreak\vskip-\parskip\noindent{\bf Aufgabe #1.}\quad}

\def\Beweis#1. {\rm\par\ifdim\lastskip<\bigskipamount\removelastskip\fi
           \smallskip\noindent {\bf Beweis.}\quad}

\nopagenumbers

\def\date{\ifcase\month\or January\or February \or March\or April\or May
\or June\or July\or August\or September\or October\or November
\or December\fi\space\number\day, \number\year}

\def\title{Title ??}
\def\author{Author ??}

\def\thanks#1{\footnote*{\eightrm#1}}

\def\rightheadline{\hfil{\eightrm\title}\hfil\tenbf\folio}
\def\leftheadline{\tenbf\folio\hfil{\eightrm\author}\hfil}
\headline={\vbox{\line{\ifodd\pageno\rightheadline\else\leftheadline\fi}}}

\def\firstheadline{}
\def\firstfootline{\cen{\rm\folio}}

\def\seite #1 {\pageno #1
               \headline={\ifnum\pageno=#1 \firstheadline
               \else\ifodd\pageno\rightheadline\else\leftheadline\fi\fi}
               \footline={\ifnum\pageno=#1 \firstfootline\else{}\fi}}

\newdimen\dimenone
 \def\checkleftspace#1#2#3#4{
 \dimenone=\pagetotal
 \advance\dimenone by -\pageshrink   
 \ifdim\dimenone>\pagegoal          
   \else\dimenone=\pagetotal
        \advance\dimenone by \pagestretch
        \ifdim\dimenone<\pagegoal
          \dimenone=\pagetotal
          \advance\dimenone by#1         
          \setbox0=\vbox{#2\parskip=0pt                
                     \hyphenpenalty=10000
                     \rightskip=0pt plus 5em
                     \noindent#3 \vskip#4}    
        \advance\dimenone by\ht0
        \advance\dimenone by 3\baselineskip   
        \ifdim\dimenone>\pagegoal\vfill\eject\fi
          \else\eject\fi\fi}


\def\subheadline #1{\nin\bigbreak\vskip-\lastskip
      \checkleftspace{0.9cm}{\bf}{#1}{\medskipamount}
          \indent\vskip0.7cm\centerline{\bf #1}\medskip}
\def\subsection{\subheadline} 

\def\lsubheadline #1 #2{\bigbreak\vskip-\lastskip
      \checkleftspace{0.9cm}{\bf}{#1}{\bigskipamount}
         \vbox{\vskip0.7cm}\cen{\bf #1}\msk \cen{\bf #2}\bsk}

\def\sectionheadline #1{\bigbreak\vskip-\lastskip
      \checkleftspace{1.1cm}{\bf}{#1}{\bigskipamount}
         \vbox{\vskip1.1cm}\cen{\bfone #1}\bsk}
\def\section{\sectionheadline} 

\def\lsectionheadline #1 #2{\bigbreak\vskip-\lastskip
      \checkleftspace{1.1cm}{\bf}{#1}{\bigskipamount}
         \vbox{\vskip1.1cm}\cen{\bfone #1}\msk \cen{\bfone #2}\bsk}

\def\lchapterheadline #1 #2{\bigbreak\vskip-\lastskip\indent\vskip3cm
                       \cen{\bftwo #1} \msk \cen{\bftwo #2} \gsk}
\def\llsectionheadline #1 #2 #3{\bigbreak\vskip-\lastskip\indent\vskip1.8cm
\cen{\bfone #1} \msk \cen{\bfone #2} \msk \cen{\bfone #3} \nobreak\bsk\nobreak}


\newtoks\literat
\def\[#1 #2\par{\literat={#2\unskip.}%
\hbox{\vtop{\hsize=.15\hsize\nin [#1]\hfill}
\vtop{\hsize=.82\hsize\nin\the\literat}}\par
\vskip.3\baselineskip}

\def\references{
\sectionheadline{\bf References}
\frenchspacing

\entries\par}

\mathchardef\emptyset="001F 
\def\address{Author: \tt$\backslash$def$\backslash$address$\{$??$\}$}

\def\abstract #1{{\narrower\baselineskip=10pt{\noindent
\eightbf Abstract.\quad \eightrm #1 }
\bigskip}}

\def\firstpage{\nin
{\obeylines \parindent 0pt }
\vskip2cm
\centerline{\bfone\title}
\gsk
\centerline{\bf\author}
\vskip1.5cm \rm}

\def\lastpage{\par\vbox{\vskip1cm\nin
\line{
\vtop{\hsize=.5\hsize{\parindent=0pt\baselineskip=10pt\nin\address}}
\hfill} }}

\def\Box #1 { \msk\par\nin 
\centerline{
\vbox{\offinterlineskip
\hrule
\hbox{\vrule\strut\hskip1ex\hfil{\smc#1}\hfill\hskip1ex}
\hrule}\vrule}\msk }

\def\adots{\mathinner{\mkern1mu\raise1pt\vbox{\kern7pt\hbox{.}}
                        \mkern2mu\raise4pt\hbox{.}
                        \mkern2mu\raise7pt\hbox{.}\mkern1mu}}


\pageno=1

\def\title{Non-abelian extensions of infinite-dimensional Lie groups} 
\def\author{Karl-Hermann Neeb}
\def\date{April 14, 2005} 
\def\rightheadline{\tenbf\folio\hfil{\tt extliegrp.tex}\hfil\eightrm\date}
\def\leftheadline{\tenbf\folio\hfil{\rm\title}\hfil\eightrm\date}

\def\ubrace{\underbrace}

\def\pr{{\rm pr}}

\def\address
{Karl-Hermann Neeb

Technische Universit\"at Darmstadt 

Schlossgartenstrasse 7

D-64289 Darmstadt 

Deutschland

neeb@mathematik.tu-darmstadt.de}

\firstpage

\def\sst{\scriptstyle} 

\abstract{In this paper we study non-abelian extensions of a 
Lie group ${\scriptstyle G}$ modeled on a locally convex space by 
a Lie group ${\scriptstyle N}$. The equivalence classes of such extension 
are grouped into those corresponding to a class of so-called smooth 
outer actions ${\scriptstyle S}$ of ${\sst G}$ on ${\sst N}$. If ${\scriptstyle S}$ is given, 
we show that the corresponding set ${\scriptstyle \Ext(G,N)_S}$ of extension 
classes is a principal homogeneous space of the locally smooth cohomology group 
${\scriptstyle H^2_{ss}(G,Z(N))_S}$. To each ${\sst S}$ 
a locally smooth 
obstruction class ${\sst \chi(S)}$ 
in a suitably defined cohomology group ${\scriptstyle H^3_{ss}(G,Z(N))_S}$ is defined. It 
vanishes if and only if there is a corresponding extension of 
${\sst G}$ by ${\sst N}$. 
A central point is that we reduce many problems concerning extensions by 
non-abelian groups to questions on extensions by abelian groups, which 
have been dealt with in previous work. An important tool is a 
Lie theoretic concept of a smooth crossed module ${\sst \alpha \: H \to G}$, which 
we view as a central extension of a normal subgroup of ${\sst G}$. 
\hfill\break 
Keywords: Lie group extension, smooth outer action, crossed module, Lie group cohomology, 
automorphisms of group extension \hfill\break 
}

\sectionheadline{Introduction} 

\nin In the present paper we continue our investigation of extensions of infinite-dimensional Lie groups. 
In this sense it is a sequel to [Ne02] and [Ne04a], where we studied central, resp., abelian 
extensions. We now turn to extensions by groups which are not necessarily 
abelian. 

The concept of a (not necessarily finite-dimensional) Lie group used here 
is that a {\it Lie group} $G$ is a manifold modeled on a locally convex space  
endowed with a group structure for which the group operations are smooth (cf.\ [Mil83], 
[Gl01]). 
An {\it extension} is an exact sequence of Lie groups 
$N \into \hat G \onto G$ which defines a locally trivial smooth 
$N$-principal bundle over the Lie group $G$. It is called 
{\it abelian}, resp., {\it central} if $N$ is abelian, resp., central in $\hat G$. 
This setup implies in particular that for each extension 
$q \: \hat G \to G$ of $G$ by $N$ there exists a section 
$\sigma \: G \to \hat G$ of $q$ mapping $\1$ to $\1$ which is smooth in 
an identity neighborhood of $G$. Starting with such a section, we consider the 
functions 
$$ S \: G \to \Aut(N), \quad S(g)(n) := \sigma(g)n\sigma(g)^{-1}, \quad 
\omega \: G \times G \to N, \ \ (g,g') \mapsto \sigma(g)\sigma(g')\sigma(gg')^{-1}. $$
Then the map $N \times G \to \hat G, (n,g) \mapsto n \sigma(g)$ is a bijection,  
and transferring the group structure from $\hat G$ to $N \times G$ yields 
$$ (n.g)(n'g') = (nS(g)(n') \omega(g,g'), gg'). \leqno(\dagger)$$

In the first half of the present paper we are mainly concerned with the appropriate 
smoothness requirements for the functions $S$ and $\omega$ under which the product set 
$N \times G$ with the group structure given by ($\dagger$) carries a Lie group structure 
for which $N \times G \to G, (n,g) \mapsto g$ is an extension of Lie groups, and 
when such Lie group extensions are equivalent. The natural context to deal with 
such questions is a non-abelian locally smooth Lie group cohomology introduced 
in Section~I. In Section~II the cohomological setup is linked to 
the classification of Lie group extensions of $G$ by~$N$. In particular 
we show that all such extensions are constructed from pairs of 
maps $(S,\omega)$ as above, where 
\litem{(1)} $S(g) S(g') = c_{\omega(g,g')} S(gg')$ for $g,g' \in G$ and 
$c_n(n') := nn'n^{-1}$ denotes conjugation. 
\litem{(2)} $S(g)(\omega(g',g'')) \omega(g,g'g'') = \omega(g,g') \omega(gg',g'')$ 
for $g,g',g'' \in G$. 
\litem{(3)} The map $G \times N \to N, (g,n) \mapsto S(g)n$ is smooth on a set 
of the form $U \times N$, where $U$ is an identity neighborhood of $G$. 
\litem{(4)} $\omega$ is smooth in an identity neighborhood. 
\litem{(5)} For each $g \in G$ the map 
$\omega_g \: G \to N, x \mapsto 
\omega(g,x) \omega(gxg^{-1},g)^{-1}$
is smooth in an identity neighborhood of $G$. 

Conditions (1) and (2) ensure that $N \times G$ is a group with the multiplication defined by 
$(\dagger)$, and (3)-(5) imply that $N \times G$ can be endowed with a manifold structure 
turning it into a Lie group. In general this manifold will not be 
diffeomorphic to the product manifold $N \times G$, but the projection onto $G$ 
will define a smooth $N$-principal bundle. A map $S \: G \to \Aut(N)$ satisfying 
(3) and for which there is an $\omega$ satisfying (1) is called a {\it smooth outer action}  
of $G$ on~$N$. 

In the setting of abstract groups, the corresponding results are due to 
Eilenberg and MacLane ([EML47]). In this context it is more convenient to work with the 
homomorphism $s \: G \to \Out(N), g \mapsto [S(g)]$ defined by $S$. They show that if there 
is some extension corresponding to $s$, then 
all $N$-extensions of $G$ corresponding to $s$ form a principal homogeneous space 
of the cohomology 
group $H^2(G,Z(N))_s$, where $Z(N)$ is endowed with the $G$-module 
structure determined by~$s$. 
In loc.cit.\ it is also shown that, for a given $s$, there is a characteristic class 
in $H^3(G,Z(N))_s$ vanishing if and only if $s$ corresponds to a group extension. 
In Section~II we adapt these results to the Lie group setting 
in the sense that for a smooth outer action $S$, we show that 
a certain characteristic class $\chi(S)$ in the Lie group cohomology $H^3_{ss}(G,Z(N))_S$ 
determines whether there exists a corresponding Lie group extension, i.e., 
there exists another choice of $\omega$ such that also (2) is satisfied. 
We further show that in the latter case the set $\Ext(G,N)_S$ of $N$-extensions of 
$G$ corresponding to $S$ is a principal homogeneous space of the cohomology group 
$H^2_{ss}(G,Z(N))_S$ defined in Section~I. 

A particular subtlety entering the picture in our Lie theoretic context is that 
in general the center $Z(N)$ of $N$ need not be a Lie group. For the results of 
Sections~I and II it will suffice that $Z(N)$ carries the structure of an initial 
Lie subgroup, which is a rather weak requirement. For the definition of the 
$Z(N)$-valued cochain spaces this property is not needed 
because we use a definition of $H^2_{ss}(G,Z(N))_S$ and 
$H^3_{ss}(G,Z(N))_S$ that does not require a manifold structure on $Z(N)$ since 
we consider $Z(N)$-valued cochains as $N$-valued functions to specify smoothness 
requirements. 

In Section~III we introduce crossed modules for Lie groups and discuss their 
relation to group extensions. A crossed module is a morphism 
$\alpha \: H \to G$ of Lie groups for which $\im(\alpha)$ and $\ker(\alpha)$ 
are split Lie subgroups, together with an action of $G$ on $H$ lifting the conjugation 
action of $G$ and extending the conjugation action of $H$ on itself. 
Then $\ker(\alpha) \into H \onto \im(\alpha)$ is a central extension of 
$N := \im(\alpha)$ by $Z := \ker(\alpha)$ and the question arises whether this 
$Z$-extension of $N$ can be enlarged in an equivariant fashion to an extension of 
$G$ by $Z$. As in Section~II, it turns out that such an enlargement exists if and only if 
a certain cohomology class in $H^3_{ss}(G,Z)_\alpha$ vanishes (Theorem~III.8). 

If the group $N$ is abelian, then the data given by a smooth outer action 
$S$ of $G$ on $N$ is the same as a smooth $G$-module structure on~$N$. 
In Section~IV we start with this data and the assumption that 
$Z(N)$ is a split Lie subgroup of $N$, so that we have a central Lie group extension 
$Z(N) \into N \onto N_{\rm ad} := N/Z(N)$. From that we construct a Lie group 
$G^S$ which is an extension of $G$ by $N_{\rm ad}$, depending on $S$, 
which has the interesting property that all $N$-extensions of $G$ corresponding to 
$S$ are extensions of the Lie group $G^S$ by $Z(N)$. Moreover, the map 
$N \to N_{\rm ad} \subeq G^S$ together with the natural action of $G^S$ on $N$ 
defines a crossed module of Lie groups whose characteristic class coincides with the 
characteristic class $\chi(S) \in H^3(G,Z(N))_S$. 

These observations reduce many 
questions on general Lie group extensions to the special case of abelian 
extensions which has been treated in detail in [Ne04a]. In particular,  
the question whether a topologically split extension 
$\n \into \hat\g \onto \g$ of Lie algebras is integrable if we assume that we 
already have 
groups $G$ and $N$, together with a compatible smooth outer action of $G$ on $N$, 
can be reduced to the question whether the corresponding abelian extension 
$\z(\n) \into \g^S \onto \g$ is integrable, which is characterized in [Ne04a] by the 
discreteness of the image of a period homomorphism 
$\pi_2(G^S) \to \z(\n)$.  In [Ne04a] we also describe an exact sequence which 
explains how the cohomology group $H^2_s(G,Z(N))_S$ is given in terms of 
topological data associated to the group $G$ and 
the corresponding Lie algebra cohomology space $H^2_c(\g,\z(\n))_S$. 

Non-abelian extensions of Lie groups occur quite naturally in the context of 
smooth principal bundles over compact manifolds. If $q \: P \to M$ is a smooth 
$K$-principal bundle over the compact manifold, where $K$ is a finite-dimensional Lie group, 
then the group $\Aut(P)$ is a Lie group extension 
$$ \Gau(P) \into \Aut(P) \onto \Diff(M)_{[P]}, $$
where $\Gau(P)$ is the gauge group of $P$ and $\Diff(M)_{[P]}$ is the open subgroup 
of $\Diff(M)$ consisting of all diffeomorphisms preserving the bundle class 
$[P]$ under pull-backs. If $K$ is abelian, then 
$\Gau(P) \cong C^\infty(M,K)$ and we have an abelian extension. 

In an appendix we collect several results on automorphisms of group extensions 
and on smooth Lie group actions on Lie group extensions. Although 
only very specific results of this appendix are used in the present paper, we felt 
that they are very useful in many different contexts, so that it makes good sense to 
collect them in such an appendix. 
For a discussion of automorphisms of abstract groups from a cohomological point of 
view, we refer to Huebschmann's paper [Hu81a]. 

\gsk

The problem to parameterize all extensions of a given group $G$ by a group $N$ is 
a core problem in group theory. It seems that the first systematic treatment 
using a parameterization of the extended group $\hat G$ as a product set 
$N \times G$ is due to Schreier ([Sch26a/b]). He mainly works with presentations of the 
groups involved. Baer discovered the first examples of outer actions for which 
no group extensions exist ([Baer34]), and Turing refined the methods introduced 
by Baer for abelian extensions to extensions by non-abelian groups ([Tu38]). 
In particular 
he discusses extensions of a finite cyclic group by a group~$N$. Cohomological methods 
to study group extensions first appear in the work of Eilenberg and MacLane 
([EML42/47]) who coin the term {\it $G$-kernel} for a pair $(N,s)$ consisting of a group 
$N$ and a homomorphism $s \: G \to \Out(N)$. They show that all extensions corresponding to a 
given kernel can be parametrized by the cohomology group $H^2(G,Z(N))_s$ and that 
the obstruction of a kernel to correspond to a group extension is an element in 
$H^3(G,Z(N))_s$. They also prove for a given $G$-module $Z$ and each cohomology 
class $c \in H^3(G,Z)$ the existence of some 
$G$-kernel $(N,s)$ with $Z(N) = Z$ whose obstruction 
class is $c$. In loc.cit.\ one even finds the reduction to abelian extensions in the sense 
that $N$-extensions of $G$ corresponding to $s$ can also be described as 
$Z(N)$-extensions of a group $G^s$, depending only on $s$. 

Non-abelian extensions of topological groups have been studied by Brown  under the 
assumptions that $G$ and $N$ are polonais ([Br71]). 
Then the appropriate group cohomology is defined by 
measurable cocycles, and there is a natural action of the group $C^1(G,N)$ of measurable 
cochains on a certain set of pairs $(S,\omega)$, also 
satisfying certain measurability conditions, 
such that the orbits $[(S,\omega)]$ of $C^1(G,N)$ are in one-to-one correspondence with 
the equivalence classes of topological extensions of $G$ by $N$ 
(cf.\ our Corollary~II.11). 

For a survey of the relations between group cohomology and non-abelian group extensions 
with twisted actions and crossed product algebras we refer to [RSW00]. Non-abelian 
group extensions also arise naturally in mathematical physics ([CFNW94]), where the 
groups involved are mostly infinite-dimensional Lie groups, although they are often 
dealt with by ad hoc methods, and not in the framework of a general theory, to which 
we contribute with the present paper. 

In a subsequent paper we plan to return to this 
topic and address the question how to describe elements of 
$H^3_s(G,Z)$, where $Z$ is a smooth $G$-module, in terms of Lie algebra cohomology classes in 
$H^3_c(\g,\z)$ and topological data associated to $G$ and $Z$, a topic 
that has been dealt with in degree $2$ in [Ne04a]. Another interesting question 
is how to realize classes in $H^3_s(G,Z)$ as characteristic classes of 
smooth crossed modules $\alpha \: H_1 \to H_2$ with $\ker(\alpha) \cong Z$ 
and $\coker(\alpha) \cong G$. Both questions are of particular interest 
if $G$ is a group of diffeomorphisms of a compact manifold $M$ and 
$\z$ is the quotient of the space $\Omega^p(M,\R)$ of smooth $p$-forms on $M$ 
modulo the subspace of exact forms. Further, it should be possible to 
describe the criteria for the integrability of the abelian Lie algebra extensions 
$\z(\n) \into \hat\g \onto \g^S$ arising in Section~IV directly in terms of 
data associated to $G$, $N$ and the smooth outer action~$S$. 

\section{I. Lie groups and their cohomology} 

\nin In this section we introduce notation and terminology used throughout this paper. 
In particular we describe a natural cohomology for Lie groups also used in [Ne04a] 
in the context of abelian extensions and provide some basic 
results and calculations that will be crucial for the extension theory. 

\subheadline{Infinite-dimensional Lie groups} 

{\bf Infinite-dimensional manifolds:} In 
this paper $\K \in \{\R,\C\}$ denotes the field of real or complex
numbers. Let $X$ and $Y$ be topological $\K$-vector spaces, $U
\subeq X$ open and $f \: U \to Y$ a map. Then the {\it derivative
  of $f$ at $x$ in the direction of $h$} is defined as 
$$ df(x)(h) := \lim_{t \to 0} {1 \over t} \big( f(x + t h) - f(x)\big)
$$
whenever the limit exists. The function $f$ is called {\it differentiable at
  $x$} if $df(x)(h)$ exists for all $h \in X$. It is called {\it
  continuously differentiable or $C^1$} if it is continuous, differentiable at all
points of $U$,  and 
$$ df \: U \times X \to Y, \quad (x,h) \mapsto df(x)(h) $$
is a continuous map. It is called a {\it $C^n$-map} if $f$ is $C^1$ and $df$ is a
$C^{n-1}$-map, and $C^\infty$ ({\it smooth}) if it is $C^n$ for all $n \in \N$. 
This is the notion of differentiability used in [Mil83] and
[Gl01], where the latter reference deals with the modifications
necessary for incomplete spaces. Since we have a chain rule for $C^1$-maps between locally convex 
spaces, we can define smooth manifolds as in
the finite-dimensional case. 

{\bf Lie groups:} A {\it Lie group} $G$ is a smooth manifold modeled on a locally convex space 
for which the group multiplication and the 
inversion are smooth maps. We refer to [Mil83] for basic results on Lie groups in this context. 
We write $\1 \in G$ for the identity
element, 
$$ \lambda_g(x) := gx, \quad \rho_g(x) := xg \quad \hbox{ and } \quad 
c_g(x) := C_G(g)(x) := gxg^{-1} $$
for left and right translations, resp., conjugation on $G$. 
The group of all Lie group
automorphisms of $G$ is denoted $\Aut(G)$, we write 
$C_G(G) = \{ c_g \: g \in G\}$ for the
normal subgroup of {\it inner automorphisms}, and $\Out(G) :=
\Aut(G)/C_G(G)$ for the quotient group of {\it outer automorphisms}. 
The {\it conjugation action of $G$} is the homomorphism $C_G \: G \to \Aut(G)$ with
range $C_G(G)$. For a normal subgroup $N \trile G$ we write $C_N \: G \to \Aut(N)$ with 
$C_N(g) := C_G(g)\res_N$. If $G$ is a connected 
Lie group, then we write $q_G \: \tilde G \to G$ for its universal covering Lie 
group and identify $\pi_1(G)$ with the kernel of $q_G$.

{\bf Lie subgroups:} A subgroup $H \leq G$ of a Lie group $G$ is called {\it an initial Lie subgroup} 
if it carries a Lie group structure for which the inclusion map 
$\iota \: H \into G$ is an {\it immersion} in the sense that the tangent map 
$T\iota \: TH \to TG$ is fiberwise injective, and for each smooth map 
$f \: M \to G$ from a smooth manifold $M$ to $G$ with $\im(f) \subeq H$ the corresponding map 
$\iota^{-1} \circ f \: M \to H$ is smooth. 
The latter property implies in particular 
that the initial Lie group structure on $H$ is unique because 
if $\iota' \: H' \into G$ also is a smooth injective immersion of a Lie group 
with $\iota'(H') = H$ and such that for each smooth map $f \: M \to G$ with $\im(f) \subeq H$ 
the map $\iota'^{-1} \circ f \: M \to H'$ is smooth, then 
in particular the maps $\iota^{-1} \circ \iota' \: H' \to H$ and 
$\iota'^{-1} \circ \iota \: H \to H'$ are smooth 
morphism of Lie groups. Since these maps are each others inverse, 
the Lie groups $H$ and $H'$ are isomorphic. 

{\bf Lie algebras of Lie groups:} Each $x \in T_\1(G)$ corresponds to
a unique left invariant vector field $x_l$ with 
$g.x := x_l(g) := d\lambda_g(\1).x, g \in G.$
The space of left invariant vector fields is closed under the Lie
bracket of vector fields, hence inherits a Lie algebra structure. In
this sense we obtain on the locally convex space $T_\1(G) \cong \g$ 
a continuous Lie bracket which is uniquely determined by the relation $[x,y]_l = [x_l, y_l]$. 

We call a Lie algebra $\g$ which is a topological vector space 
such that the Lie bracket is continuous a {\it topological Lie
  algebra $\g$}. In this sense the Lie algebra $\L(G) := (T_\1(G),[\cdot,\cdot])$ of a Lie
group $G$ is a locally convex topological Lie algebra.  

{\bf Topological splitting:} A linear subspace $W$ of a topological vector space $V$ is called 
{\it (topologically) split} if it is closed and there is a continuous linear map 
$\sigma \: V/W \to V$ for which the map 
$$ W \times V/W \to V, \quad (w,x) \mapsto w + \sigma(x) $$
is an isomorphism of topological vector spaces. Note that the closedness of 
$W$ guarantees that the quotient topology turns $V/W$ into a Hausdorff space which 
is a topological vector space with respect to 
the induced vector space structure. 
A continuous linear map $f \: V \to W$ between topological vector spaces  
is said to be {\it (topologically) split} if the subspaces 
$\ker(f) \subeq V$ and $\im(f) \subeq W$ are topologically split.

\subheadline{Locally smooth cohomology of Lie groups} 

\Definition I.1. If $G$ and $N$ are Lie groups, then we 
put $C^0_s(G,N) := N$ and for 
$p \in \N$ we call a map $f \: G^p \to N$ {\it locally smooth}
if there exists an open identity neighborhood $U \subeq G^p$ such that
$f\res_U$ is smooth. We say that $f$ is {\it normalized} if 
$$ (\exists j)\ g_j = \1 \quad \Rarrow \quad f(g_1,\ldots, g_p) =\1. $$
We write 
$C^p_s(G,N)$ for the space of all normalized locally smooth maps $G^p
\to N$, the so-called (locally smooth) {\it $p$-cochains}. 

We shall see below that for $p  = 2$ and non-connected Lie groups, we sometimes have to 
require additional smoothness: We write $C^2_{ss}(G,N)$ for the 
set of all elements $\omega \in C^2_s(G,N)$ with the additional property that for 
each $g \in G$ the map 
$$ \omega_g \: G \to N, \quad x \mapsto 
\omega(g,x) \omega(gxg^{-1},g)^{-1} $$
is smooth in an identity neighborhood of $G$. Note that $\omega_g(\1) = \1$. 

For $f \in C^1_s(G,N)$ we define 
$$ \delta_f \: G \times G \to N, \quad \delta_f(g,g') := f(g)f(g') f(gg')^{-1} $$
and observe that 
$$ \eqalign{ (\delta_f)_g(x) 
&=  \delta_f(g,x) \delta_f(gxg^{-1},g)^{-1} 
=  f(g)f(x)f(gx)^{-1} f(gx)f(g)^{-1}f(gxg^{-1})^{-1} \cr
&=  f(g)f(x)f(g^{-1}) f(gxg^{-1})^{-1} \cr} $$
is smooth in an identity neighborhood of $G$, so that 
$\delta_f \in C^2_{ss}(G,N).$

If $G$ and $N$ are abstract groups, we 
write $C^p(G,N)$ for the set of all functions $G^p \to N$. 
\qed

For the following we also observe that for $h, f \in C^1_s(G,N)$ we have 
$$ \delta_{hf}(g,g') 
= h(g) f(g) h(g') f(g') f(gg')^{-1} h(gg')^{-1} 
= h(g) c_{f(g)}(h(g')) \delta_f(g,g') h(gg')^{-1}. \leqno(1.1) $$ 

For cochains with values in abelian groups, resp., smooth modules, the following 
definition recalls the natural Lie group cohomology setting from [Ne04a]. 

\Definition I.2. Let $G$ be a Lie group and $A$ a {\it smooth $G$-module}, 
i.e., a pair $(A,S)$ of an  abelian Lie group $A$ and 
a smooth $G$-action $G \times A \to A$ given by the homomorphism $S \: G \to \Aut(A)$. 
We call the elements of 
$C^p_s(G,A)$ {\it (locally smooth normalized) $p$-cochains} and write 
$$ d_S \: C^p_s(G,A) \to C^{p+1}_s(G,A) $$
for the group differential given by 
$$ \eqalign{ 
&\ \ \ \  (d_S f)(g_0, \ldots, g_p) \cr
&:= S(g_0)(f(g_1, \ldots, g_p)) 
+ \sum_{j = 1}^p (-1)^j f(g_0, \ldots, g_{j-1} g_j,\ldots, g_p) 
+ (-1)^{p+1} f(g_0, \ldots, g_{p-1}). \cr} $$
It is easy to verify that $d_S(C^p_s(G,A)) \subeq
C^{p+1}_s(G,A)$. We thus obtain a sub-complex of the standard group cohomology 
complex $(C^\bullet(G,A),d_S)$. In view of $d_S^2 = 0$, the space 
$Z^p_s(G,A)_S := \ker d_S\res_{C^p_s(G,A)}$
of {\it $p$-cocycles} contains the space 
$B^p_s(G,A)_S :=  d_S(C^{p-1}_s(G,A))$
of {\it $p$-coboundaries}. The quotient 
$$ H^p_s(G,A)_S := Z^p_s(G,A)_S/B^p_s(G,A)_S $$
is the {\it $p^{th}$ locally smooth cohomology group of $G$ with values in
the $G$-module $A$}. We write $[f] \in H^p_s(G,A)_S$ for the cohomology class of a
cocycle $f \in Z^p_s(G,A)_S$. 
\qed

In the following we also need generalizations of Lie group cohomology to settings 
where the target groups are non-abelian. For $1$-cocycles the generalization is quite 
straight forward: 

\Definition I.3. Let $G$ and $N$ be Lie groups and 
$S \: G \to \Aut(N)$ a homomorphism defining a smooth action of $G$ on $N$. 
A function $f \in C^1_s(G,N)$ is called a {\it $1$-cocycle} (or a {\it crossed 
homomorphism}) if 
$$ f(gg') = f(g) S(g)(f(g')) \quad \hbox{ for all } \quad g,g' \in G. $$ 
We write $Z^1_s(G,N)_S$ for the set of all $1$-cocycles $G \to N$ with respect to the 
action of $G$ on $N$ given by $S$. The group $N = C^0_s(G,N)$ acts on 
$Z^1_s(G,N)_S$ by 
$$ (n.f)(g) := n f(g) S(g)(n)^{-1}, $$
and the set of orbits of $N$ in $Z^1_s(G,N)_S$ is the {\it cohomology set} 
$$ H^1_s(G,N)_S := Z^1_s(G,N)_S/G. 
\qeddis 

For the extension theory of Lie groups we also need non-abelian $2$-cocycles, 
but for the appropriate concept we have to weaken the setting of a smooth 
action of $G$ on $N$ as follows. 

\Definition I.4. Let $G$ and $N$ be Lie groups. 
We define 
$C^1_s(G, \Aut(N))$ as the set of all maps $S \: G \to \Aut(N)$ with 
$S(\1) = \id_N$ and for which there exists an open identity neighborhood 
$U \subeq G$ such that the map 
$$U \times N \to N,\quad (g,n) \mapsto S(g)(n) $$ 
is smooth. 

We call a  map $S \in C^1_s(G,\Aut(N))$ 
a {\it smooth outer action of $G$ on $N$} if  there exists 
$\omega \in C^2_{s}(G,N)$ with $\delta_S = C_N \circ \omega$. 
It is called {\it strongly smooth} if $\omega$ can be chosen in the smaller set 
$C^2_{ss}(G,N)$. 

On the set of smooth outer actions we define an equivalence relation by 
$$  \eqalign{ 
S \sim S'  \quad \Longleftrightarrow \quad 
\big(\exists h \in C^1_s(G,N)\big)\ 
 S' = (C_N \circ h) \cdot S. \cr} $$
In Lemma I.10 below we shall see that for each (strongly) smooth outer action 
$S$ the outer action $S_h := (C_N \circ h)\cdot S$ also is a (strongly) smooth outer action 
because $\delta_S = C_N \circ \omega$ implies 
$$ \delta_{S_h}(g,g') 
= C_N\Big( h(g)S(g)(h(g')) \omega(g,g') h(gg')^{-1}\Big) \quad \hbox{ for } \quad g,g' \in G. $$

We write $[S]$ for the equivalence class of $S$ and   
 call $[S]$ a {\it smooth $G$-kernel}. 
\qed

\Remark I.5. (a) A smooth outer action $S \: G \to \Aut(N)$ need not be a group 
homomorphism, but $\delta_S = C_N \circ \omega$ 
implies that $S$ induces a group homomorphism $s := Q_N \circ S \: G \to \Out(N)$, 
where 
$Q_N \: \Aut(N) \to \Out(N)$ denotes the quotient homomorphism. 
Clearly $s$ depends only on the class $[S]$ of $S$ because 
$Q_N \circ C_N = \1$. 

 (b) If $G$ and $N$ are discrete groups, then for each homomorphism 
$s \: G \to \Out(N)$ 
there exists a map $S \: G \to \Aut(N)$ with $S(\1) = \id_N$ and 
$Q_N \circ S = s$ and
a map $\omega \: G \times G \to N$ with 
$\delta_S = C_N \circ \omega.$
In this case all outer actions are smooth, 
the smoothness conditions on $S$ and $\omega$ are vacuous, and 
$S \sim S'$ is equivalent to $Q_N \circ S = Q_N \circ S'$, so that a smooth 
$G$-kernel is nothing but a homomorphism $s \: G \to \Out(N)$. 

 (c) If $S \in C^1_s(G,\Aut(N))$ 
is a homomorphism of groups, then we may choose $\omega = \1$ and 
$S$ defines a smooth action of $G$ on $N$. In fact, the action is smooth on a 
set of the form $U \times N$, so that the assertion follows from the fact that 
all automorphisms $S(g)$ are smooth. 

 (d) Suppose that $S_1, S_2 \in C^1_s(G,\Aut(N))$ and that the map 
$$ U \times N \to N^2,\quad (g,n) \mapsto (S_1(g)(n), S_2(g)(n)) $$ 
is smooth. Then the map 
$U \times N \to N,(g,n) \mapsto S_1(g)S_2(g)(n)$ 
is smooth because it is a composition of smooth maps. Therefore 
$C^1_s(G,\Aut(N))$ is closed under pointwise products. 

 (e) If $S \: G \to \Aut(N)$ is a smooth outer action, then we obtain 
in particular 
$$ C_N(\omega(g,g^{-1})) = \delta_S(g,g^{-1}) = S(g) S(g^{-1}), $$
which leads to 
$$ S(g)^{-1} = S(g^{-1}) C_N(\omega(g,g^{-1})^{-1}). $$
Hence there also exists an identity neighborhood $U \subeq G$ for which the map 
$$ U \times N \to N, \quad (g,n) \mapsto S(g)^{-1}(n) = 
 S(g^{-1}) \omega(g,g^{-1})^{-1} n \omega(g,g^{-1}) $$
is smooth. 
\qed

\Lemma I.6. Each smooth outer action $S \: G \to \Aut(N)$ defines on 
$Z(N)$ the structure of a $G$-module by $g.z := S_Z(g)(z) := S(g)(z)$. 
If, in addition, $Z(N)$ is an initial Lie subgroup of $N$, then $Z(N)$ is a smooth $G$-module. 

\Proof. First we observe that inner automorphisms act trivially on the center 
$Z(N)$ of $N$, which implies that $g.z := S(g)(z)$ defines indeed a $G$-action on $Z(N)$. 

Let $U \subeq G$ be an open identity neighborhood such that the map 
$U \times N \to N, (g,n) \mapsto S(g).n$
is smooth. If $Z(N)$ is an initial Lie subgroup of $N$, restriction and 
co-restriction define a smooth map 
$$ U \times Z(N) \to Z(N), \quad (g,z) \mapsto g.z = S_Z(g)(z). \leqno(1.2) $$
As each element of $g$ acts as a Lie group 
automorphism on $Z(N)$, 
the smoothness of the action $G \times Z(N) \to Z(N)$ follows from the smoothness of the map in (1.2). 
\qed

The following lemma sheds some more light on the concept of an initial subgroup of a 
Lie group. In the present paper we shall mostly use it for the center of a Lie group 
$N$ which is not always a Lie group, but mostly carries an initial Lie subgroup structure.

\Lemma I.7. Let $G$ be a Lie group and $H \leq G$ a subgroup. Then the following assertions hold: 

\litem{(a)} If $H$ is an initial normal  subgroup, 
then the conjugation action of $G$ on $H$ is smooth. 

\litem{(b)} If $H$ is totally disconnected, 
then the discrete group $H_d$ (the group $H$ endowed with the discrete topology) defines 
on $H$ via the inclusion map $\iota \: H_d \to H \subeq G$ the structure of an initial 
Lie subgroup. The same conclusion holds if all smooth arcs in $H$ are constant, i.e., if 
the smooth arc-components in $H$ are trivial. 

\litem{(c)} If $\dim G < \infty$, then any subgroup $H \leq G$ carries an initial Lie group 
structure. 

\litem{(d)} Let $G$ be a Lie group with a smooth exponential function 
and $$\z(\g)^G := \{ x \in \z(\g) \: (\forall g \in G)\ \Ad(g).x = x\}.$$ 
We assume that the 
group $\Gamma_Z := \exp^{-1}(\1) \cap \z(\g)$ is discrete and that $\z(\g)^G$ is Mackey 
complete. Then $Z(G)$ carries an initial Lie group structure $Z(G)^L$ with 
$\L(Z(G)^L) = \z(\g)^G.$

\Proof. (a) The conjugation map $G \times H \to G, (g,h) \mapsto ghg^{-1}$ is smooth
with values in $H$, hence smooth as a map $G \times H \to H$. 

(b) Since 
every discrete space carries a natural $0$-dimensional smooth manifold structure, 
the group $H_d$ is a Lie group with $\L(H_d) = \{0\}$. 
If $M$ is a smooth manifold and $f \: M \to G$ a smooth map with 
$f(M) \subeq H$, then $f$ maps each connected component of $M$ into a connected component of 
$H$, so that $f$ is locally constant because $H$ is totally disconnected. 
Hence the map $\iota^{-1} \circ f \: M \to H_d$ is locally constant and therefore smooth. 

(c) Let $H_a \subeq H$ be the arc-component of $H$, viewed as a topological subgroup of $G$. 
According to Yamabe's Theorem ([Go69]), the arc-component $H_a$ of $G$ is analytic, 
i.e., there exists a Lie subalgebra $\h \subeq \g$ such that 
$H_a = \la \exp_G \h \ra$. Let $H_a^L$ denote the group $H_a$ endowed with its intrinsic 
Lie group topology for which $\exp \: \h \to H_a^L$ is a local diffeomorphism in~$0$. 
Then $H \subeq \{ g \in G \: \Ad(g)\h = \h\}$ implies that $H$ acts by conjugation 
smoothly on $H_a^L$, so that $H$ carries a Lie group structure for which $H_a^L$ is an 
open subgroup. Let $H^L$ denote this Lie group. Now the 
inclusion map $\iota \: H^L \to G$ is an immersion whose differential in $\1$ is the 
inclusion $\h \into \g$. 

We claim that $\iota \: H^L \to G$ is an initial Lie subgroup. 
In fact, let $f \: M \to G$ be a smooth map from the smooth manifold $M$ to $G$ 
with $f(M) \subeq H$. We have to show that $f$ is smooth, i.e., that $f$ is smooth in a  
neighborhood of each point $m \in M$. Replacing $f$ by $f \cdot f(m)^{-1}$ and observing 
that the group operations in $H$ and $G$ are smooth, we may w.l.o.g.\ assume that 
$f(m) = \1$. 

Let $U \subeq \h$ be an open $\0$-neighborhood, $\m \subeq \g$ be a vector space complement 
to $\h$ and and $V \subeq \m$ an open $\0$-neighborhood for which the map 
$$ \Phi \: U \times V \to G, \quad (x,y) \mapsto \exp x \exp y $$
is a diffeomorphism onto an open subset of $G$. Then 
$$ H_a \cap (\exp U \exp V) = \bigcup_{\exp y \in H_a} \exp U \exp y $$
and each set $\exp U \exp y$ contained in $H_a$ also is an open subset of $H^L$. 
Since the topology of $H^L_a$ is second countable, the set 
$\exp^{-1}(H_a) \cap \m$
is countable. Every smooth arc 
$\gamma \: I \to \exp U \exp V$ is of the form 
$\gamma(t) = \exp \alpha(t) \exp \beta(t)$ with smooth arcs 
$\alpha \: I \to U$ and $\beta \: I \to V$, and   
for every smooth arc contained in $H_a$ the arc $\beta$ is constant. 

We conclude that if $W \subeq M$ is an open connected neighborhood of $m$ with 
$f(W) \subeq \exp U \exp V$, then $f(W) \subeq \exp U$. Then the map 
$$ \exp\res_U^{-1} \circ f\res_W \: W \to \h $$
is smooth, so that the corresponding map 
$$ \iota^{-1} \circ f \res_W = \exp_{H^L} \circ \exp\res_U^{-1} \circ f\res_W \: W \to H^L $$
is also smooth. This proves that the map $\iota^{-1} \circ f \: M \to H^L$ is smooth, 
and hence that $\iota \: H^L \to H$ is an initial Lie subgroup of $G$. 

(d) First we observe that the exponential function 
$\exp_Z \: \z(\g)^G \to Z(G)$
is a group homomorphism, and our condition implies that 
$\Gamma := \ker \exp_Z$ is a discrete subgroup of $\z(\g)^G$, so that 
$Z := \z(\g)^G/\Gamma$ carries a natural Lie group structure, and the smooth function 
$\exp_Z \: \z(\g)^G \to G$ 
factors through an injective immersion of Lie groups 
$\iota \: Z \to G$
whose image is the subgroup $\exp(\z(\g)^G)$ of $Z(G)$.
We now endow $Z(G)$ with the Lie group structure for which 
$\exp \z(\g)^G$ is an open subgroup isomorphic to $Z$ (Theorem~B.1) and write 
$Z(G)^L$ for this Lie group. Then the inclusion map 
$\iota \: Z(G)^L \to G$
is an injective immersion of Lie groups. 

To see that this defines an initial Lie group structure on $Z(G)$, let 
$f \: M \to G$ be a smooth map with $f(M) \subeq Z(G)$. 
We may w.l.o.g.\ assume that $M$ is connected and fix some $m \in M$. 
Then the left logarithmic derivative 
$$ \delta(f) \in \Omega^1(M,\z(\g)^G) \quad \hbox{ with } \quad 
\delta(f)(m) = f(m)^{-1} df(m) $$
is a smooth $\z(\g)^G$-valued $1$-form on $M$. In fact, for each smooth curve 
$\gamma \: I \to Z(G)$ we have $\gamma(t)^{-1}\gamma'(t) \in \z(\g)^G$ for each 
$t \in I$ because for each $g \in G$ the relation 
$c_g \circ \gamma = \gamma$ leads to 
$$ \gamma(t)^{-1}.\gamma'(t) 
= (c_g \circ \gamma)(t)^{-1}.(c_g \circ \gamma)'(t) 
= \Ad(g).(\gamma(t)^{-1}.\gamma'(t)). $$
The $1$-form $\delta(f)$ satisfies the Maurer--Cartan equation 
$$ d(\delta(f)) + {1\over 2} [\delta(f), \delta(f)] = 0 $$
([KM97]), so that  the fact that $\z(\g)$ is abelian implies that $\delta(f)$ is closed. 
Let $q_M \: \tilde M \to M$ be a universal covering map. 
As $\z(\g)$ is Mackey complete, the Poincar\'e Lemma holds for 
$\z(\g)$-valued $1$-forms, so that there exists a smooth function 
$f_1 \: \tilde M \to \z(\g)$
with $d f_1 = q_M^* \delta(f)$ (cf.\ [Ne02, Prop.~III.4]). 
If $\tilde m \in \tilde M$ is a point with 
$q_M(\tilde m) = m$, we may assume, in addition, that $f_1(\tilde m) = 0$. 
Then the smooth function 
$$ f_2 := \exp_Z \circ f_1 \: \tilde M  \to Z(G)^L  
\quad \hbox{ satisfies} \quad 
\delta(f_2) = d f_1 = q_M^* \delta(f)= \delta(f \circ q_M), $$
which in turn leads to 
$$ f \circ q_M = f(m) \cdot f_2 $$
since the solutions of the Maurer-Cartan equation are uniquely determined 
by their values in $\tilde m$. This implies that the function 
$f \circ q_M \: \tilde M \to Z(G)^L$ is smooth, and hence that the function 
$\iota^{-1} \circ f \: M \to Z(G)^L$ is smooth. 
\qed

The following example shows that in general $Z(G)$, considered as a topological group, 
is not a Lie group because it may be a totally disconnected non-discrete subgroup.  
Since all known Lie groups modeled on complete locally convex spaces 
satisfy the requirements from Lemma~I.7(d), for all these Lie groups $G$ the center 
$Z(G)$ carries an initial Lie group structure. Therefore the assumption that $Z(G)$ 
is an initial Lie subgroup is not restrictive. 

\Example I.8. We consider the Lie group 
$$ N := \C^\N \rtimes \R^\N \quad \hbox{ with } \quad 
((z_n), (t_n)) * ((z_n'), (t_n')) := ((z_n + e^{t_n}z_n'), (t_n + t_n')), $$
where the locally convex structure on $\C^\N$ and $\R^\N$ is the product topology. 
We observe that $N$ is an infinite topological product of groups isomorphic to 
$\C \rtimes \R$, endowed with the multiplication 
$$ (z,t) * (z',t') = (z + e^t z', t + t'). $$
Therefore the center of $N$ is 
$$ Z(N) = \{0\} \times (2\pi \Z)^\N \cong \Z^\N, $$
and this group is totally disconnected, but not discrete, hence not a Lie group. Therefore the center 
of a Lie group need not be a Lie group. 

Since $Z(N)$ is totally disconnected, Lemma~I.7 implies that 
the inclusion map $Z(N)_d \to N$ of the discrete group $Z(N)_d$ defines on  
$Z(N)$ the structure of an initial Lie subgroup.  
\qed

\Definition I.9. If $S$ is a smooth outer action of $G$ on the Lie group $N$, then we define 
$$ C^p_s(G,Z(N)) := \{ f \in C^p_s(G,N) \: \im(f) \subeq Z(N)\}, \quad 
C^2_{ss}(G,Z(N)) := C^2_{ss}(G,N) \cap C^2_s(G,Z(N)) $$ 
and accordingly 
$$Z^p_s(G,Z(N))_{S} := \{ f \in C^p_s(G, Z(N)) \: d_S f = \1\} $$
and  $Z^2_{ss}(G,Z(N))_{S} := Z^2_s(G,Z(N))_{S} \cap C^2_{ss}(G,Z(N))$ 
(cf.\ Lemma~I.6). 
Note that on $f \in C^p_s(G,Z(N))$ the group differential $d_S$ for the $G$-module structure 
can be written, in additive notation, as 
$$ \eqalign{ 
&\ \ \ \ (d_S f)(g_0, \ldots, g_p) \cr 
&= S(g_0)(f(g_1,\ldots, g_n)) 
 +  \sum_{j=1}^n (-1)^j f(g_0, \ldots, g_{j-1} g_j,\ldots, g_n) 
+ (-1)^{n+1} f(g_0, \ldots, g_{n-1}). \cr} $$
This implies that $d_S(C^p_s(G,Z(N))) \subeq C^{p+1}(G,Z(N))$, so that we obtain a differential
complex $(C^\bullet_s(G,Z(N)),d_S)$ whose cohomology groups are denoted by 
$$ H^p_s(G,Z(N))_S := Z^p_s(G,Z(N))_S/d_S(C^p_s(G,Z(N))). $$

For $\alpha \in C^1_s(G,Z(N))$ and $g \in G$ we then have 
$$ \eqalign{ 
&\ \ \ \ (d_S\alpha)_g(x) \cr
&= (d_S \alpha)(g,x) - (d_S \alpha)(gxg^{-1},g) 
= g.\alpha(x) - \alpha(gx) +  \alpha(g) - (gxg^{-1}).\alpha(g) 
+ \alpha(gx) - \alpha(gxg^{-1}) \cr
&= g.\alpha(x) + \alpha(g) - (gxg^{-1}).\alpha(g) - \alpha(gxg^{-1}),\cr} $$ 
which is smooth in an identity neighborhood of $G$. Hence 
$$B^2_s(G,Z(N))_{S} := d_S(C^1_s(G,Z(N))) \subeq Z^2_{ss}(G,Z(N))_{S}, $$ 
and we define 
$$ H^2_{ss}(G,Z(N))_{S} := Z^2_{ss}(G,Z(N))_{S}/B^2_s(G,Z(N))_{S}, $$
and 
$$ H^3_{ss}(G,Z(N))_{S} := Z^3_{s}(G,Z(N))_{S}/d_S(C^2_{ss}(G,Z(N))). $$

We observe that the action of $G$ on $Z(N)$ defined by a smooth outer action 
$S$ only depends on the equivalence class $[S]$. In this sense we also write 
$$ Z^p_s(G,Z(N))_{[S]} := Z^p(G,Z(N))_{S}, \quad 
H^p_s(G,Z(N))_{[S]} := H^p_s(G,Z(N))_{S} \quad \hbox{ etc. } \quad 
\qeddis 

If $Z(N)$ is an initial Lie subgroup of $N$, then it carries the structure of a 
smooth $G$-module 
(Lemma~I.6) and the definitions above are consistent with 
Definition~I.2, but in general we do not 
want to assume that $Z(N)$ is an initial Lie subgroup 
and it is not necessary for the results in Section~II below. 

In Definition~I.3 above  
we have defined $H^1_s(G,N)$ as the set of orbits of $N = C^0_s(G,N)$ on the 
set $Z^1_s(G,N)$. The following lemma describes an action of $C^1_s(G,N)$ that in some 
sense is a generalization of the $N$-action on $Z^1_s(G,N)$ for degree $1$ to degree $2$. 

\Lemma I.10. Let $G$ and $N$ be Lie groups and consider the group 
$C^1_s(G, N)$ with respect to pointwise multiplication. 
This group acts on the set 
$C^1_s(G, \Aut(N))$ by 
$h.S  := (C_N \circ h) \cdot S$
and on the product set 
$$ C^1_s(G, \Aut(N)) \times C^2_s(G,N) \quad \hbox{ by } \quad 
h.(S,\omega) := (h.S, h *_S \omega) \leqno(1.3) $$
with 
$$ (h *_S \omega)(g,g') := h(g) S(g)(h(g')) \omega(g,g') h(gg')^{-1}. $$
The stabilizer of $(S,\omega)$ is given by 
$$ C^1_s(G,N)_{(S,\omega)} = Z^1_s(G,Z(N))_S $$
which depends only on $[S]$, but not on $\omega$. 
Moreover, the following assertions hold: 
\litem{(1)} The subset $\{(S,\omega) \: \delta_S = C_N \circ \omega\}$ 
is invariant. 
\litem{(2)} The subset $\{(S,\omega) \in C^1_s(G, \Aut(N)) \times C^2_{ss}(G,N) \: 
\delta_S = C_N \circ \omega\}$ is also invariant. 
\litem{(3)} If $\delta_S = C_N \circ \omega$, then $\im(d_S \omega) \subeq Z(N)$ holds for 
$$ (d_S\omega)(g,g',g'') := S(g)\big(\omega(g',g'')\big) \omega(g, g'g'')
\omega(gg',g'')^{-1} \omega(g,g')^{-1}. $$
\litem{(4)} If $\delta_S = C_N \circ \omega$ and 
$h.(S,\omega) = (S',\omega')$, then $d_{S'} \omega' = d_S \omega.$
\litem{(5)} If $\delta_S = C_N \circ \omega$, then 
$d_S \omega \in Z^3_s(G,Z(N))_S$. 
\litem{(6)}  For a smooth outer action $S$ 
the cohomology class 
$[d_S \omega] \in H^3_{s}(G,Z(N))_S$
depends only on the equivalence class $[S]$. 
Similarly, for a strongly smooth outer action and 
$\omega \in C^2_{ss}(G,N)$ with $\delta_S = C_N \circ \omega$,  
the cohomology class 
$[d_S \omega] \in H^3_{ss}(G,Z(N))_S$ 
depends only on the equivalence class $[S]$. 

\Proof. That (1.3) defines a group action follows from 
the trivial relation $(hh').S = h.(h'.S)$,  and 
$$ \eqalign{ 
&\ \ \ \ \big(h *_{h'.S} (h' *_S \omega)\big)(g,g')  
= h(g) (h'.S)(g)(h(g')) (h' *_S \omega)(g,g') h(gg')^{-1} \cr
&= h(g) h'(g) S(g)(h(g')) h'(g)^{-1} h'(g) S(g)(h'(g')) \omega(g,g') h'(gg')^{-1} 
h(gg')^{-1} \cr
&= (hh')(g) S(g)(hh'(g')) \omega(g,g') (hh')(gg')^{-1} 
= ((hh') *_S \omega)(g,g'). \cr} $$

To calculate the stabilizer of the pair $(S,\omega)$, we observe that the 
condition $h.S = S$ is equivalent to $h \in C^1_s(G,Z(N))$. Then 
$h *_S \omega = \omega \cdot d_S h$, and this equals $\omega$ if and only if 
$h$ is a cocycle. Therefore $C^1_s(G,N)_{(S,\omega)} = Z^1_s(G,Z(N))_S$. 

 (1) If the pair $(S,\omega)$ satisfies $\delta_S = C_N \circ \omega$, then we obtain 
for $h \in C^1_s(G,N)$ with Formula (1.1) (after Definition~I.1), the relation 
$$ \eqalign{ \delta_{h.S}(g,g') 
&= \delta_{(C_N \circ h) \cdot S}(g,g')  
= C_N(h(g)) c_{S(g)}(C_N(h(g'))) \delta_S(g,g') C_N(h(gg')^{-1})\cr
&= C_N(h(g))C_N(S(g)(h(g')))C_N(\omega(g,g')) C_N(h(gg')^{-1})\cr
&= C_N\big(h(g) S(g)(h(g')) \omega(g,g') h(gg')^{-1}\big) 
= (C_N \circ (h *_S \omega))(g,g'). \cr} $$

 (2) Suppose that $\omega \in C^2_{ss}(G,N)$ satisfies 
$\delta_S = C_N \circ \omega$. We have to show that 
$h *_S \omega \in C^2_{ss}(G,N)$. In fact, for $g,x \in G$ we have 
$$ \eqalign{ (h *_S \omega)_g(x) 
&= (h*_S \omega)(g,x) (h*_S \omega)(gxg^{-1},g)^{-1} \cr
&= h(g) S(g)(h(x)) \omega(g,x) h(gx)^{-1} 
\big(h(gxg^{-1}) S(gxg^{-1})(h(g)) \omega(gxg^{-1},g) h(gx)^{-1}\big)^{-1} \cr 
&= h(g) S(g)(h(x)) \omega(g,x) \omega(gxg^{-1},g)^{-1} S(gxg^{-1})(h(g)^{-1})
h(gxg^{-1})^{-1} \cr 
&= h(g) S(g)(h(x)) \omega_g(x) S(gxg^{-1})(h(g)^{-1}) h(gxg^{-1})^{-1}. \cr } $$
This expression is a product of factors which are smooth functions of $x$ in an identity 
neighborhood of $G$. Therefore $h *_S \omega \in C^2_{ss}(G,N)$. 

 (3) This follows from 
$$ \eqalign{ 
&\ \ \ \ C_N((d_S\omega)(g,g',g'')) \cr
&= C_N\Big(S(g)(\omega(g',g'')) \omega(g,g'g'')\omega(gg',g'')^{-1} 
\omega(g,g')^{-1}\Big) \cr 
&= S(g) (C_N \circ \omega)(g',g'') S(g)^{-1} 
(C_N \circ \omega)(g,g'g'') (C_N \circ \omega)(gg',g'')^{-1} 
(C_N \circ \omega)(g,g')^{-1} \cr 
&= S(g) \delta_S(g',g'') S(g)^{-1} \delta_S(g,g'g'') \delta_S(gg',g'')^{-1} 
\delta_S(g,g')^{-1} \cr 
&=  S(g)\big(S(g')S(g'')S(g'g'')^{-1}\big)S(g)^{-1} 
S(g)S(g'g'')S(gg'g'')^{-1}S(gg'g'')\cr
&\ \ \ \
S(g'')^{-1}S(gg')^{-1}S(gg')S(g')^{-1}S(g)^{-1}
= \1.\cr} $$

 (4) In the following calculation we frequently use that 
$d_S \omega$ and $d_{S'} \omega'$ have values in $Z(N)$ to see that 
they commute with all elements of $N$. We use under-braces to indicate which parts 
of the formula cancel each other. 

We calculate 
$$ \eqalign{ 
&\ \ \ \ (d_{S'}\omega')(g,g',g'') \cr
&= S'(g)(\omega'(g',g''))
\omega'(g,g'g'') \omega'(gg',g'')^{-1}
\omega'(g,g')^{-1} \cr
&= \big(C_N(h(g))\circ S(g)(\omega'(g',g''))\big)
\omega'(g,g'g'') \omega'(gg',g'')^{-1}
\omega'(g,g')^{-1} \cr
&= \ubrace{h(g)}_1 S(g)[h(g')S(g')(h(g''))
\omega(g',g'') \ubrace{h(g'g'')}_5] \ubrace{h(g)^{-1}}_2 \cr
&\ \ \ \  \ubrace{h(g)}_2 \ubrace{S(g)(h(g'g''))}_5 
\omega(g,g'g'') \ubrace{h(gg'g'')^{-1}}_3 
\ubrace{h(gg'g'')}_3 \omega(gg',g'')^{-1} 
S(gg')(h(g''))^{-1} \ubrace{h(gg')^{-1}}_4 \cr 
&\ \ \ \ \ubrace{h(gg')}_4 \omega(g,g')^{-1} 
S(g)(h(g'))^{-1} \ubrace{h(g)^{-1}}_1 \cr 
&= \ubrace{S(g)(h(g'))}_1 S(g)\circ S(g')(h(g''))
\big[S(g)(\omega(g',g'')) \omega(g,g'g'') \omega(gg',g'')^{-1}\big] \cr
&\ \ \ \ S(gg')(h(g''))^{-1} \omega(g,g')^{-1} 
\ubrace{S(g)(h(g'))^{-1}}_1 \cr 
&= [\omega(g,g')S(gg')(h(g'')) \omega(g,g')^{-1}] 
\big[ (d_S\omega)(g,g',g'')\omega(g,g')\big] S(gg')(h(g'')^{-1})
\omega(g,g')^{-1} \cr
&= (d_S \omega)(g,g',g'') \cr} $$
because the values of $d_S \omega$ are central.

 (5) It is easily seen that $d_S \omega$ vanishes if one of its three
arguments is~$\1$, so that, in view of (3), we have $d_S \omega \in C^3_s(G,Z(N))$. 
It remains to show that it is a cocycle. 

In the following calculation we write the group structure on
$Z(N)$ multiplicatively. We shall use several times that the values of
$d_S\omega$ are central, so that they commute with all values of~$\omega$. 
We simply write $g.z := S(g)(z)$ for the action of $G$ on $Z(N)$ induced by the smooth outer 
action~$S$. We have to show that for 
$g,g',g'',g''' \in G$ the following expression vanishes: 
$$ \eqalign{ &\ \ \ \ d_{G}(d_S\omega)(g,g',g'',g''') \cr
&=  g.(d_S\omega)(g',g'',g''') (d_S\omega)(gg',g'',g''')^{-1} 
(d_S\omega)(g,g'g'',g''')(d_S\omega)(g,g', g''g''')^{-1} (d_S\omega)(g,g',g'') \cr
&=  g.(d_S\omega)(g',g'',g''') (d_S\omega)(gg',g'',g''')^{-1} 
(d_S\omega)(g,g'g'',g''')(d_S\omega)(g,g',g'') (d_S\omega)(g,g', g''g''')^{-1} \cr
&=  S(g)
[S(g')(\omega(g'',g''')) \omega(g',g''g''') \omega(g'g'',g''')^{-1}
\omega(g',g'')^{-1}] \cr
&\ \ \ \ \omega(gg',g'') \omega(gg'g'',g''') \omega(gg', g''g''')^{-1}
S(gg')(\omega(g'',g''')^{-1}) \cr
&\ \ \ \ S(g)(\omega(g'g'',g''')) \omega(g,g'g''g''')
\omega(gg'g'',g''')^{-1} \omega(g,g'g'')^{-1} \cr
&\ \ \ \ S(g)(\omega(g',g'')) \omega(g,g'g'') \omega(gg',g'')^{-1}
\omega(g,g')^{-1} \cr
&\ \ \ \ \omega(g,g') \omega(gg',g''g''')\omega(g,g'g''g''')^{-1} 
S(g)(\omega(g',g''g''')^{-1})\cr
&=  \big[\omega(g,g') S(gg')(\omega(g'',g''')) \omega(g,g')^{-1}\big]  
S(g)[\omega(g',g''g''') \omega(g'g'',g''')^{-1}
\omega(g',g'')^{-1}] \cr
&\ \ \ \ \big[\omega(gg',g'') \omega(gg'g'',g''') \omega(gg',
g''g''')^{-1}
S(gg')(\omega(g'',g''')^{-1}) \big] \cr
&\ \ \ \ \big[S(g)(\omega(g'g'',g''')) \omega(g,g'g''g''')
\omega(gg'g'',g''')^{-1} \omega(g,g'g'')^{-1}\big] \cr
&\ \ \ \ S(g)(\omega(g',g'')) \omega(g,g'g'') \omega(gg',g'')^{-1}
\omega(gg',g''g''')\omega(g,g'g''g''')^{-1} 
S(g)(\omega(g',g''g''')^{-1})\cr
&= \omega(g,g') \big[\omega(gg',g'') \omega(gg'g'',g''') \omega(gg',
g''g''')^{-1}S(gg')(\omega(g'',g''')^{-1}) \big] \cr
&\ \ \ \ S(gg')(\omega(g'',g''')) \omega(g,g')^{-1} 
S(g)[\omega(g',g''g''') \omega(g'g'',g''')^{-1}
\omega(g',g'')^{-1}] \cr
&\ \ \ \ S(g)(\omega(g',g'')) \omega(g,g'g'') \omega(gg',g'')^{-1}
\omega(gg',g''g''')\omega(g,g'g''g''')^{-1} \cr
&\ \ \ \ S(g)(\omega(g',g''g''')^{-1}) \big[S(g)(\omega(g'g'',g''')) \omega(g,g'g''g''')
\omega(gg'g'',g''')^{-1} \omega(g,g'g'')^{-1}\big] \cr
}$$ $$\eqalign{ 
&= \omega(g,g') \omega(gg',g'') \omega(gg'g'',g''') \omega(gg',g''g''')^{-1}
\omega(g,g')^{-1} 
S(g)[\omega(g',g''g''') \omega(g'g'',g''')^{-1}]\cr
&\ \ \ \ \omega(g,g'g'') \omega(gg',g'')^{-1}
\omega(gg',g''g''')\omega(g,g'g''g''')^{-1} 
S(g)(\omega(g',g''g''')^{-1})  \cr
&\ \ \ \ \big[S(g)(\omega(g'g'',g''')) \omega(g,g'g''g''')
\omega(gg'g'',g''')^{-1} \omega(g,g'g'')^{-1}\big] \cr
&= \omega(g,g') \omega(gg',g'') \omega(gg'g'',g''') \omega(gg',g''g''')^{-1}
\omega(g,g')^{-1} 
S(g)[\omega(g',g''g''') \omega(g'g'',g''')^{-1}]\cr
&\ \ \ \ \big[S(g)(\omega(g'g'',g''')) \omega(g,g'g''g''')
\omega(gg'g'',g''')^{-1} \omega(g,g'g'')^{-1}\big]\cr
&\ \ \ \ \omega(g,g'g'') \omega(gg',g'')^{-1}
\omega(gg',g''g''')\omega(g,g'g''g''')^{-1} 
S(g)(\omega(g',g''g''')^{-1})  \cr
&=  \omega(g,g')\omega(gg',g'') \omega(gg'g'',g''') \omega(gg',g''g''')^{-1}\omega(g,g')^{-1} \cr
&\ \ \ \ \big\{S(g)(\omega(g',g''g'''))\omega(g,g'g''g''')
\big[\omega(gg'g'',g''')^{-1} \omega(gg',g'')^{-1}
\omega(gg',g''g''')\big] \cr 
&\ \ \ \ \omega(g,g'g''g''')^{-1}S(g)(\omega(g',g''g''')^{-1}) \big\}\cr
&= \omega(g,g') \omega(gg',g'') \omega(gg'g'',g''') \omega(gg',g''g''')^{-1}\omega(g,g')^{-1} \cr
&\ \ \ \ c_{S(g)[\omega(g',g''g''')]\omega(g,g'g''g''')}
\big[\omega(gg'g'',g''')^{-1} \omega(gg',g'')^{-1}
\omega(gg',g''g''')\big]  \cr
&=  \omega(g,g')\omega(gg',g'') \omega(gg'g'',g''') \omega(gg',g''g''')^{-1}\omega(g,g')^{-1} \cr
&\ \ \ \ c_{\omega(g,g')\omega(gg',g''g''')}
\big[\omega(gg'g'',g''')^{-1} \omega(gg',g'')^{-1}
\omega(gg',g''g''')\big]  \cr
&= \omega(g,g')  \omega(gg',g'') \omega(gg'g'',g''') \omega(gg',g''g''')^{-1}\omega(g,g')^{-1} 
\omega(g,g')\omega(gg',g''g''') \cr
&\ \ \ \ \big[\omega(gg'g'',g''')^{-1} \omega(gg',g'')^{-1}
\omega(gg',g''g''')\big] \omega(gg',g''g''')^{-1} \omega(g,g')^{-1}\cr
&= \1. \cr} $$

 (6) We prove the assertion for strongly smooth outer actions. The other 
case is proved similarly by using elements $\omega \in C^2_s(G,N)$ instead. 

We fix some $\omega \in C^2_{ss}(G,N)$ with 
$\delta_S = C_N \circ \omega$. If $\omega' \in C^2_{ss}(G,N)$ also satisfies 
$\delta_S = C_N \circ \omega'$, 
then $\beta := \omega' \cdot \omega^{-1} \in C^2_{ss}(G,Z(N))$ and thus  
$d_S \omega' = d_S \omega + d_S \beta = d_S \omega + d_S \beta,$
hence $[d_S\omega'] = [d_S \omega]$ and the cohomology class 
$[d_S \omega]$ does not depend on the choice of $\omega$. 

If $S' \sim S$, then there exists an $h \in C^1_s(G,N)$ with 
$S' = (C_N \circ h) \cdot S$, and we have seen in (2) that  
$\omega' := h *_S \omega \in C^2_{ss}(G,N)$ satisfies 
$\delta_{S'} = C_N \circ \omega'$. Now $d_{S'} \omega' = d_S \omega$ implies in particular 
that $[d_{S'} \omega'] =[d_S \omega]$ and that this 
cohomology class neither depends on the representative 
$S'$ of $[S]$ nor on the choice of $\omega$. 
\qed

\Definition I.11. Let $S \in C^1_s(G,\Aut(N))$ be a strongly smooth outer action and 
pick $\omega \in C^2_{ss}(G,N)$ with 
$\delta_S = C_N \circ \omega$. 
In view of Lemma~I.10(6), the cohomology class 
$$ \chi_{ss}(S) := [d_S \omega] \in H^3_{ss}(G,Z(N))_{S} $$
does not depend on the choice of $\omega$ and is constant on the equivalence class of 
$S$, so that we may also write $\chi_{ss}([S]) := \chi_{ss}(S)$. We call 
$\chi_{ss}(S)$ the {\it characteristic class of $S$}. 

For a smooth outer action we likewise obtain a characteristic class 
$$ \chi_s(S) = \chi_s([S]) = [d_S\omega] \in H^3_s(G,Z(N))_{S}. 
\qeddis

\sectionheadline{II. Extensions of Lie groups} 

\nin In this section we discuss the basic concepts related to extensions of 
Lie groups. Since every discrete group can be viewed as a Lie group, our discussion 
includes in particular the algebraic context of group extensions 
([MacL63], [EML47]).

\Definition II.1. Let $G$ be a Lie group. A subgroup $H$ of $G$ 
is called a {\it split Lie subgroup} if it carries a Lie group structure for which the 
inclusion map $i_H \: H \into G$ is a morphism of Lie groups and the right action of 
$H$ on  $G$ defined by restricting the multiplication map of $G$ to a map 
$G \times H \to G$ defines a smooth $H$-principal bundle. This means that the coset space 
$G/H$ is a smooth manifold and that the quotient map $\pi \: G \to G/H$ has smooth local 
sections. 
\qed

\Examples II.2. Since the Lie algebra $\h$ of a Lie subgroup $H$ of a Lie group
$G$ need not have a closed complement in 
$\g$, not every Lie subgroup is split. A simple example is the
subgroup $H := c_0(\N,\R)$ in $G := \ell^\infty(\N,\R)$ 
(cf.\ [Wei95, Satz IV.6.5]). 
\qed

\Lemma II.3. Each split Lie subgroup of $G$ is initial. 

\Proof. The condition that $H$ is a split Lie subgroup implies that 
there exists an open subset $U$ of some locally convex space $V$ 
and a smooth map $\sigma \: U \to G$ such that the map 
$$ U \times H \to G, \quad (x,h) \mapsto \sigma(x) h $$
is a diffeomorphism onto an open subset of $G$. 
We conclude that for every manifold $M$ a map 
$f \: M \to H$ is smooth whenever it is smooth as a map 
to $G$, resp., the open subset $\sigma(U)H$ of~$G$. 
\qed 

\Definition II.4. An {\it extension of Lie groups} is a short exact sequence 
$$ \1 \to N \sssmapright{\iota} \hat G \sssmapright{q} G \to \1 $$
of Lie group morphisms, for which 
$N \cong \ker q$ is a split Lie subgroup. 
This means that $\hat G$ is a smooth $N$-principal 
bundle over $G$ and $G \cong \hat G/N$. In the following we shall identify 
$N$ with the subgroup $\iota(N) \trile \hat G$. 

We call two extensions 
$N \into \hat G_1 \onto G$ and 
$N \into \hat G_2 \onto G$ of the Lie group $G$ by the 
Lie group $N$ {\it equivalent} if there exists a Lie group morphism 
$\phi \: \hat G_1 \to \hat G_2$ such that the following diagram commutes: 
$$ \matrix{
 N & \into& \hat G_1 & \onto & G \cr 
\mapdown{\id_N} & & \mapdown{\phi} & & \mapdown{\id_G} \cr 
 N & \into& \hat G_2 & \onto & G. \cr } $$
It is easy to see that any such $\phi$ is in particular an isomorphism of
 Lie groups. We write $\Ext(G,N)$ for the set of equivalence classes of
Lie group extensions of $G$ by~$N$. 

We call an extension $q \: \hat G \to G$ with $\ker q = N$ 
{\it split} if there exists a Lie group morphism 
$\sigma \: G \to \hat G$ with $q \circ \sigma = \id_G$. 
In this case the map 
$N \rtimes_S G \to \hat G, (n,g) \mapsto n \sigma(g)$
is an isomorphism, where the semidirect product is defined by the 
homomorphism 
$$ S := C_N \circ \sigma \: G \to \Aut(N), $$
where $C_N \: \hat G \to \Aut(N)$ is the conjugation action of $\hat G$ on $N$. 
In view of Lemma~I.7, $S$ defines a smooth action of $G$ on $N$. 
\qed 

\Remark II.5. We give a description of Lie group extensions 
$N \into \hat G \onto G$ in terms of data associated to $G$ and $N$. 
Let $q\: \hat G \to G$ be a Lie group extension of $G$ by $N$.
By assumption, the map $q$ has a smooth local section. 
Hence there exists a locally smooth global section 
$\sigma \: G \to \hat G$ which is {\it normalized} in the sense that 
$\sigma(\1) = \1$, i.e., $\sigma \in C^1_s(G,\hat G)$. Then the map 
$$ \Phi \: N \times G \to \hat G, \quad (n,g) \mapsto n \sigma(g)$$ is
a bijection which restricts to a local diffeomorphism on an identity neighborhood. 
In general $\Phi$ is not continuous, but we may nevertheless use it to identify 
$\hat G$ with the product set $N \times G$ endowed with the multiplication 
$$ (n,g) (n',g') = (n S(g)(n') \omega(g,g'), gg'), \leqno(2.1) $$
where 
$$ S := C_N \circ \sigma \: G \to \Aut(N) \quad \hbox{ for } \quad 
C_N \: \hat G \to \Aut(N),\ C_N(g) = gng^{-1}, \leqno(2.2) $$
and 
$$ \omega \: G \times G \to N, \quad 
(g,g') \mapsto \sigma(g) \sigma(g')\sigma(gg')^{-1}. \leqno(2.3) 
$$
Note that $\omega$ is smooth in an identity neighborhood and the map 
$G \times N \to N, (g,n) \mapsto S(g).n$ is smooth in a set of the form 
$U \times N$, where $U$ is an identity neighborhood of $G$. 
The maps $S$ and $\omega$ satisfy the relations 
$$  \sigma(g) \sigma(g') = \omega(g,g')\sigma(gg'), \leqno(2.4) $$
$$ S(g)S(g') = C_N(\omega(g,g')) S(gg'), \quad \hbox{ i.e.,} \quad 
\delta_S = C_N \circ \omega, \leqno(2.5) $$
and  
$$ \omega(g,g') \omega(gg',g'') 
= S(g)\big(\omega(g',g'')\big) \omega(g, g'g''). \leqno(2.6) $$
\qed

\Definition II.6. The elements of the set 
$$ Z^2_{ss}(G,N) := \{ (S,\omega) \in C^1_s(G,\Aut(N)) 
\times C^2_{ss}(G,N) \: 
\delta_S = C_N \circ \omega, d_S \omega = \1\} $$ 
are called {\it smooth factor systems for the pair $(G,N)$} or 
{\it locally smooth $2$-cocycles}. 
\qed

\Lemma II.7. Let $S \in C^1(G,\Aut(N))$ and $\omega \in C^2(G,N)$ with 
$\delta_S = C_N \circ \omega$. We define a product on $N \times G$ by 
$$ (n,g) (n',g') = (n S(g)(n') \omega(g,g'), gg'). $$
Then $N \times G$ is a group if and only if $d_S \omega = \1$. 
Inversion in this group is given by 
$$ (n,g)^{-1} = (S(g)^{-1}\big(n^{-1} \omega(g,g^{-1})^{-1}\big),
g^{-1}) =  (\omega(g^{-1},g)^{-1} S(g^{-1})(n^{-1}), g^{-1}) $$
and in particular 
$$ (\1,g)^{-1} = (S(g)^{-1}.\big(\omega(g,g^{-1})^{-1}\big),g^{-1}) 
=  (\omega(g^{-1},g)^{-1}, g^{-1}). $$
Conjugation is given by 
$$ \eqalign{ 
&\ \ \ \ (n,g)(n',g')(n,g)^{-1}
= \Big(n S(g)(n') \omega_g(g') S(gg'g^{-1})^{-1}(n^{-1}), gg'g^{-1}\Big), \cr} $$
so that 
$(n,\1)(n',g')(n,\1)^{-1} 
= (n n' S(g')(n^{-1}), g')$
and 
$$ (\1,g)(n',g')(\1,g)^{-1} 
= \big( S(g)(n') \omega_g(g'), gg'g^{-1}\big). \leqno(2.7) $$

\Proof. The associativity of the multiplication on $N \times G$ is equivalent to the
equality of 
$$ \eqalign{\big((n,g) (n',g')\big)(n'',g'')
&= (n S(g)(n') \omega(g,g'), gg') (n'',g'') \cr
&= (n S(g)(n') \omega(g,g') S(gg')(n'') \omega(gg',g''), gg'g'') \cr} $$
and 
$$ \eqalign{ (n,g) \big((n',g') (n'',g'')\big) 
&= (n,g)(n' S(g')(n'') \omega(g',g''), g'g'')  \cr
&= (n S(g)\big(n' S(g')(n'') \omega(g',g'')\big)\omega(g,g'g''), gg'g'')\cr
&= (n S(g)(n') (S(g)S(g')(n'')) S(g)(\omega(g',g''))\omega(g,g'g''), gg'g'') \cr} $$
for all $g,g',g'' \in G$ and $n,n',n'' \in N$. 
This means that 
$$  \omega(g,g') S(gg')(n'') \omega(gg',g'') = (S(g)S(g')(n''))
S(g)(\omega(g',g''))\omega(g,g'g''). $$
For $n'' = 1$ this leads to $d_S \omega = \1$. 
If, conversely, 
$d_S \omega = \1$, then the associativity condition is equivalent to 
$$  \omega(g,g') S(gg')(n'') \omega(gg',g'') 
= (S(g)S(g')(n''))\omega(g,g') \omega(gg',g'') $$ 
and hence to $\delta_S = C_N \circ \omega$. Therefore the conditions 
$\delta_S = C_N \circ \omega$ 
and $d_S \omega = \1$ are equivalent to the associativity of the multiplication 
on $N \times G$. 

To see that we actually obtain a group, we first observe that $S(\1) = \id_N$ 
implies that $\1 := (\1,\1)$ is an identity element of 
$\hat G := N \times G$, so that $(\hat G,\1)$ is a monoid. 
For $(n,g) \in \hat G$ the element 
$$ (S(g)^{-1}\big(n^{-1} \omega(g,g^{-1})^{-1}\big), g^{-1}) $$
is a right inverse and likewise 
$$ (\omega(g^{-1},g)^{-1} S(g^{-1})(n^{-1}), g^{-1}) $$
is a left inverse. 
Now the associativity of $\hat G$ implies that left and right inverse
are equal, hence an inverse of $(n,g)$. Therefore $\hat G$ is a group. 

The formula for the inversion has already been obtained above. 
For the conjugation we first oberve that for $g,g' \in G$ we have 
$$ \eqalign{ \1 
&= (d_S \omega)(gg',g^{-1},g)  
= S(gg')(\omega(g^{-1},g)) \omega(gg',\1) \omega(gg'g^{-1},g)^{-1} \omega(gg',g^{-1})^{-1}\cr
&= S(gg')(\omega(g^{-1},g)) \omega(gg'g^{-1},g)^{-1} \omega(gg',g^{-1})^{-1}, \cr}$$
so that 
$$ S(gg')(\omega(g^{-1},g))^{-1} \omega(gg',g^{-1}) 
= \omega(gg'g^{-1},g)^{-1}. $$
We now obtain 
$$ \eqalign{
&\ \ \ \ (n,g)(n',g')(n,g)^{-1} \cr
&= \Big(n S(g)(n') \omega(g,g'),gg'\Big) 
\Big(\omega(g^{-1},g)^{-1} S(g^{-1})(n^{-1}), g^{-1}\Big) \cr
&= \Big(n S(g)(n') \omega(g,g') S(gg')\big(\omega(g^{-1},g)^{-1} S(g^{-1})(n^{-1})\big) 
\omega(gg',g^{-1}) , gg'g^{-1}\Big) \cr
&= \Big(n S(g)(n') \omega(g,g') \omega(gg'g^{-1},g)^{-1}\omega(gg',g^{-1})^{-1} 
\big(S(gg') S(g^{-1})(n^{-1})\big) \omega(gg',g^{-1}) , gg'g^{-1}\Big) \cr
&= \Big(n S(g)(n') \omega_g(g')\delta_S(gg',g^{-1})^{-1} 
\big(S(gg') S(g^{-1})(n^{-1})\big), gg'g^{-1}\Big) \cr
&= \Big(n S(g)(n') \omega_g(g') S(gg'g^{-1})^{-1}(n^{-1}), gg'g^{-1}\Big). \cr} $$
\qed

For a smooth factor system $(S,\omega)$ we write 
$N \times_{(S,\omega)} G$ for the set $N \times G$ endowed with the group  
multiplication (2.1).   

\Proposition II.8. If $(S,\omega)$ is a smooth factor system, 
then $N \times_{(S,\omega)} G$ carries a unique structure of a Lie group for which the map 
$$ N \times G \to N \times_{(S,\omega)} G, \quad 
(n,g) \mapsto (n,g) $$
is smooth on a set of the form $N \times U$, where $U$ is an open identity neighborhood 
in $G$. 
Then 
$$ q \: N \times_{(S,\omega)} G \to G, \quad (n,g) \mapsto g $$
is a Lie group extension of $G$ by $N$. 
Each Lie group extension $q \: \hat G \to G$ of $G$ by 
$N$ gives rise to smooth factor systems by choosing a locally smooth
normalized section $\sigma \: G \to \hat G$ and defining $(S,\omega) 
:= (C_N \circ \sigma, \delta_\sigma)$. Then the extension 
$q \: \hat G \to G$ is equivalent to $N \times_{(S,\omega)} G$. 

\Proof. (1) Let $U_G \subeq G$ be an open symmetric $\1$-neighborhood such that the maps 
$$ U_G \times N \to N, \quad (g,n) \mapsto S(g).n 
\quad \hbox{ and } \quad 
\omega\res_{U_G \times U_G} $$
are smooth. 
We consider the subset 
$$ U := N \times U_G = q^{-1}(U_G) \subeq \hat G := N \times_{(S,\omega)} G. $$
Then $U = U^{-1}$ because $q$ is a group homomorphism. 
We endow $U$ with the product manifold structure
from $N \times U_G$. Since the multiplication 
$m_G\res_{U_G \times U_G} \: U_G \times U_G \to
G$ is continuous, there exists an open symmetric identity neighborhood $V_G
\subeq U_G$ with $V_G V_G \subeq U_G$. Then the set 
$V := N \times V_G$ is an open subset of $U$ such that the multiplication map 
$$ V \times V \to U, \quad \big((n,g),(n',g')\big) \to 
(n S(g)(n') \omega(g,g'), gg') $$
is smooth. The inversion 
$$ U \to U, \quad (n,g) \mapsto 
(\omega(g^{-1},g)^{-1} S(g^{-1})(n^{-1}), g^{-1}) $$
(cf.~Lemma II.7) is also smooth. 

For $(n,g) \in \hat G$ let $V_g \subeq U_G$ be an open identity neighborhood with 
$c_g(V_g) \subeq U_G$. Then $c_{(n,g)}(q^{-1}(V_g)) \subeq U$. 
That the conjugation map 
$c_{(n,g)} \: q^{-1}(V_g) \to U$ 
is smooth in an identity neighborhood follows immediately from 
formula (2.7) in Lemma~II.7 because 
$\omega \in C^2_{ss}(G,N)$ ensures that the function $\omega_g$ is smooth in an identity 
neighborhood of~$G$. 

Eventually Theorem B.1 implies that $\hat G$ 
carries a unique Lie group structure for
which the inclusion map $U = N \times U_G \into \hat G$ is a local
diffeomorphism onto an identity neighborhood. 
It is clear that with respect to the Lie group structure on $\hat G$, the
map $q \: \hat G \to G$ is a smooth $N$-principal bundle because the 
map $V_G \to \hat G, g \mapsto (0,g)$ defines a section of $q$ which
is smooth on an identity neighborhood in $G$ that might be smaller
than $V_G$. 

 (2) Assume, conversely, that $q \: \hat G \to G$ is a Lie group 
extension of $G$ by $N$. Then there exists an open $\1$-neighborhood $U_G \subeq G$ and a
smooth section $\sigma \: U_G \to \hat G$ of the map $q \: \hat G
\to G$. We extend $\sigma$ to a global section $G \to \hat G$ 
and put $(S,\omega):= (C_N \circ \sigma, \delta_\sigma)$. 
We have already seen in Definition~I.1 that 
$\omega = \delta_\sigma \in C^2_{ss}(G,N)$. 
Therefore $(S,\omega)$ is a smooth factor 
system, so that we can use (1) to obtain a Lie group structure on 
$N \times_{(S,\omega)} G$ for which 
there exists an 
identity neighborhood $V_G \subeq G$ for which the product map 
$$N \times V_G \to N \times_{(S,\omega)} G, \quad (n,v) \mapsto (n,\1)(0,v) = (n,v) $$
is smooth. This implies that the group isomorphism 
$N \times_{(S,\omega)} G \to \hat G$ is a local diffeomorphism, hence an
isomorphism of Lie groups.
\qed

If the group $G$ is connected, then we can weaken the condition 
on smooth factor systems by replacing $\omega \in C^2_{ss}(G,N)$ 
by the weaker condition $\omega \in C^2_s(G,N)$: 

\Lemma II.9. If $G$ is connected and 
$S \in C^1_s(G, \Aut(N))$ and $\omega \in C^2_s(G,N)$ satisfy 
$$ \delta_S = C_N  \circ \omega \quad \hbox{ and } \quad 
d_S \omega = \1, $$
then $\omega \in C^2_{ss}(G,N)$ and in particular 
$(S,\omega) \in Z^2_{ss}(G,N)$. 

\Proof. That the map 
$\omega_g \: G \to N$ 
is smooth in an identity neighborhood 
means that conjugation with $(\1,g)$ is smooth in an 
identity neighborhood because 
$$ \eqalign{ (\1,g)(n',g')(\1,g)^{-1} 
&= \Big( S(g)(n') \omega_g(g'), gg'g^{-1}\Big) \cr} $$
by formula (2.7). 

The set of all elements $x \in N \times_{(S,\omega)} G$ for which the 
conjugation $c_x$ is smooth in an identity neighborhood is a subsemigroup 
containing $N$ (cf.\ 
Lemma~II.7) and $\{\1\} \times W$ for some identity neighborhood $W$ in $G$.
If $G$ is connected, it is generated by $W$ as a semigroup, so that 
$N \times_{(S,\omega)} G$ is generated by $N \times W$. Therefore the local smoothness 
requirement for the functions $\omega_g$ is redundant. 
\qed

In the following proposition we show that equivalence classes of extensions correspond to 
orbits of the group $C^1_s(G,N)$ in the set of smooth factor systems. 

\Proposition II.10. For two smooth factor systems 
$(S, \omega), (S', \omega') \in Z^2_{ss}(G,N)$ the following are equivalent 
\litem{(1)} $N \times_{(S,\omega)} G$ and $N \times_{(S',\omega')} G$ are equivalent 
extensions of $G$ by $N$. 
\litem{(2)} There exists an $h \in C^1_s(G,N)$ with 
$h.(S,\omega) = (S',\omega')$. 

If these conditions are satisfied, then the map 
$$ \phi \: N \times_{(S',\omega')} G \to N \times_{(S,\omega)} G, \quad 
(n,g) \mapsto (n h(g), g) $$
is an equivalence of extensions and all equivalences of extensions 
$N \times_{(S',\omega')} G \to N \times_{(S,\omega)} G$ are of this form.   

\Proof. Let 
$\phi \: N \times_{(S',\omega')} G \to N \times_{(S,\omega)} G$
be an equivalence of extensions. Then there exists a function 
$h \in C^1_s(G,N)$ such that $\phi$ has the form 
$\phi(n,g) = (n h(g), g)$
with respect to the product coordinates on 
$N \times_{(S,\omega)} G$, resp., $N \times_{(S',\omega')} G$. 
We then have 
$$ \eqalign{ 
\phi(n,g) \phi(n',g') 
&= (n h(g), g) (n' h(g'), g') 
= (n h(g) S(g)(n' h(g')) \omega(g,g'), gg') \cr
&= \big(n [C_N(h(g))(S(g)(n'))] h(g) S(g)(h(g')) \omega(g,g'), gg'\big) \cr
&= \big(n (h.S)(g)(n') (h *_S \omega)(g,g') h(gg'), gg'\big) \cr 
&= \phi(n (h.S)(g)(n') (h *_S \omega)(g,g'), gg'). \cr} $$
Since $\phi$ is an equivalence of extensions, it is an injective morphism 
of groups, which implies that the product in $N \times_{(S',\omega')} G$ is given by 
$$ (n S'(g)(n') \omega'(g,g'),gg') = (n,g)(n',g') = (n (h.S)(g)(n') (h *_S \omega)(g,g'), 
gg'). $$
This implies that $\omega' = h *_S \omega$ and further that $S' = h.S$, so that 
``(1) $\Rarrow$ (2)'' is proved. 

If, conversely, $(S',\omega') = h.(S,\omega) = (h.S, h *_S \omega)$, 
then we define 
$$\phi \: N \times_{(S',\omega')} G \to N \times_{(S,\omega)} G, \quad 
\phi(n,g) := (n h(g),g), $$ 
and the calculation above shows that $\phi$ 
is a homomorphisms of groups. Since $\phi$ is smooth in an identity 
neighborhood, it is a morphisms of Lie groups. Clearly $\phi\res_N = \id_N$ 
and $\phi$ induces the identity map on $G$, so that it is an equivalence of 
extensions. This completes the proof. 
\qed

\Corollary II.11. The map 
$$ Z^2_{ss}(G,N) \to \Ext(G,N), \quad (S,\omega) \mapsto [N \times_{(S,\omega)} G] $$
induces a bijection 
$$ H^2_{ss}(G,N) := Z^2_{ss}(G,N)/C^1_s(G,N) \to \Ext(G,N). 
\qeddis 

The preceding proposition implies in particular that 
$N \times_{(S,\omega)} G \sim N \times_{(S',\omega')} G$ implies 
$[S] = [S']$, i.e., equivalent extensions correspond to the same 
smooth $G$-kernel. In the following we write $\Ext(G,N)_{[S]}$ for the set of 
equivalence classes of $N$-extensions of $G$ corresponding to the $G$-kernel $[S]$. 

\Theorem II.12. Let $S$ be a smooth outer action of 
$G$ on $N$ with $\Ext(N, G)_{[S]} \not=\eset$. 
Then each extension class in $\Ext(G, N)_{[S]}$ can be represented by a Lie group of the form 
$ N \times_{(S,\omega)} G$. Any other Lie group extension
$ N \times_{(S,\omega')} G$ representing an element of $\Ext(G,N)_{[S]}$ 
satisfies 
$$ \omega' \cdot \omega^{-1} \in Z^2_{ss}(G,Z(N))_{[S]}, $$
and the Lie groups  
$ N \times_{(S,\omega)} G$ and $N \times_{(S,\omega')} G$ define equivalent 
extensions if and only if 
$$ \omega' \cdot \omega^{-1} \in B^2_{s}(G,Z(N))_{[S]}. $$

\Proof. From Proposition II.8 we know that each Lie group 
extension $\hat G$ of $ G$ by $N$ is equivalent to one of the form
$ N \times_{(S',\omega')} G$. 
If $[S] = [S']$ and $h \in C^1_s(G,N)$ satisfies $h.S = S'$, then 
$h^{-1}.(S',\omega') = (S, h^{-1} *_{S'} \omega')$, so that 
$\omega'' := h^{-1} *_{S'} \omega'$ satisfies 
$N \times_{(S',\omega')} G \cong N \times_{(S,\omega'')} G$. This means that each extension class 
in $\Ext(G,N)_{[S]}$ can be represented by a group of the form 
$N \times_{(S,\omega'')} G$. 

If $(S,\omega)$ and $(S',\omega')$ are smooth factor systems, 
then 
$C_N \circ \omega = \delta_S = C_N \circ \omega'$
implies $\beta := \omega' \cdot \omega^{-1} \in C^2_{s}(G,Z(N))$. 
Further $\im(\beta) \subeq Z(N)$ leads to 
$\beta_g = \omega_g' \omega_g^{-1}$ for each $g \in G$, which shows that $\beta_g$ 
is smooth in an identity neighborhood, so that $\beta \in C^2_{ss}(G,Z(N))$. 
Moreover, $d_S \omega' = d_S \omega = \1$, so that 
$$ \1 = d_S \omega' = d_S \omega \cdot d_S \beta = d_S\beta $$
implies $\beta \in Z^2_{ss}(G,Z(N))_{S}$. This means that 
$\omega' \in \omega \cdot Z^2_{ss}(G,Z( N))_{S}$. 

In view of Proposition~II.10, the equivalence of the extensions 
$N \times_{(S,\omega')} G$ and $N \times_{(S,\omega)} G$ is equivalent to the existence of 
an $h \in C^1_s(G,N)$ with 
$$ S = h.S = (C_N \circ h) \cdot S \quad \hbox{ and } \quad 
\omega' = h *_S \omega. $$
Then $C_N \circ h = \id_N$ implies $h \in C^1_s(G,Z(N))$ which further leads to 
$\omega' = h *_S \omega = (d_S h) \cdot \omega$, i.e., 
$\omega' \omega^{-1} \in B^2_s(G,Z(N))_S$. If, conversely, 
$\omega' \omega^{-1} = d_S h$ for some $h \in C^1_s(G,Z(N))$, then 
$h.S = S$ and $h *_S \omega = \omega'$. 
\qed

\Corollary II.13. For 
a smooth $G$-kernel $[S]$ with $\Ext(G,N)_{[S]} \not=\eset$, the map 
$$ H^2_{ss}(G,Z(N))_{[S]} 
\times \Ext(G, N)_{[S]} \to  \Ext(G,N)_{[S]}, 
\quad (\beta, [N \times_{(S,\omega)} G]) \mapsto [N \times_{(S,\omega \cdot\beta)} G] $$
is a well-defined simply transitive action. 
\qed

\Remark II.14. (Abelian extensions) Suppose that $N = A$ is an abelian Lie group. 
Then the adjoint representation of $A$ is trivial and 
a smooth factor system $(S,\omega)$ for $(G,A)$ consists of a smooth module 
structure $S \:  G \to \Aut(A)$ and an element $\omega \in
C^2_{ss}(G,A)$. In this case $d_S$ is the Lie group differential 
corresponding to the module structure on $A$ (Definition~I.2). 
Therefore $(S,\omega)$ defines a 
Lie group $N \times_{(S,\omega)} G$ if and only if 
$d_S \omega = \1$, i.e., $\omega \in Z^2_{ss}(G,A)$. 
In this case we write $A \times_\omega  G$ for this Lie group, which is 
$A \times  G$, endowed with the multiplication 
$$ (a,g)(a',g') = (a + g.a' + \omega(g,g'), gg'). $$

Further $S \sim S'$ if and only if 
$S = S'$. Hence a smooth $G$-kernel $[S]$ 
is the same as a smooth $G$-module structure $S$ on $A$ 
and $\Ext(G,A)_S := \Ext(G,A)_{[S]}$ is the class of all $A$-extensions of 
$G$ for which the associated $G$-module structure on $A$ is $S$. 

According to Corollary II.13, the equivalence classes of extensions 
correspond to cohomology classes of cocycles, so that the map 
$$ H^2_{ss}(G,A)_S \to \Ext(G,A)_S, \quad [\omega] \mapsto [A \times_\omega  G] $$
is a well-defined bijection. This was also 
shown directly in [Ne04a] in the context of abelian extensions.   
\qed

\Remark II.15. The result of Corollary~II.13 can also be visualized on the level of 
extensions as follows. For that we assume that $Z(N)$ is an initial Lie subgroup 
of $N$, so that it carries its own Lie group structure and the smooth outer action 
$S$ induces on $Z(N)$ the structure $S_Z$ of a smooth $G$-module (Lemma~I.6). 
In view of Remark~II.14, we then have 
$$ H^2_{ss}(G,Z(N))_{S}  \cong \Ext(G,Z(N))_{S_Z}. $$
Let 
$Z(N) \to \hat G_1 \sssmapright{q_1} G$
be an extension of $G$ by the smooth $G$-module $Z(N)$ and 
$N \to \hat G_2 \sssmapright{q_2} G$
 an extension of $G$ by $N$ corresponding to the $G$-kernel $[S]$. 
Then we consider the group 
$$  H := q_1^*\hat G_2 \cong 
\{ (g_1, g_2) \in \hat G_1 \times \hat G_2 \: q_1(g_1) = q_2(g_2)\} $$
which is a Lie group extension of $G$ by $Z(N) \times N$, and the subgroup 
$\Delta_Z := \{ (z,z^{-1}) \: z \in Z(N)\}$ is a central split Lie subgroup of $H$. 
Therefore the {\it Baer product}
$\hat G := H/\Delta_Z$
is a Lie group extension of $G$ by $N \cong (Z(N) \times N)/\Delta_Z$. 

If $\hat G_1 \cong Z(N) \times_f G$ and $\hat G_2 \cong N \times_{(S,\omega)} G$, 
then 
$$ H \cong (Z(N) \times N) \times_{((S_Z,S), (f, \omega))} G $$
and 
$\hat G \cong N \times_{(S, f \cdot \omega)} G.$ This implies that the action of the abelian 
group $H^2_{ss}(G, Z(N))_S$ on $\Ext(G,N)_{[S]}$ corresponds to the Baer multiplication of 
extensions of $G$ by $Z(N)$ with $N$-extensions of $G$. 
\qed

\Remark II.16. The preceding results show that the fibers of the map 
$[(S,\omega)] \mapsto [S]$ 
that assigns to an extension of $G$ by $N$, represented by a smooth factor system 
$(S,\omega)$, the class $[S]$ of the corresponding smooth outer action, 
are given by 
$Q^{-1}([S]) = \Ext(G,N)_{[S]}$, and the group $H^2_{ss}(G,Z(N))_{S}$ acts simply 
transitively on the fibers. Clearly the action of $G$ on $Z(N)$ depends on $[S]$ and so does 
the cohomology group $H^2_{ss}(G,Z(N))_{S}$. 

In the example described below, we shall see that for an abelian smooth $G$-module 
$A = N = Z(N)$ the group $H^2_{ss}(G,A)_S$ very much depends on 
the action $S$ of $G$ on $A$. 
For $G = \R^2$ and $A = \R$ we have 
$H^2_s(G,A) \cong H^2_c(\g,\a)$
(cf.\ [Ne04a, Th.~7.2]). On the other hand $C^2_c(\g,\a)$ is the space of alternating bilinear maps 
$\g \times \g \to \a$, hence $1$-dimensional. Further $\dim \g = 2$ implies that 
each $2$-cochain is a cocycle. Since $B^2_c(\g,\a) = d_\g \a$ 
vanishes if the module $\a$ is trivial and coincides with 
$Z^2_c(\g,\a)$ otherwise, we have 
$$ H^2_{ss}(G,A) \cong H^2_s(G,A) \cong H^2_c(\g,\a) \cong \cases{ 
\R & for $\g.\a = \{0\}$ \cr 
\{0\} & for $\g.\a \not= \{0\}$. \cr} 
\qeddis 

The following theorem provides a Lie theoretic criterion for the non-emptiness of the 
set $\Ext(G,N)_{[S]}$. 

\Theorem II.17. If $S$ is a strongly smooth outer action of $G$ on $N$, then 
$$ \chi_{ss}(S)= \1 \quad \Longleftrightarrow \quad \Ext(G,N)_{[S]} \not= \eset. $$
If $G$ is connected and $S$ is a smooth outer action, then 
$$ \chi_{s}(S) = \1 \quad \Longleftrightarrow \quad \Ext(G,N)_{[S]} \not= \eset. $$

\Proof. If there exists a Lie group extension $\hat G$ corresponding to 
$[S]$, then we may w.l.o.g.\ assume that it is of the form 
$N \times_{(S,\omega)} G$ for a smooth factor system $(S,\omega)$ 
(Theorem~II.12). This implies in particular that $S$ is strongly smooth, and 
$d_S \omega = \1$ (Lemma~II.7) leads to $\chi_{ss}(S) = [d_S \omega] = \1$. 

Suppose, conversely, that $S$ is strongly smooth with 
$\chi_{ss}(S) = \1$. Then there exists 
$\omega \in C^2_{ss}(G,N)$ with $\delta_S  = C_N \circ \omega$ and 
some $h \in C^2_{ss}(G,Z(N))$ with $d_S\omega = d_S h^{-1}$, so that 
$\omega' := \omega \cdot h \in C^2_{ss}(G,N)$  
satisfies $d_S \omega' = d_S \omega \cdot d_S h = \1$ and 
$\delta_S = C_N \circ \omega'$. Hence $(S,\omega')$ is a smooth factor system, 
and Proposition~II.8 implies the existence of a Lie group extension 
$N \times_{(S,\omega)} G$ of $G$ by $N$ corresponding to $[S]$. 

If $G$ is connected and $S$ is a smooth outer action with 
$\chi_s(S) = \1$, then a similar argument provides some 
$\omega' \in C^2_s(G,N)$ with $\delta_S = C_N \circ \omega'$ and 
$d_S \omega' = \1$. In view of Lemma~II.9, $\omega' \in C^2_{ss}(G,N)$,  and we 
thus obtain a Lie group extension 
$N \times_{(S,\omega')} G$ of $G$ by $N$. 
\qed

\subheadline{Split extensions} 

We conclude this section with a characterization of the {\it split extensions}, i.e., 
those extensions $q \: \hat G \to G$ for which there exists a smooth homomorphism 
$\sigma \: G \to \hat G$ with $q \circ \sigma = \id_G$, which implies that 
$\hat G \cong N \rtimes_S G$ with respect to the smooth action of $G$ on $N$ 
given by $S := C_N \circ \sigma$. 

\Proposition II.18. The extension $\hat G = N \times_{(S,\omega)} G$ is 
split if and only if there exists an $h \in C^1_s(G,N)$ with 
$$ h *_S \omega = \1. \leqno(2.8) $$ 
If, in addition, $\omega \in C^2_{ss}(G,Z(N))$, then $S$ is a representation, and 
{\rm(2.8)} is equivalent to 
$$ d_S h = \omega^{-1}. $$
Let 
$$ Z^1_s(G,N,Z(N))_S := \{ f \in C^1_s(G,N) \: \im(d_S f) \subeq Z(N)\} $$
and observe that 
$$ Z^1_s(G,N,Z(N))_S/C^1_s(G,Z(N)) \to Z^1(G,\Aut(N))_S, \quad [f] \mapsto C_N \circ f $$
is injective onto the set of all $1$-cocycles $G \to C_N(N) \subeq \Aut(N)$ that can be 
lifted to elements of $C^1_s(G,N)$, where 
the action of $C^1_s(G,Z(N))$ on $C^1_s(G,N)$ is given by pointwise multiplication. 
Moreover, 
$$ [\omega].[N \times_{(S,\1)} G] \quad \hbox{ splits } \quad 
\Longleftrightarrow \quad [\omega^{-1}] \in \im(\delta), $$
where 
$$ \delta \: Z^1_s(G,N,Z(N))/ C^1_s(G,Z(N)) \to H^2_{ss}(G,Z(N)) \quad 
[f] \mapsto [d_S f]. $$

\Proof. Since the orbit of the factor system $(S,\omega)$ under 
the group $C^1_s(G,N)$ consists of all other factor systems 
describing the same extension, i.e., corresponding to different choices of 
sections $G \to \hat G$ (Proposition~II.10), 
the extension $\hat G$ splits if and only if 
this orbit contains a factor system $(S',\omega') = h.(S,\omega)$ for which the 
corresponding section $\sigma' \: G \to \hat G, g \mapsto (h(g),g)$ is a 
group homomorphism, i.e., $\delta_{\sigma'} = \1$. 
Hence the first assertion follows from $\delta_{\sigma'} = h *_S \omega$. 

If, in addition, $\omega$ has values in $Z(N)$, then $\delta_S = C_N \circ \omega = \1$ 
shows that $S$ is a representation. Moreover, 
$h *_S \omega = d_S h \cdot \omega$ vanishes if and only if $d_S h = \omega^{-1}$. 
\qed

\Remark II.19. In general the image of the map 
$\delta$ is a quite complicated subset of the cohomology group 
$H^2_{ss}(G,Z(N))_S$. 
To understand this subtle point, we note that $d_S h = \omega^{-1}$ does not seem to imply 
$d_S h^{-1} = \omega$. It is true that for $S_h = (C_N \circ h) \cdot S$ we have 
$$ d_{S_h} h^{-1} = 
h^{-1} *_{S_h} \1 = h^{-1} *_{S_h} (h *_S \omega) 
= \omega $$
because $(S,\omega) = h^{-1}.(h.(S,\omega)) = h^{-1}.(S_h,\1)$. 
The condition $d_S h^{-1} = d_{S_h} h^{-1}$ is equivalent to 
$$S(g)(h(g')) = h(g) S(g)(h(g')) h(g)^{-1} $$ 
for $g,g' \in G$, which means that the two functions $h$ and $S.h$ commute pointwise. 

For the special case $S = \1$ we have $d_S h = \delta_h$ and 
$\im(\delta_h) \subeq Z(N)$ means that the map $C_N \circ h \: G \to \Aut(N)$ 
is a homomorphism. Fix some $g > 1$ and let us consider the special case, where 
$N = \tilde\SL_2(\R)$ is the universal covering group of $\SL_2(\R)$ 
and $G$ is the discrete group with 
$2g$ generators 
$\alpha_1,\ldots, \alpha_{2g}$ subject to the commutator relation 
$$ [\alpha_1, \alpha_2] \cdots [\alpha_{2g-1}, \alpha_{2g}] 
= \alpha_1 \alpha_2 \alpha_1^{-1} \alpha_2^{-1} \cdots 
\alpha_{2g-1} \alpha_{2g} \alpha_{2g-1}^{-1} \alpha_{2g}^{-1} = \1. $$
In this case $Z(N) \cong \Z$ is discrete and $N_{\rm ad} := N/Z(N) \cong \PSL_2(\R)$. 
Then each homomorphism $h \: G \to N_{\rm ad}$ corresponds to a 
$(2g)$-tuple of points $(x_1, \ldots, x_{2g}) \in \PSL_2(\R)^{2g}$ satisfying 
$$ [x_1, x_2] \cdots [x_{2g-1}, x_{2g}] = \1. $$

Since $G$ and $Z(N)$ are discrete, we have 
$$ H^2_{ss}(G,Z(N))_S = H^2(G,Z(N)). $$
Using the fact that $G$ is the fundamental group $\pi_1(\Sigma_g)$ 
of a compact orientable surface $\Sigma_g$ 
of genus $g$, it is known that $H^2(G,Z(N)) \cong H^2(\Sigma_g, Z(N)) 
\cong \Hom(H_2(\Sigma_g),Z(N)) \cong Z(N)$ because $H_2(\Sigma_g) \cong \Z$ ([Mil58]). 
The isomorphism 
$$ \Phi \: H^2(G,Z(N)) \to Z(N) $$
can be obtained by choosing for each central extension 
$Z(N) \into \hat G \to G$ a section $\sigma \: G \to \hat G$ and then putting 
$$ \delta([\hat G]) := 
[\sigma(\alpha_1), \sigma(\alpha_2)] \cdots [\sigma(\alpha_{2g-1}), \sigma(\alpha_{2g})]. $$
From this observation it follows that 
$\im(\delta) \subeq Z(N)$ 
coincides with the set 
$$ \{ [y_1,y_2] \cdots [y_{2g-1}, y_{2g}] \: 
(y_1, \ldots, y_{2g}) \in N, [y_1,y_2] \cdots [y_{2g-1}, y_{2g}] \in Z(N)\}. $$
In [Mil58] it is shown that, as a subset of $\Z \cong Z(N)$, this set coincides with 
$\{ n \in \Z \: |n| < g\}$. In particular it is not a subgroup. 
\qed

\Remark II.20. (a) Suppose that the extension 
$\hat G = N \times_{(S,\omega)} G$ splits, so that we may w.l.o.g.\ assume that 
$\omega = \1$. For $h \in C^1_s(G,N)$ we then have 
$$ h.(S,\1) = (h.S, h *_S \1) = (h.S, d_S h) 
\quad \hbox{ with } \quad 
d_S h(g,g') = h(g) S(g)(h(g')) h(gg')^{-1}. $$
We conclude that 
$$ \{ h \in C^1_s(G,N) \: h *_S \1 = \1\} 
= Z^1_s(G,N)_S $$
is the set of $1$-cocycles for the smooth action of $G$ on $N$ defined by the homomorphism 
$S \: G \to \Aut(N)$. Therefore the orbit of $(S,\1)$ under $C^1_s(G,N)$ 
may contain different elements of the form $(S',\1)$. The representations 
$S' \: G \to \Aut(N)$ arising that way have the form 
$$ S_f := (C_N \circ f) \cdot S, \quad f \in Z^1_s(G,N)_S. $$

If $N$ is abelian, all representations $S_f$ are the same, so that 
the $C^1_s(G,N)$-orbit of $(S,\1)$ contains no other element of the form 
$(S',\1)$. But if $f \in Z^1_s(G,N)_S$ is a cocycle whose values do not 
lie in the center $Z(N)$ of $N$, then $S_f \not= S$. Nevertheless, 
the two split extensions 
$$ N \rtimes_S G  \quad \hbox{ and } \quad N \rtimes_{S_f} G $$
are equivalent. 

 (b) A condition that is slightly weaker than the condition 
$d_S h = \1$ is that $S_h = (C_N \circ h) \cdot S$ is a homomorphism, which, in view of 
$\delta_{h.S} = C_N \circ (d_S h)$ is equivalent to $\im(d_S h) \subeq Z(N)$, 
which in turn means that $C_N \circ h \in Z^1(G,\Aut(N))$ with respect to the 
action of $G$ on $\Aut(N)$ given by $g.\phi := S(g) \phi S(g)^{-1}$. 
Then $d_S d_S h = d_{S_h} d_S h = d_{S_h}(h *_S \1) = \1$ (Lemma~I.10(4)) implies that 
$d_S h \in Z^2_{ss}(G,Z(N))_S$, so that on the level of extension classes we have 
$$ [N \rtimes_S G] 
= [N \times_{(h.S, d_S h)} G] 
= [d_S h].[N \times_{(S_h, \1)} G]
= [d_S h].[N \rtimes_{S_h} G]. $$
This shows that the semidirect products 
$$ N \rtimes_S G  \quad \hbox{ and } \quad N \rtimes_{S_h} G $$
are equivalent extensions of $G$ by $N$ if and only if 
the class $[d_S h] \in H^2_{ss}(G,Z(N))_S$ vanishes. 
\qed

\Corollary II.21. Suppose that $Z(N) \trile N$ is a split Lie subgroup, so that 
$N_{\rm ad} := N/Z(N)$ carries a canonical Lie group structure for which 
$q_N \: N \to N_{\rm ad} = N/Z(N)$ defines a central extension of $N_{\rm ad}$ by $Z(N)$. 
Further let $S \: G \to \Aut(N)$ be a smooth action of $G$ on $N$ 
and consider the induced action $\oline S$ of $G$ on $N_{\rm ad}$. 
Then we have an exact sequence of groups, resp.\ pointed spaces, 
$$ Z(N)^G \into N^G \to N_{\rm ad}^G \to 
H^1_s(G,Z(N))_S \to  
H^1_s(G,N)_S \to  H^1_s(G,N_{\rm ad})_{\oline S} 
\sssmapright{\delta} H^2_{ss}(G,Z(N))_{S},  $$
where for 
$f \in Z^1_s(G,N_{\rm ad})_S$ we have 
$$ \delta([f]) =[d_S \hat f] $$
for some $\hat f \in Z^1_s(G,N, Z(N))_S$ satisfying $q_N \circ \hat f = f$ 
and 
$$ \Ext(G,N)_{[S],{\rm split}} = -\im(\delta).[N \rtimes_S G]. $$

\Proof. For the $7$-term exact sequence and the connecting map $\delta$ 
corresponding to the central Lie 
group extension $Z(N) \into N \onto N_{\rm ad}$, we refer to [Ne05, Sect.~II]. 

If the cocycle 
$f \: G \to N$ has values in $Z(N)$, then $f.(S,\1) = (S,\1)$, so that 
$Z^1_s(G,Z(N))_S$ corresponds to the stabilizer of $(S,\1) \in Z^2_{ss}(G,N)$ and the map 
$$ Z^1_s(G,N, Z(N))_S \to Z^2_{ss}(G,N), \quad f \mapsto f.(S,\1) $$
factors through an injective map 
$$ Z^1_s(G,N, Z(N))_S/Z^1_s(G,Z(N))_S 
\cong Z^1_s(G, N_{\rm ad}), \quad [f] \mapsto f.(S,\1) $$
because for each $h \in Z^1_s(G,N_{\rm ad})_S$ there is an 
$\hat h \in Z^1_s(G,N, Z(N))_S$ with $q_N \circ \hat h = h$ and 
$\hat h.(S,\1) = (\hat h.S, d_S \hat h)$. 

Since 
$d_S \hat h \in Z^2_{ss}(G,Z(N))_S$, this factor system corresponds to the extension class 
$$ [d_S \hat h].[N \rtimes_S G] = \delta(h).[N \rtimes_S G] \in \Ext(G,N)_{[S]}. $$
This class equals $[N \rtimes_S G]$ if and only if $\delta(h)  = 0$ which is equivalent to 
$[h]$ being contained in the image of $H^1_s(G,N)_S$ under the connecting map $\delta$. 

The remaining assertion follows from Remark~II.20(b) and 
Proposition~II.18. 
\qed

In [Bo93], Borovoi deals with non-abelian $H^2$-sets in the context of Galois 
cohomology, where he calls the cohomology classes in $H^2_{ss}(G,N)$ corresponding to split 
extensions {\it neutral}. Our Corollary~II.21 is a Lie group 
version of Borovoi's Proposition~2.3.

\sectionheadline{III. Smooth crossed modules} 

\nin In this section we introduce the notion of a smooth crossed module for a Lie group $G$. 
Our point of view is that a smooth crossed module is a 
central extension $\hat N \to N$ of a normal split Lie subgroup 
$N \trile G$ for which the conjugation action 
of $G$ on $N$ lifts to a smooth action on $\hat N$. 
It turns out that smooth crossed modules of Lie groups 
provide the natural framework to reduce problems related to general extensions 
of Lie groups to abelian extensions, which is carried out in Section~IV. 

\Definition III.1. A morphism $\alpha \: H \to G$ of Lie groups, 
together with a homomorphism $\hat S \: G \to \Aut(H)$ defining a smooth action 
$\hat S \: G \times H \to H, (g,h) \mapsto g.h = \hat S(g)(h)$ 
of $G$ on $H$, is called a {\it smooth crossed module} if the following conditions 
are satisfied: 
\litemindent=1.0cm
\litem{(CM1)} $\alpha \circ \hat S(g) = c_{\alpha(h)} \circ \alpha$ for $g \in G$. 
\litem{(CM2)} $\hat S \circ \alpha = C_H \: H \to \Aut(H)$ is the conjugation action. 
\litem{(CM3)} $\ker(\alpha)$ is a split Lie subgroup of $H$ and 
$\im(\alpha)$ is a split Lie subgroup of $G$ for which 
$\alpha$ induces an isomorphism $H/\ker(\alpha) \to \im(\alpha)$. 
\litemindent=0.7cm

\nin The conditions (CM1/2) express the compatibility of the $G$-action on 
$H$ with the conjugation actions of $G$ and $H$. 
\qed

\Lemma III.2. If $\alpha \: H \to G$ is a smooth crossed module, then 
the following assertions hold: 
\litem{(1)} $\im(\alpha)$ is a normal subgroup of $G$. 
\litem{(2)} $\ker(\alpha) \subeq Z(H)$. 
\litem{(3)} $\ker(\alpha)$ is $G$-invariant and $G$ acts smoothly on $\ker(\alpha)$. 

\Proof. (1) follows from (CM1), (2) from (CM2), and 
(3) from (CM1), Lemma~I.7(a) and Lemma II.3. 
\qed

Crossed modules for which $\alpha$ is injective are 
inclusions of normal split Lie subgroups 
and surjective crossed modules are central extensions. 
In this sense the concept of a crossed module interpolates between split 
normal subgroups and central extensions. 

In the following we shall adopt the following perspective on crossed modules. 
Let $\alpha \: H \to G$ be a smooth crossed module. 
Then $N := \im(\alpha)$ is a split normal Lie subgroup of 
$G$ and $\alpha \: H \to N$ is a 
central Lie group extension of $N$ by $Z := \ker(\alpha)$. 
In this sense a smooth crossed module can be viewed as a 
central Lie group extension 
$\alpha \: H \to N$ of a split normal subgroup $N$ of $G$ 
for which there 
exists a smooth $G$-action on  $H$ satisfying (CM1/2). 

If, conversely, $N$ is a split normal subgroup of 
$G$ and $\alpha \: \hat N \to N$ a central Lie group extension 
of $N$ by $Z$, then we have a natural 
smooth $N$-module structure on $\hat N$ given by 
$\hat C_N(n)(n') := nn'n^{-1}$. To obtain the structure of a smooth 
crossed module for $\alpha \: \hat N \to G$, one has to 
extend the action of $N$ on $\hat N$ 
to a smooth $G$-action on  $\hat N$ for which $\alpha$ is 
equivariant. In the following we shall adopt this point of view.

We may w.l.o.g.\ assume that $\hat N = Z \oplus_f N$ for some 
$f \in Z^2_{ss}(N,Z)$ (Remark II.14) and write the 
$ G $-module structure on $ Z$ as $(g,z) \mapsto g.z$. Then the 
$ G $-module structure on $\hat N  =  Z \oplus_f  N $ can be written
as 
$$ g.(z,n) = (g.z + \theta(g)(g.n), c_g(n)), \leqno(3.1) $$
where the function 
$\theta \:  G  \to C^1_s(N,Z)$
is a map for which $\tilde\theta \: G \times N \to Z, (g,n) \mapsto \theta(g)(n)$
is smooth on an open neighborhood of $G \times \{\1\}$ in $G \times N$. 
That (3.1) defines a group action of $G$ on $\hat N$ 
means that $\theta$ is a cocycle with respect to the natural action of $G$ 
on $C^1_s(N,Z)$ by $(g.\alpha)(n) := g.\alpha(g^{-1}ng)$, and that with
respect to the natural action of $G$ on $C^2_s(N,Z)$ by 
$(g.\beta)(n_1,n_2) := g.\beta(g^{-1}n_1g, g^{-1}n_2g)$ we have 
$$ d_N(\theta(g)) = g.f - f = d_S(f)(g) \quad \hbox{ for all } g \in
G. \leqno(3.2) $$
(cf.\ Corollary~A.4 and Lemma~A.11). 

\Lemma III.3. For a cocycle 
$\theta \: G \to C^1_s(N,Z)$ satisfying {\rm(3.2)}, the following are equivalent: 
\litem{(1)} The corresponding action of $G$ on $\hat N$ is smooth. 
\litem{(2)} The function $\tilde\theta : G \times N \to Z, (g,n)
\mapsto \theta(g)(n)$ is smooth in an identity neighborhood and for each 
$n \in  N$ the function 
$$ \theta_n \: G \to Z, \quad g \mapsto \theta(g)(g.n) + f(g.n,n^{-1}) \leqno(3.3) $$
is smooth in an identity neighborhood of $G$.  

\Proof. This is an immediate consequence of Proposition A.13. 
\qed

\Definition III.4. In the following we write $C^1_s(G, C^1_s(N,Z))$ for the set of all maps 
$\alpha \: G \to C^1_s(N,Z)$ for which the function 
$\tilde\alpha : G \times N \to Z, (g,n)\mapsto \alpha(g)(n)$ is smooth in an identity neighborhood. 
We write $C^1_{ss}(G,C^1_s(N,Z))_f$ for the set of all elements with the additional property 
that for each $n \in  N$ the function 
$$ \alpha_n \: G \to Z, \quad g \mapsto \alpha(g)(g.n) + f(g.n, n^{-1})$$
is smooth in an identity neighborhood. 
The set of all cocycles in $C^1_s(G,C^1_s(N,Z))$ is denoted 
$Z^1_s(G,C^1_s(N,Z))$ and likewise 
$Z^1_{ss}(G, C^1_s(N,Z))_f := 
C^1_{ss}(G, C^1_s(N,Z))_f \cap Z^1_{s}(G,C^1_s(N,Z))$. 

The natural action of $N$ on $\hat N$ is given by 
$$ (z,n)(z',n') (z,n)^{-1} 
= \big(z' + f(n,n')-f(nn'n^{-1},n),nn'n^{-1}\big). $$
We define 
$$ \tilde f(n)(n') := f(n,n^{-1}n'n)-f(n',n), \leqno(3.4) $$
and observe that (3.1) defines a smooth action extending the
given action of $N$ if $\theta \res_N = \tilde f$, i.e.,
that the pair $(f,\theta)$ is an element of the group 
$$ D := \{ (f,\theta) \in Z^2_{ss}(N,Z) \times Z^1_{s}(G,C^1_s(N,Z)) \: 
d_N \circ \theta = d_S(f), \theta\res_N = \tilde f, 
\theta \in Z^1_{ss}(G,C^1_s(N,Z))_f\}. $$
That the set of all these elements is indeed a subgroup 
follows directly from Lemma~III.3. We also write 
$$ D(f) := \{ \theta \in Z^1_{ss}(G,C^1_s(N,Z))_f \: (f,\theta) \in D \} = \pr_2(D) $$
for the set of all smooth actions of $G$ on $\hat N$ extending the given action of $N$. 
\qed

The following lemma provides a convenient description of the set $D(f)$ as a homogeneous 
space. 

\Lemma III.5. For any $(f,\theta) \in D$ and $\alpha \in 
Z^1_{ss}(G/N,Z^1_s(N,Z))_0$  we have 
$$ \alpha * (f,\theta) := (f, \theta + q^*\alpha) \in D. $$
This defines for each $f \in \pr_1(D)\subeq Z^2_{ss}(N,Z)$ 
a simply transitive action of $Z^1_{ss}(G/N,Z^1_s(N,Z))_0$ on the set~$D(f)$. 

\Proof. That $\alpha * (f,\theta) \in D$ follows from 
$$ d_N \circ (\theta + q^*\alpha) 
= d_N \circ \theta + d_N \circ q^*\alpha 
= d_N \circ \theta + q^*(d_N \circ \alpha) = d_N \circ \theta = d_S(f). $$

If, conversely, $(f,\theta), (f,\theta') \in D$, then the difference 
$\theta' -\theta$ has values in $Z^1_{s}(N,Z)$ and vanishes on $N$, hence is of the 
form $q^*\alpha$ for some $\alpha \in Z^1_{ss}(G/N, Z^1_s(N,Z))_0$ because 
for each $n  \in N$ the map 
$$ g \mapsto (\theta' - \theta)(g)(g.n) 
= \theta'(g)(g.n) + f(g.n,n^{-1})  - \big(\theta(g)(g.n) - f(g.n, n^{-1})\big) $$
is smooth in an identity neighborhood of $G$. This proves the lemma. 
\qed

We want to construct a map 
$Q \: D  \to H^3_{s}(G/N,Z)$
such that $Q(f,\theta) = 0$ characterizes the extendability of $f$ 
to a cocycle $\tilde f \in Z^2_{s}(G,Z)$ satisfying 
$\tilde f\res_{G \times N} = \tilde\theta$. 
In the following we fix a pair $(f,\theta) \in D$ 
and consider the corresponding smooth 
action of $G$ on $\hat N$ given by 
$$  \hat S(g)(z,n) := (g.z + \theta(g)(g.n), g.n) $$
(cf.\ Corollary~A.4 and Proposition A.13). To define the map $Q$, let 
$\sigma \: G/N \to G$ be a normalized locally smooth section of the
quotient map and define 
$$ S := \hat S \circ \sigma \: G/N \to \Aut(\hat N). $$

There exists an $\omega \in C^2_{s}(G/N,\hat N)$ with 
$q_N \circ \omega = \delta_\sigma$, 
where $q_N \: \hat N \to N$ is the quotient map. Then 
$$ \delta_S = \hat S \circ \delta_\sigma = \hat S \circ q_N \circ \omega 
= C_{\hat N} \circ \omega, \leqno(3.5) $$
which shows that $S$ is a smooth outer action of $G/N$ on $\hat N$. 

If we even find an $\omega \in C^2_{ss}(G/N,\hat N)$ with 
$\delta_S = C_N \circ \omega$, we have a strongly smooth outer action, 
but this does not follow from the definition of a smooth crossed 
module. Using Lemma~I.10, it is easy to see that this additional smoothness condition 
for $S$ does not depend on the choice of the section $\sigma$, hence is a property of 
the pair $(f,\theta)$. We write $D_{ss} \subeq D$ for the set of all pairs 
$(f,\theta) \in D$ for which $S$ is a strongly smooth outer action. 
In terms of the product coordinates $\hat N = Z \times_f N$, this means that there 
exists $\omega_Z \in C^2_{ss}(G/N, Z)$ such for 
$\omega := (\omega_Z, \delta_\sigma) \in C^2_s(G/N,\hat N)$ the functions 
$\omega_g \: G/N \to \hat N$ 
given by 
$$ \eqalign{ 
\omega_g(x) &= \omega(g,x) \omega(gxg^{-1},g)^{-1} \cr
&= (\omega_Z(g,x), \delta_\sigma(g,x)) 
(\omega_Z(gxg^{-1},g), \delta_\sigma(gxg^{-1},g))^{-1} \cr
&= (\omega_Z(g,x), \delta_\sigma(g,x)) 
\big(-\omega_Z(gxg^{-1},g) - f(\delta_\sigma(gxg^{-1},g)^{-1}, \delta_\sigma(gxg^{-1},g)),  
\delta_\sigma(gxg^{-1},g)^{-1}\big) \cr
&= (\omega_{Z,g}(x) - f(\delta_\sigma(gxg^{-1},g)^{-1}, \delta_\sigma(gxg^{-1},g)) 
+ f(\delta_\sigma(g,x), \delta_\sigma(gxg^{-1},g)^{-1}), \delta_{\sigma,g}(x)\big)\cr} $$
are smooth in an identity neighborhood. For $\delta_{\sigma,g}$ this follows as in 
Definition~I.1, but for the $Z$-component this depends on the values of the section 
$\sigma$ in elements far from the identity, so that the smoothness assumptions in 
the identity do not imply anything. 

\Lemma III.6. The function 
$$ \tilde Q(f,\theta,\sigma,\omega) := d_S \omega \: (G/N)^3 \to \hat N, 
\quad (g,g',g'') \mapsto 
S(g)(\omega(g',g'')) \omega(g,g'g'') \omega(gg',g'')^{-1} \omega(g,g')^{-1}
$$
has values in $Z$ and is a $3$-cocycle in $Z^3_s(G/N,Z)$. 
Its cohomology class $[d_S\omega] \in H^3_{s}(G/N,Z)$ does not
depend on the choices of $\omega$ and $\sigma$. Likewise,  for 
$(f,\theta) \in D_{ss}$ the cohomology class of $d_S\omega$ in $H^3_{ss}(G/N,Z)$ 
does not depend on $\omega$ and $\sigma$. 

\Proof. First we observe that $(d_S\omega)(g,g',g'') \in Z$ follows from 
$$ \eqalign{ 
q_N((d_S\omega)(g,g',g'')) 
&= C_N(\sigma(g))(\delta_\sigma(g',g'')) \delta_\sigma(g,g'g'')\delta_\sigma(gg',g'')^{-1} 
\delta_\sigma(g,g')^{-1} \cr 
&= \sigma(g)\big(\sigma(g')\sigma(g'')\sigma(g'g'')^{-1}\big)\sigma(g)^{-1} 
\sigma(g)\sigma(g'g'')\sigma(gg'g'')^{-1}\sigma(gg'g'')\cr
&\ \ \ \
\sigma(g'')^{-1}\sigma(gg')^{-1}\sigma(gg')\sigma(g')^{-1}\sigma(g)^{-1}
= \1.\cr} $$

It is easy to verify that $d_S \omega$ vanishes if one of its three
arguments is~$\1$, so that $d_S \omega \in C^3_s(G/N,Z)$. 
Since $S$ is a smooth outer action of $G/N$ on $\hat N$ and 
$\delta_S = C_{\hat N} \circ \omega$, Lemma~I.10(5) implies that 
$d_S \omega \in Z^3_s(G/N, Z(\hat N))$. Since its values lie in the subgroup 
$Z$ of $Z(\hat N)$, we have $d_S \omega \in Z^3_s(G/N, Z)$. 

It remains to show that the cohomology class of $[d_S \omega]$ in 
$H^3_{s}(G/N, Z)$, resp., $H^3_{ss}(G/N,Z)$ 
does not depend on the choices of $\omega$ and $\sigma$. 
We only discuss the case $(\theta, f) \in D_{ss}$. 
First we show that is does not depend on the
choice of $\omega$ if $\sigma$ is fixed. 
So let $\tilde\omega \in C^2_{ss}(G/N,\hat N)$ be another map 
with $q_N \circ \tilde\omega = \delta_\sigma$. Then there exists a function
$\beta \in C^2_{ss}(G/N,Z)$ with 
$\tilde\omega = \beta \cdot \omega$. As all values of $\beta$ commute
with all values of $\omega$, this leads to 
$$ d_S \tilde\omega 
= d_S \beta + d_S \omega 
= d_{G/N} \beta + d_S \omega \in d_S \omega + d_{G/N}(C^2_{ss}(G/N,Z)). $$
Therefore the choice of $\omega$ for a fixed section $\sigma$ 
has no effect on the cohomology class $[d_S\omega]$ in $H^3_{ss}(G/N,Z)$. 

Now we consider another section $\tilde\sigma \in C^1_s(G/N,G)$ of $q$. 
Then there exists a function $h\in C^1_s(G/N,N)$ with 
$\tilde\sigma = h \cdot \sigma$ and 
$$ \delta_{\tilde\sigma}(g,g') 
= h(g) \sigma(g)h(g') \sigma(g')\big(h(gg') \sigma(gg')\big)^{-1}
= h(g) c_{\sigma(g)}(h(g')) \delta_\sigma(g,g')
h(gg')^{-1}. $$
Let $\hat h \in
C^1_s(G/N,\hat N)$ with $q_N \circ \hat h = h$. 
Then the prescription 
$$ \tilde\omega(g,g') 
= \hat h(g) S(g)(\hat h(g'))
\omega(g,g')\hat h(gg')^{-1} = (\hat h *_S \omega)(g,g') $$
defines an element $\tilde \omega \in C^2_{ss}(G/N,\hat N)$ 
for which $\tilde S := \hat S \circ \tilde \sigma = (C_{\hat N} \circ \hat h) \cdot S$ 
satisfies $\delta_{\tilde S} = C_{\hat N} \circ \tilde \omega$ 
(Lemma~I.10(1)), and we conclude from Lemma~I.10(4) that 
$d_{\tilde S} \tilde\omega = d_S \omega$. 
\qed

The preceding lemma implies that the cohomology class 
$$ Q(f,\theta) := [\tilde Q(f,\theta,\sigma,\omega)] = [d_S\omega]$$
in $H^3_{s}(G/N,Z)$ only depends on the pair $(f,\theta)$, hence defines a map 
$Q \: D \to H^3_{s}(G/N,Z).$
In the strongly smooth case we obtain a map 
$Q_{ss} \: D_{ss} \to H^3_{ss}(G/N,Z).$

\Lemma III.7. The maps 
$Q \: D \to H^3_{s}(G/N,Z)$ and 
$Q_{ss} \: D_{ss} \to H^3_{ss}(G/N,Z)$ are group homomorphisms. 

\Proof. Again we consider only the case $(\theta,f ) \in D_{ss}$. We  
choose $\omega \in C^2_{ss}(G/N,\hat N)$ with $q_N \circ \omega = \delta_\sigma$. 
We have to see how $Q(f,\theta, \sigma, \omega)$ 
depends on $(f,\theta)$. Writing $\hat N := Z \times_f N$, any $\omega$ 
with $q_N \circ \omega = \delta_\sigma$ can be written 
as $\omega = (\omega_Z, \delta_\sigma)$, independently of $f$ and $\theta$. 

Writing $g.n = gng^{-1}$ for the conjugation action of $G$ on $N$, we now obtain 
$$ \eqalign{ 
S(g).\omega(g',g'') 
&= \hat S(\sigma(g))(\omega_Z(g',g''), \delta_\sigma(g',g'')) \cr
&= (g.\omega_Z(g',g'') + \theta(\sigma(g))(\sigma(g).\delta_\sigma(g',g'')), 
\sigma(g).\delta_\sigma(g',g'')). \cr}$$
Therefore 
$$ \eqalign{ 
&\ \ \ \ \big(S(g).\omega(g',g'')\big)\omega(g,g'g'') \cr
&= \Big(g.\omega_Z(g',g'') + \omega_Z(g,g'g'') + 
\theta(\sigma(g))(\sigma(g).\delta_\sigma(g',g'')) 
+ f(\sigma(g).\delta_\sigma(g',g''), \delta_\sigma(g,g'g'')), \cr
&\ \ \ \ \ \ \ \ \ \ \ \ \ \big(\sigma(g).\delta_\sigma(g',g'')\big) \delta_\sigma(g,g'g'')\Big), \cr} $$
and similarly 
$$ \eqalign{ 
\omega(g,g') \omega(gg',g'') 
&= (\omega_Z(g,g') + \omega_Z(gg', g'')  
+ f(\delta_\sigma(g,g'), \delta_\sigma(gg',g'')), 
\delta_\sigma(g,g')\delta_\sigma(gg',g'')).\cr} $$
This implies that 
$$ \eqalign{ 
(d_S \omega)(g,g',g'')
&= (d_S \omega_Z)(g,g',g'') 
+ \theta(\sigma(g))(\sigma(g).\delta_\sigma(g',g'')) \cr
&\ \ \ \ + f(\sigma(g).\delta_\sigma(g',g''), \delta_\sigma(g,g'g'')) 
- f(\delta_\sigma(g,g'), \delta_\sigma(gg',g'')). \cr} $$
Hence the cohomology class of $d_S \omega$ depends additively on 
the pair $(f,\theta)$, i.e., $Q_{ss}$ is a group homomorphism $D_{ss} \to 
H^3_{ss}(G/N,Z)$. 
\qed

The following theorem explains the meaning of the condition 
$Q(f,\theta) = \1$. It is the main result of the present section and establishes 
an important link between smooth crossed modules and the extendability of central 
extensions of normal subgroups to the whole group. 

\Theorem III.8. Let $f \in Z^2_{ss}(N,Z)$ with $\hat N = Z \times_f N$
and $\theta \in Z^1_{ss}(G,C^1_s(N,Z))_f$ describe a smooth action $\hat S$ of
$G$ on $\hat N$ such that the homomorphism 
$\alpha \: \hat N  \to  G, (z,n) \mapsto n$ is a smooth crossed module. Further let 
$\sigma \in C^1_s(G/N, G)$ with $q \circ \sigma = \id_{G/N}$ be a locally smooth 
section, $\omega \in
C^2_{ss}(G/N,\hat N)$ with $\delta_\sigma = q_N \circ \omega$, and $S = \hat S \circ \sigma$. 
Then the following are equivalent:
\litem{(1)} $Q_{ss}(f,\theta) = [d_S\omega]= \1$ in $H^3_{ss}(G/N,Z)$. 
\litem{(2)} There exists a Lie group extension 
$ Z \into \hat G  \sssmapright{q} G$ and a $G$-equivariant 
equivalence \break $\hat N \to q^{-1}(N)$ of $Z$-extensions of $N$. 
\litem{(3)} There exists a cocycle $f_G \in Z^2_{ss}(G,Z)$ with 
$f_G \res_{N \times N} = f$ and $\tilde f_G \res_{G \times N} = \theta.$ 

\Proof. (1) $\Rarrow$ (2): Suppose that $[d_S\omega] =\1$, i.e., 
$d_S\omega = d_S \beta$ for some $\beta \in C^2_{ss}(G/N,Z)$. 
Then $\omega' := \omega \cdot \beta^{-1}$ satisfies $q_N \circ \omega'
= \delta_\sigma$ and 
$d_S \omega' = d_S \omega - d_S\beta = 0.$
Since $\omega' \in C^2_{ss}(G/N, \hat N)$
 follows from $\beta \in C^2_{ss}(G/N,Z)$, the pair 
$(S,\omega')$ is a smooth factor system for $G/N$ and $\hat N$, so that we obtain a group 
$\hat G = \hat N \times_{(S,\omega')} G/N$ with the product 
$$ (n,g)(n',g') := (n S(g)(n') \omega'(g,g'), gg') $$
and a natural Lie group structure (Proposition~II.8).
Then $\hat G$ is a Lie group extension of $G/N$ by $\hat N$. On the other hand,  
we can view $(C_N \circ \sigma, \delta_\sigma)$ as a smooth factor
system for $G/N$ and $N$, which leads to the description of $G$ as 
$N \times (G/N)$, endowed with the product 
$$ (n,g)(n',g') := (n C_N(\sigma(g))(n') \delta_\sigma(g,g'), gg'). $$
Hence the map 
$$ q_G \: \hat G = \hat N \times (G/N) = (Z \times_f N) \times (G/N)
\to G = N \times (G/N), \quad ((z,n),g) \mapsto (n,g) $$
defines a Lie group extension of $G$ by $Z$ containing 
$\hat N$. It remains to verify that the action of $G$ on $\hat N$
induced by the conjugation action of $\hat G$, for which
$Z$ acts trivially, coincides with the action given by $\theta$. 

It follows from the construction, that the 
conjugation action of $\hat G$ on $\hat N = Z \times_f N$ is given by 
$$ \eqalign{ (n,g)(z,n') (n,g)^{-1} 
&= c_{\hat N}(n)S(g) (z,n') 
= c_{\hat N}(n) (g.z + \theta(g)(C_N(g)(n')), C_N(g)(n')) \cr
&= (g.z + \theta(g)(C_N(g)(n')) + \theta(n)(C_N(n)C_N(g)(n')), C_N(n)C_N(g)(n')) \cr
&= (g.z + \theta((n,g))(C_N(n)C_N(g)(n')), C_N((n,g))(n')) 
= (n,g).(z,n'). \cr}$$ 
This completes the proof of (2). 

 (2) $\Rarrow$ (3): Let $q \: \hat G \to G$ be an extension as
in (2). We may w.l.o.g.\ assume that $\hat G = Z \times_{f_G} G$ for some 
$f_G \in Z^2_{ss}(G,Z)$, where $f_G$ extends $f$, i.e., $\hat N = Z
\times_f N$ is contained in~$\hat G$. The condition that the inclusion
$\hat N \into \hat G$ is $G$-equivariant means that 
$\tilde f_G\res_{G \times N} = \theta$. 

 (3) $\Rarrow$ (1): Let $f_G \in Z^2_{ss}(G,Z)$ be a cocycle as
in (3) and 
$$ q_G \: \hat G := Z \times_{f_G} G \to G$$
the corresponding Lie group extension. The condition$f_G \res_{N \times N} =
f$ means that the group $\hat N = Z \times_f N$ is a split Lie subgroup of $\hat
G$, and the second condition $\tilde f_G \res_{G \times N} = \theta$
ensures that the conjugation action $C_{\hat N}$ of $\hat G$ on $\hat N$ induces
the  action of $G$ on $\hat N$ defined by $\theta$. 

We consider the locally smooth normalized section 
$\hat\sigma \: G/N \to \hat G, g \mapsto (0, \sigma(g)).$
Then 
$\delta_{\hat\sigma}$ satisfies $q_N \circ \delta_{\hat \sigma}= \delta_\sigma$. 
In view of $d_S \omega = \1$  and $S = C_{\hat N} 
\circ \hat \sigma$, we now have 
$$ Q(f,\theta,\sigma,\delta_{\hat\sigma}) = d_S \delta_{\hat\sigma} = \1. $$
It follows in particular that $Q(f,\theta) = [d_S\delta_{\hat\sigma}] = \1$ 
(Lemma~III.6). 
\qed

\Remark III.9. Suppose that $G$ is connected and that $(f,\theta) \in D$. 
If $Q(f,\theta) = \1$ in $H^3_s(G/N,Z)$, then the same argument as in the 
proof of Theorem~III.8 provides 
$\omega' \in C^2_s(G/N, \hat N)$ with  
$\delta_S = C_{\hat N} \circ \omega'$ and $d_S \omega' = \1$. 
Now Lemma~II.9 implies that $\omega' \in C^2_{ss}(G/N,\hat N)$ and that 
(2) still holds. All other implications do not require the strong 
smoothness of $S$. Hence Theorem~III.8 remains true for connected groups $G$ if 
we replace the condition in (1) by the weaker requirement 
$Q(f,\theta) = \1$ in $H^3_s(G/N,Z)$. 
\qed

\sectionheadline{IV. Applications to non-abelian extensions of Lie groups} 

\nin Let $G$ and $N$ be Lie groups, and assume that $Z(N)$ is a split Lie subgroup 
of $N$, so that $N_{\rm ad} := N/Z(N)$ carries a natural Lie group structure. 
Further let $S \in C^1_s(G,\Aut(N))$ be a strongly smooth outer action. 
In this section we describe the set $\Ext(G,N)_{[S]}$ of $N$-extensions 
of $G$ corresponding to $[S]$ in terms of abelian Lie group extensions. 
The main point is that $S$ defines a Lie group extension
$$ \1 \to N_{\rm ad} := N/Z(N) \to G^S \ssmapright{q_S} G \to \1 $$
and that for each $N$-extension $\hat G$ of $G$ corresponding to $S$, we have 
$\hat G/Z(N) \cong G^S$, so that we may consider $\hat G$ as an extension of
$G^S$ by $Z(N)$. This relation reduces several problems concerning general 
Lie group extensions to the abelian case for which we refer to [Ne04a]. 

Pick $\omega \in C^2_{ss}(G,N)$ with $\delta_S = C_N \circ \omega$. 
Since $N_{\rm ad} \cong C_N(N) \subeq \Aut(N)$, we write the elements of
this group as $C_N(n)$ and not as cosets $nZ(N)$.
On the product set 
$$ G^S := N_{\rm ad} \times G$$
we define the product 
$$ \eqalign{(C_N(n),g)(C_N(n'),g') 
&:= 
\Big( C_N\big( n S(g)(n') \omega(g,g')\big), gg'\Big) 
= \Big(C_N(n)c_{S(g)}(C_N(n')) \delta_S(g,g'), gg'\Big).\cr} $$
Note that the second form of the product shows in particular 
that it does not depend on $\omega$. We observe that 
$$ S_1 : G \to \Aut(N_{\rm ad}), \quad S_1(g)(C_N(n)) = C_N(S(g).n) $$
is a smooth outer action with 
$\delta_{S_1} = C_{N_{\rm ad}} \circ \delta_S$ if we consider 
$\delta_S$ as $N_{\rm ad}$-valued function. 

\Lemma IV.1. $G^S$ is a group with the following properties: 
\litem{(1)} It carries a natural Lie group structure for which 
$q_S \: G^S \to G, (C_N(n),g) \mapsto g $
is a Lie group extension of $G$ by $N_{\rm ad}$. 
The equivalence class of this extension only depends on the $G$-kernel $[S]$. 
\litem{(2)} The map 
$\rho \: G^S \to \Aut(N), (C_N(n), g) \mapsto C_N(n) S(g)$
defines on $N$ a smooth $G^S$-action. 
\litem{(3)} The map 
$\alpha \: N  \to G^S, n \mapsto (C_N(n),\1)$
defines a smooth crossed module with 
$$ \ker \alpha = Z(N) \quad \hbox{ and } \quad 
\coker(\alpha) = G^S/N_{\rm ad} \cong G, $$
and $S$ is the corresponding strongly smooth outer action of $G \cong G^S/N_{\rm ad}$ on $N$. 

\Proof. (1) First we have to show that $(S_1, \delta_S)$ is a smooth
factor system for $G$ and $N_{\rm ad}$: Since we define the Lie group structure on $N_{\rm ad} 
= C_N(N)$ in such a way that $C_N \: N \to N_{\rm ad}$ is the quotient map, 
the relation $\delta_S = C_N \circ \omega$ implies  
$\delta_S \in C^2_{ss}(G,N_{\rm ad})$, and we have
$\delta_{S_1} = C_{N_{\rm ad}} \circ \delta_S$. 

If $U \subeq G$ is an open identity neighborhood such that 
the map $U \times N \to N$, $(g,n) \mapsto S(g).n$ is smooth, then the
map 
$$ U \times N_{\rm ad} \to N_{\rm ad}, \quad (g,C_N(n)) \mapsto C_N(S(g).n) =
S(g)C_N(n) S(g)^{-1} = S_1(g)(C_N(n)) $$
is also smooth. 

The local smoothness of $\delta_S = C_N \circ \omega$ follows from the local
smoothness of $\omega$. Moreover, 
$$ \delta_{S,g}(g') = \delta_S(g,g') \delta_S(gg'g^{-1},g)^{-1} 
= C_N\big(\omega(g,g')\omega(gg'g^{-1},g)^{-1}\big) 
= C_N\big(\omega_{g}(g')) $$
is smooth in an identity neighborhood because $\omega \in C^2_{ss}(G,N)$. 

For $g,g',g'' \in G$ we have 
$$ \eqalign{ 
&\ \ \ \ S_1(g)\big(\delta_S(g',g'')\big) \delta_S(g,g'g'') \cr
&= S(g) \Big( S(g')S(g'')S(g'g'')^{-1}\Big) S(g)^{-1}
S(g)S(g'g'')S(gg'g'')^{-1} \cr
&= S(g) S(g')S(g'')S(gg'g'')^{-1} = S(g)S(g')S(gg')^{-1} S(gg')S(g'')S(gg'g'')^{-1} \cr
&= \delta_S(g,g') \delta_S(gg',g'').\cr} $$

This shows that $(S_1, \delta_S)$ is a smooth factor system, so that 
$G^S$ is a group and Proposition~II.8 provides a natural Lie group
structure for which $q_S \: G^S \to G$ is a Lie group extension of $G$
by $N_{\rm ad}$. 

To see that the Lie group $G^S = N_{\rm ad} \times_{(S_1, \delta_S)} G$ depends, 
as an extension of $G$ by $N_{\rm ad}$, only on the equivalence class $[S]$, let 
$\alpha \in C^1_s(G,N)$ and $S' := (C_N \circ \alpha) \cdot S$. 
Then 
$$ S_1' := (C_{N_{\rm ad}} \circ h) \cdot S_1 = h.S_1 \quad \hbox{ for } \quad 
h := C_N \circ \alpha \in C^1_s(G,N_{\rm ad}), $$
and formula (1.1) after Definition~I.1 leads to 
$$ \delta_{S'}(g,g') 
= C_N\Big(\alpha(g) S(g)\big(\alpha(g')\big)\Big) \delta_S(g,g') C_N(\alpha(gg')^{-1}) 
= (h *_{S_1} \delta_S)(g,g'). $$
According to Proposition~II.10, the map 
$$ \psi \: G^{S'} \to G^{S}, \quad 
(C_N(n),g) \mapsto (C_N(n \alpha(g)), g) = (C_N(n) h(g),g) $$
is an equivalence of extensions of $G$ by $N_{\rm ad}$. 

 (2) First we show that $\rho$ is a representation: 
$$ \eqalign{ 
\rho(C_N(n),g)\rho(C_N(n'),g')) 
&= C_N(n) S(g) C_N(n') S(g') 
= C_N(n) C_N(S(g).n') S(g) S(g') \cr
&= C_N\big(n (S(g).n')\big) \delta_S(g,g') S(gg') 
= \rho\Big( C_N\big(n (S(g).n')\big) \delta_S(g,g'), gg'\Big)\cr
&= \rho( (C_N(n),g)(C_N(n'),g')). \cr} $$
Since $\rho(G^S)$ consists of smooth automorphisms of $N$, the
smoothness of the representation follows from the local smoothness and 
the smoothness of the orbit maps. Local smoothness is a consequence of
the local smoothness of the map 
$$ (N_{\rm ad} \times G) \times N \to N, \quad ((C_N(n), S(g)),n') \to C_N(n)
S(g).n' $$
which, in addition, shows that the action of $G^S$ on $N$ is a map 
$G^S \times N \to N$ that is smooth on a set of the form 
$U^S \times N$, where $U^S$ is an identity neighborhood in $G^S$. 
Since $G^S$ acts on $N$ by smooth maps, this implies the 
smoothness of the action. 

 (3) is an immediate consequence of (1) and (2). 
\qed

If the group $G$ is connected, then we can use Lemma~II.9, 
instead of Proposition~II.8, in the preceding proof 
and draw the same conclusions if $S$ is not strongly 
smooth. 

\Lemma IV.2. The map 
$$ \psi = (\rho,q_S) \: G^S \to \Aut(N) \times G $$
is injective and yields an isomorphism of groups 
$$ G^S \cong \{(\phi,g) \in \Aut(N) \times G \: S(g) \in \phi \cdot C_N(N)\}.$$ 

\Proof. Since $\ker q_S = N_{\rm ad}$ and $(\ker \rho) \cap N_{\rm ad} = \{\1\}$, the map 
$\psi$ is an injective homomorphism of groups. 

For each element $(C_N(n),g) \in G^S$ we have 
$\psi(C_N(n),g) = (C_N(n) S(g), g),$
which proves ``$\subeq$,'' and for any pair 
$(\phi,g) \in \Aut(N) \times G$ with $\phi \in S(g)\cdot C_N(N)$ we find an element 
$n \in N$ with $\phi = C_N(n)S(g)$, which means that 
$(\phi,g) = \psi(C_N(n),g)$. This proves the lemma. 
\qed

\Proposition IV.3. Let $q \: \hat G \to G$ be a Lie group extension of 
$G$ by $N$, corresponding to the smooth $G$-kernel $[S]$,  and 
$C_N$ the corresponding representation of $\hat G$ on $N$. Assume further that 
$Z(N)$ is a split Lie subgroup of $N$. Then $\hat G/Z(N) \cong G^S$ and the map 
$$ \gamma = (C_N \cdot (S \circ q)^{-1}, q) \: \hat G \to \Aut(N)\times G, \quad 
g \mapsto (C_N(g) S(q(g))^{-1}, q(g)) $$ 
defines a Lie group extension 
$$ Z(N) \into \hat G \sssmapright{\gamma} G^S. $$
This assignment has the following properties: 
\litem{(1)} If $q_j \: \hat G_j \to G$, $j =1,2$, are equivalent
extensions of $G$ by $N$, then 
$\gamma_j \: \hat G_j \to G^S$ are equivalent extensions of $G^S$ by
the smooth $G^S$-module $Z(N)$. We thus obtain a map 
$$ \tau \: \Ext(G,N)_{[S]} \to \Ext(G^S,Z(N))_{[\rho]}. $$
\litem{(2)} A Lie group extension $\gamma \: \hat G^S \to G^S$ of $G^S$ by
$Z(N)$ comes from an extension of $G$ by $N$ corresponding to $[S]$ 
 if and only if there exists a $G^S$-equivariant equivalence 
$$ \alpha \: N \to \gamma^{-1}(N_{\rm ad}) \subeq \hat G^S $$
of central extensions of $N_{\rm ad}$ by $Z(N)$. Note that $G^S \cong \hat G^S/Z(N)$ 
acts on $\gamma^{-1}(N_{\rm ad})$ because $Z(N)$ acts trivially on this group. 

\Proof. Since the extension $\hat G$ corresponds to $[S]$, it is equivalent
to an extension of the type $N \times_{(S,\omega)} G$ with 
$q(n,g) = g$, where $(S,\omega)$ is a
smooth factor system (Proposition II.8). This means that 
$$ C_N(n,g) = C_N(n) S(g)  = C_N(n) S(q(n,g)), \leqno(4.1) $$
so that 
$\gamma(n,g) = (C_N(n),g).$
Now the explicit formulas for the multiplications 
in $N \times_{(S,\omega)} G$ and $G^S$
imply that $\gamma$ is a surjective morphism of Lie groups. 
Its kernel is $Z(N)$, so that the assumption that this is a split Lie
subgroup implies that $\gamma$ defines a Lie group 
extension of $G^S$ by $Z(N)$. 

 (1) If 
$\phi \: \hat G_1 \to \hat G_2$
is an equivalence of $N$-extensions of $G$, then 
the representations 
$C_N^j$, $j =1,2$, of $\hat G_j$ on $N$ satisfy 
$C_N^2 \circ \phi = C_N^1$
because $\phi\res_N = \id_N$. Therefore the quotient maps 
$$ \gamma_j = (C_N^j \cdot (S \circ q_j)^{-1}, q_j) \: \hat G_j  \to G^S $$
satisfy 
$\gamma_2 \circ \phi 
= ((C_N^2 \circ \phi) \cdot (S \circ q_2 \circ \phi)^{-1}, q_2 \circ \phi) 
= (C_N^1 \cdot (S \circ q_1)^{-1}, q_1) = \gamma_1.$
This means that $\phi \: \hat G_1 \to \hat G_2$ also is an equivalence of
extensions of $G^S$ by $Z(N)$. 

 (2) Suppose first that the extension 
$\gamma \: \hat G^S \to G^S$ by $Z(N)$ comes from the $N$-extension 
$q \: \hat G \to G$ corresponding to $[S]$ in the sense that 
$\hat G = \hat G^S$ with $\hat G/Z(N) \cong G^S$. 
We may assume that $\hat G = G_{(S,\omega)}$ (Proposition II.8). Then (4.1) shows that 
$(n,g) \in \hat G$ acts on $N$ by $C_N(n) S(g) =  \rho(C_N(n),g)$. 
Therefore the inclusion $N \into \hat G$ onto the subset 
$\gamma^{-1}(N_{\rm ad})$ is equivariant with respect to the 
action of $G^S$, and therefore in particular for the action of 
$N_{\rm ad}$, so that it is an equivalence of central extensions of
$N_{\rm ad}$ by $Z(N)$. 

Suppose, conversely, that 
$\gamma \: \hat G^S \to G^S$ is an extension of $G^S$ by $Z(N)$ for
which there exists a $G^S$-equivariant
equivalence 
$\alpha \: N \to \gamma^{-1}(N_{\rm ad})$
of central extensions of $N_{\rm ad}$ by $Z(N)$. 
Then 
$$ \hat G^S/\alpha(N) = \hat G^S/\gamma^{-1}(N_{\rm ad}) \cong G^S/N_{\rm ad}
\cong G, $$
so that we obtain an extension of $G$ by~$N$ by the quotient map $q = q_S \circ \gamma 
\: \hat G^S \to G$ with kernel 
$\gamma^{-1}(N_{\rm ad})$. The action of $G^S \cong \hat G^S/Z(N)$
on $N$ induced by the conjugation action of $\hat G^S$ on $N$ 
coincides with the given action 
$$ \rho \: G^S \to \Aut(N), \quad 
(C_N(n),g) \mapsto C_N(n) S(g)$$ 
 of $G^S$ on $N$ because $\alpha$ is $G^S$-equivariant. Therefore 
the $G$-kernel of the extension $q \: \hat G^S \to G$ is~$[S]$. 

For the conjugation action $C_N$ of $\hat G^S$ on $N$ we have 
$C_N = \rho \circ \gamma$ and $q = q_S \circ \gamma$, so that
the corresponding map $\hat G^S \to G^S$ coincides with
$\gamma$. This means that $\gamma \: \hat G^S\to G^S$ is associated
to the extension $q \: \hat G^S \to G$ by the process described above. 
\qed

\Remark IV.4. In Lemma~IV.1 we have seen that the map 
$\alpha \: N \to G^S$ defines a smooth crossed module with smooth outer action 
$S \: G \to \Aut(N)$. According to Theorem~III.8, the central extension 
$$ Z(N) \into N \sssmapright{\alpha} N_{\rm ad} $$
can be embedded in a $G^S$-equivariant way into an abelian extension 
$$ Z(N) \into \hat G^S \sssmapright{\alpha} G^S $$
if and only if $\chi_{ss}(S) \in H^3_{ss}(G,Z(N))_S$ vanishes. If $G$ is connected, 
we have the simpler criterion  $\chi_{s}(S) = \1$ in $H^3_{s}(G,Z(N))_S$. 
Comparing with Proposition~IV.3, we see that this happens if and only if 
$\Ext(G,N)_{[S]} \not= \eset$ (cf.\ also Theorem~II.17). 
\qed

\section{Appendix A. Automorphisms of Lie group extensions} 

\nin In this appendix we collect several useful results of automorphisms of 
group extensions. 

Let $q \: \hat G \to G$ be a Lie group extension of $G$ by $N$ 
corresponding to the $G$-kernel $[S]$. In this appendix we discuss the group 
$$\Aut(\hat G, N) := \{ \phi \in \Aut(\hat G) \: \phi(N) = N\} $$ 
of all Lie group automorphisms of $\hat G$ preserving the 
split normal Lie subgroup $N$. We have a group homomorphism 
$$ \Phi \: \Aut(\hat G, N) \to \Aut(N) \times \Aut(G), \quad 
\phi \mapsto (\phi_N, \phi_G) := (\Phi_N(\phi), \Phi_G(\phi)) $$
with 
$$ \phi_N  := \phi\res_N \quad \hbox{ and } \quad 
\phi_G \circ q = q \circ \phi  $$
(cf.\ Lemmas~I.7 and~II.3). The kernel of $\Phi$ consists of those automorphisms of $\hat G$ inducing
the identity on the subgroup $N$ and the quotient group $G$. One of the main 
results of this appendix is an exact sequence 
$$ \0 \to Z^1_s(G,Z(N))_S \sssmapright{\Psi} 
\Aut(\hat G,N) \sssmapright{\Phi} (\Aut(N) \times \Aut(G))_{[S]}
\sssmapright{I} H^2_{ss}(G,Z(N))_{S}, $$
where $I$ is a $1$-cocycle for the action of the stabilizer 
$(\Aut(N) \times \Aut(G))_{[S]}$ of $[S]$ on the group 
$H^2_{ss}(G,Z(N))_{S}$. 
In the context of abstract groups, a similar result has been obtained by 
Wells in [Wel71] (cf.\ also [Ro84]). 

\subheadline{Automorphisms of Lie group extensions} 

First we take a closer look at the action of $\Aut(N) \times \Aut(G)$ on 
factor systems and extensions. 
In the following we shall often refer to the natural action of 
$\Aut(N) \times \Aut(G)$ on the sets $C^p_s(G,N)$ by 
$$ \big((\phi, \psi).f\big)(g_1, \ldots, g_p) 
:= \phi(f(\psi^{-1}(g_1), \ldots, \psi^{-1}(g_p))). $$

\Lemma A.1. The group $\Aut(N) \times \Aut(G)$ acts naturally on 
$C^1_s(G,\Aut(N)) \times C^2_{ss}(G,N)$ by 
$$ (\phi, \psi).(S,\omega) := (c_\phi \circ S \circ \psi^{-1}, 
\phi \circ \omega \circ (\psi^{-1} \times \psi^{-1})). $$
For $(\phi,\psi).(S,\omega) = (S',\omega')$ the following assertions hold: 
\litem{(1)} If $\delta_S = C_N \circ \omega$, then $\delta_{S'} = C_N \circ \omega'$. 
\litem{(2)} $d_{S'} \omega' = \phi \circ (d_S \omega) \circ (\psi \times \psi \times \psi)^{-1}$. 
In particular, the subset $Z^2_{ss}(G,N)$ is invariant under the action of 
$\Aut(N) \times \Aut(G)$. 

\litem{(3)} For $h \in C^1_s(G,N)$ and $(\phi, \psi) \in \Aut(N) \times \Aut(G)$ we have 
$$ (\phi, \psi).(h.(S,\omega)) = \big((\phi,\psi).h\big).\big((\phi, \psi).(S,\omega)\big). $$
In particular, we obtain 
an action of the semidirect product group 
$$ C^1_s(G,N) \rtimes (\Aut(N) \times \Aut(G)) $$
and hence an action of $\Aut(N) \times \Aut(G)$ on the set 
$H^2_{ss}(G,N)$ of $C^1_s(G,N)$-orbits. 

\litem{(4)} Let ${\cal E}\: N \sssmapright{\iota} \hat G \sssmapright{q} G$ be a 
short exact sequence describing a Lie group extension of $G$ by $N$ and consider the exact 
sequence $$ {\cal E}' := (\phi,\psi).{\cal E}: \quad 
\1 \to N \smapright{\iota \circ \phi^{-1}} \hat G \smapright{\psi \circ q} G \to\1. $$
For $\hat G= N \times_{(S,\omega)} G$ the sequence ${\cal E}'$ 
describes an extension equivalent to 
$N \times_{(\phi, \psi).(S,\omega)} G$ and the map 
$$ \mu_{(\phi, \psi)} \: N \times_{(S,\omega)} G \to 
 N \times_{(\phi, \psi).(S,\omega)} G, \quad 
(n,g) \mapsto (\phi(n), \psi(g))$$ 
is an isomorphism of Lie groups. 

\Proof. (1) For $S' := (\phi, \psi).S$ and $\delta_S = C_N \circ \omega$ we have 
$$ \eqalign{ \delta_{S'}  
&= c_\phi \circ \delta_S \circ (\psi \times \psi)^{-1} 
= c_\phi \circ C_N \circ \omega \circ (\psi \times  \psi)^{-1} \cr
&= c_\phi \circ C_N \circ \phi^{-1} \circ (\phi,\psi).\omega 
= C_N \circ (\phi,\psi).\omega. \cr} $$

 (2) follows from 
$$ \eqalign{ (d_{S'} \omega')(g,g',g'') 
&= S'(g)(\omega'(g',g''))
\omega'(g,g'g'') 
\omega'(gg',g'')^{-1} \omega'(g,g')^{-1} \cr
&= \phi\Big(S(g)\big(\omega(\psi^{-1}g', \psi^{-1}g'')\big)\cdots 
\omega(\psi^{-1}g, \psi^{-1}g')^{-1}\Big)\cr 
&= \phi \circ d_S\omega \circ (\psi \times \psi \times \psi)^{-1}(g,g',g''). \cr} $$
In view of (1), this formula shows that the 
action of $\Aut(N) \times \Aut(G)$ on the product set 
$C^1_s(G, \Aut(N)) \times C^2_{ss}(G,N)$ preserves the subset $Z^2_{ss}(G,N)$ 
of smooth factor systems.  

 (3) follows from an easy calculation. 
For the identification of $\Ext(G,N)$ with the set 
$H^2_{ss}(G,N)$ of $C^1_s(G,N)$-orbits in $Z^2_{ss}(G,N)$, we refer to Proposition~II.10. 

 (4) Suppose that $\hat G = N \times_{(S,\omega)} G$ holds for a smooth 
factor system $(S,\omega)$ (Proposition~II.8). 
If $\sigma \: G \to N \times_{(S,\omega)} G, g \mapsto (\1,g)$ 
is the canonical section, then 
$\sigma' := \sigma \circ \psi^{-1} \in C^1_s(G,\hat G)$ satisfies 
$$ (\psi \circ q) \circ \sigma' = \psi \circ q \circ \sigma \circ \psi^{-1} 
= \id_G, $$
so that we can interprete $\sigma'$ as a section for the extension 
${\cal E}'$. 
We now have 
$$ \delta_{\sigma'} = \delta_\sigma \circ \psi^{-1} = \iota \circ 
\omega \circ \psi^{-1},$$ 
and the action of $\hat G$ on $N$ corresponding to the extension ${\cal E}'$ is given by 
$$ C_N'(x)(g) 
:= \iota'^{-1}(x \iota'(g)x^{-1})
= \phi \circ \iota^{-1}(x \iota(\phi^{-1}(g))x^{-1}) 
= \phi \circ C_N(x)(\phi^{-1}(g)) 
= (c_\phi \circ C_N(x))(g). $$
Therefore 
$$ S' := C_N' \circ \sigma' 
= c_\phi \circ C_N \circ \sigma \circ \psi^{-1} 
= c_\phi \circ S \circ \psi^{-1} = (\phi, \psi).S. $$
Moreover, the relation $\iota \circ \omega = \delta_\sigma$ leads to 
$$ \eqalign{  \delta_{\sigma'} 
&= \delta_\sigma \circ (\psi \circ \psi)^{-1}  
= \iota \circ \omega \circ (\psi \circ \psi)^{-1}  
= \iota \circ \phi^{-1} \circ \phi \circ \omega \circ (\psi \circ \psi)^{-1}  
= \iota' \circ (\phi,\psi).\omega.\cr} $$ 
This implies that 
$\omega' = (\phi,\psi).\omega$, so that the sequence ${\cal E}'$ is represented by the 
extension $N \times_{(S',\omega')}~G$ with $(S',\omega') = (\phi,\psi).(S,\omega)$. 
Furthermore, the map $\mu_{(\phi,\psi)}$ is a group homomorphism:  
$$ \eqalign{ \mu_{(\phi,\psi)}((n,g)(n',g')) 
&= \mu_{(\phi,\psi)}(n S(g)(n') \omega(g,g'), gg') \cr
&=  \Big(\phi(n) \big(\phi \circ S(g)(n')\big) \big(\phi \circ \omega(g,g')\big), 
\psi(gg') \Big) \cr
&=  (\phi(n) (c_\phi \cdot S(g))(\phi(n')) \omega'(\psi(g), \psi(g')), \psi(gg') ) \cr
&=  (\phi(n) S'(g)(\phi(n')) \omega'(\psi(g), \psi(g')), \psi(gg') ) \cr
&=  (\phi(n), \psi(g)) (\phi(n'), \psi(g')) 
= \mu_{(\phi,\psi)}(n,g) \mu_{(\phi,\psi)}(n',g').\cr} $$
That it is an isomorphism of Lie groups follows from the fact that 
it is a local diffeomorphism.  
\qed

The following two lemmas provide descriptions of kernel and range of~$\Phi$. 

\Lemma A.2. On $Z(N)$ we consider the action $S_Z$ of $G$ induced by the 
conjugation action of $\hat G$ on $N$. For each 
$f \in Z^1_s(G,Z(N))_{S_Z}$ 
we obtain an element $\phi_f := (f \circ q) \cdot \id_{\hat G} 
\in \ker \Phi$, and the map 
$$ \Psi \: Z^1_s(G,Z(N))_{S_Z} \to \ker \Phi, \quad f \mapsto \phi_f $$
is a bijective group homomorphism. 

\Proof. Let $f \in Z^1_s(G,Z(N))$. Clearly $\phi_f \: \hat G \to \hat G$ 
is a diffeomorphism of 
$\hat G$ with $q \circ \phi_f = q$ whose inverse is given by 
$\phi_{f^{-1}}$. That it is a group homomorphism follows for 
$\gamma = q(\hat\gamma)$ and 
$\gamma' = q(\hat\gamma')$ from 
$$ \phi_f(\hat\gamma \hat\gamma') 
= f(\gamma\gamma') \hat\gamma \hat\gamma' 
= f(\gamma) \gamma(f(\gamma')) \hat\gamma \hat\gamma' 
= f(\gamma) \hat\gamma f(\gamma') \hat\gamma^{-1} \hat\gamma \hat\gamma' 
= \phi_f(\gamma) \phi_f(\gamma'). \leqno(A.1) $$ 

From $q \circ \phi_f = q$ we immediately obtain that 
$\phi_f \circ \phi_{f'} = \phi_{ff'}$
for $f,f' \in Z^1_s(G,Z(N))$, showing that $\Psi$ 
is a group homomorphism. This homomorphism is obviously injective. 
To see that is also is surjective, let $\phi \in \ker \Phi$. 
Then $\phi = \hat f \cdot \id_{\hat G}$ with a smooth function 
$\hat f \: \hat G \to N$. The fact that $\phi$ is a homomorphism implies 
that $\hat f \in Z^1_s(\hat G,N)$ with respect to the conjugation action of 
$\hat G$ on $N$ (cf.\ (A.1)). Moreover, $\hat f(N) = \{\1\}$ implies 
that $\hat f$ is constant on the $N$-cosets: 
$$ \hat f(gn) = \hat f(g) c_{g}(\hat f(n)) = \hat f(g), \quad g \in \hat G, n \in N. $$
Hence $\hat f$ can be written as 
$f \circ q$ for some smooth function $f \: G \to N$. 
Since $N$ is normal in $\hat G$, we also have for each $g \in \hat G$ and $n \in N$ 
the relation 
$$ \hat f(g) = \hat f\big(g\cdot (g^{-1}ng)\big) = \hat f(ng) = \hat f(n) c_n(\hat f(g)) 
= c_n(\hat f(g)), $$
and hence that $\im(f) \subeq Z(N)$. This proves that 
$f \in Z^1_s(G,Z(N))$, so that $\Psi$ is surjective. 
\qed

\Proposition A.3. Let $(S,\omega) \in Z^2_{ss}(G,N)$, 
$\hat G =  N \times_{(S,\omega)} G$ be the 
corresponding Lie group extension of $G$ by $N$ and 
$(\phi, \psi) \in \Aut(N) \times \Aut(G)$. 

\par {\rm(a)} We write the extension $\hat G$ as the exact sequence 
${\cal E} \: \1\to N \sssmapright{\iota} \hat G \sssmapright{q} G\to\1. $
Then $(\phi,\psi) \in \im(\Phi)$ 
if and only $(\phi, \psi).{\cal E} \sim {\cal E}$ holds for 
the extension 
$$ (\phi,\psi).{\cal E}\:  \quad \1\to N \smapright{\iota \circ \phi^{-1}} \hat G 
\smapright{\psi \circ q} G\to\1.  $$

\par {\rm(b)} An automorphisms $\nu \in \Aut(\hat G,N)$ satisfies 
$\Phi(\nu) = (\phi, \psi)$ if and only if it is of the form 
$$ \nu(n,g) = (\phi(n) h(\beta(g)), \beta(g)) \leqno(A.2) $$
with $h \in C^1_s(G,N)$ satisfying 
$(\phi, \psi).(S,\omega) = h.(S,\omega).$
In particular, $(\phi, \psi) \in \im(\Phi)$ if and only if we have 
in $\Ext(G,N) \cong Z^2_{ss}(G,N)/C^1_s(G,N)$ the relation  
$$(\phi, \psi).[(S,\omega)] = [(S,\omega)].$$

\Proof. (a) For $\lambda \in \Aut(\hat G,N)$ we consider the extension 
${\cal E}' := (\lambda_N, \lambda_G).{\cal E}$ and put 
$\iota' := \iota \circ \lambda_N^{-1}$ and 
$q' := \lambda_G \circ q$. Then the map 
$\lambda \: \hat G \to \hat G$ yields an equivalence of extensions 
$$ \matrix{
 N & \smapright{\iota \circ \lambda_N^{-1}} & \hat G & 
\smapright{\lambda_G \circ q} & G \cr 
\mapdown{\id_N} & & \mapdown{\lambda} & & \mapdown{\id_G} \cr 
 N & \smapright{\iota} & \hat G & \smapright{q}& G. \cr } $$
Therefore $\Phi(\lambda).{\cal E} \sim {\cal E}$. 
If, conversely, $(\phi, \psi).{\cal E} \sim {\cal E}$, 
there exists an equivalence of extensions 
$$ \matrix{
 N & \sssmapright{\iota \circ \phi^{-1}} & \hat G & \sssmapright{\psi \circ q} & 
G \cr 
\mapdown{\id_N} & & \mapdown{\lambda} & & \mapdown{\id_G} \cr 
 N & \sssmapright{\iota} & \hat G & \sssmapright{q}& G. \cr } $$
This means that $\lambda_N = \phi$ and $\lambda_G = \psi$. 

(b) Let $\nu \in \Aut(\hat G, N)$ with 
$\Phi(\nu) = (\phi, \psi)$. In the product coordinates of 
$\hat G = N \times_{(S,\omega)} G$ we then have 
$$ \nu(n,g) = (\phi(n) h(\psi(g)), \psi(g)) $$
for some $h \in C^1_s(G,N)$. 
Let $(S',\omega') := (\phi, \psi).(S,\omega)$ 
and consider the isomorphism 
$$ \mu \: N \times_{(S,\omega)} G \to N \times_{(S',\omega')} G, \quad 
(n,g) \mapsto (\phi(n), \psi(g)) $$
(Lemma~A.1(4)). Then we obtain an isomorphism 
$$ \lambda := \nu \circ \mu^{-1} \: N \times_{(S',\omega')} G \to N \times_{(S,\omega)} G, 
\quad (n,g) \mapsto (n h(g), g) $$
which is an equivalence of extensions, so that Proposition~II.10 implies that 
$(S',\omega') = h.(S,\omega)$. 

If, conversely, $(S',\omega') = h.(S,\omega)$, then 
$\lambda$ is an equivalence of extensions and 
$\nu := \lambda \circ \mu$ an automorphism of Lie groups. 
\qed

\Corollary A.4. Suppose that $N$ is a smooth abelian $G$-module, where the action of $G$ on 
$N$ is given by $S \: G \to \Aut(N)$, $\omega \in Z^2_{ss}(G,N)$,  
and that $\hat G =  N \times_\omega G$ is the corresponding abelian extension of 
$G$ by $N$. Further let $(\phi, \psi) \in \Aut(N) \times \Aut(G)$. 
Then an automorphisms 
$\nu \in \Aut(\hat G,N)$ satisfies 
$\Phi(\nu) = (\phi, \psi)$ if and only if it is of the form 
$$ \nu(n,g) := (\phi(n) h(\beta(g)), \beta(g)) \leqno(A.3) $$
with $h \in C^1_s(G,N)$ satisfying 
$$ (\phi, \psi).S = S \quad \hbox{ and } \quad 
(\phi, \psi).\omega - \omega = d_S h. 
\qeddis

The preceding observations provide the exact sequence 
$$ \1 \to Z^1_s(G,Z(N)) \sssmapright{\Psi} \Aut(\hat G,N) \sssmapright{\Phi} 
\big(\Aut(N) \times \Aut(G)\big)_{[{\cal E}]} \to \1. $$

The following lemma expresses that the actions of $\hat G$ on $N$ and $G$ preserve 
the extension class ${\cal E}$. 

\Lemma A.5. Let $N$ be a $G$-Lie group and 
$q_G \: \hat G \to G$ be an $N$-extension of $G$ for which the corresponding outer 
action is contained in the class $[S]$. Then for each 
$g \in G$ there is an automorphism 
$\phi \in \Aut(\hat G,N)$ with $\phi_N = S(g)$ and 
$\phi_G = c_g$. 

\Proof. Let $\hat g \in \hat G$ with $q(\hat g) = g$. Then 
$\psi := C_{\hat G}(\hat g) \in \Aut(\hat G,N)$ satisfies 
$\psi_G = c_g$ and $\psi_N$ is of the form $c_n \circ S(g)$ for 
some $n \in N$ because 
$[\hat G] \in \Ext(G,N)_{[S]}$. On the other hand 
$C_{\hat G}(n)$  satisfies 
$C_{\hat G}(n)_G = \id_G$ and $C_{\hat G}(n)_N = c_n$, so that 
$\phi := C_{\hat G}(n)^{-1} \circ \psi$ meets all requirements. 
\qed

\Remark A.6. The group $N$ also acts on the set $Z^2_{ss}(G,N)$ of factor systems 
by the homomorphism $(C_N, \1) \: N \to \Aut(N) \times \Aut(G)$. 
This means that 
$$ (n.S)(g) = C_N(n) \cdot S(g) \cdot C_N(n)^{-1} = C_N\big(n S(g)(n^{-1})\big) S(g) 
\quad \hbox{ and } \quad 
n.\omega = C_N(n) \circ \omega. $$
We consider the function $h \in C^1_s(G,N)$ given by 
$h(g) := n S(g)(n^{-1}) = d_S(n)(g)$ and claim that 
$$ h *_S \omega = n.\omega, $$
which implies that $n.(S,\omega) = h.(S,\omega)$. 
In fact, we have 
$$ \eqalign{ 
(h *_S \omega)(g,g') 
&= h(g) S(g)(h(g')) \omega(g,g') h(gg')^{-1} \cr
&= n S(g)(n)^{-1} S(g)(n S(g')(n^{-1})) \omega(g,g') S(gg')(n) n^{-1}\cr
&= n S(g)S(g')(n^{-1}) \omega(g,g') S(gg')(n) n^{-1}\cr
&= n \big(C_N(\omega(g,g')) S(gg')(n^{-1})\big) \omega(g,g') S(gg')(n) n^{-1}\cr
&= n \omega(g,g') S(gg')(n^{-1}) \omega(g,g')^{-1} \omega(g,g')S(gg')(n) n^{-1}\cr
&= n \omega(g,g') n^{-1} = C_N(n)(\omega(g,g')).\cr} $$
\qed

\Remark A.7. For an extension $q \: \hat G \to G$ it is also interesting 
to consider its {\it gauge automorphism group} 
$$ \Gau(\hat G) := \{ \phi \in \Aut(\hat G,N) \: q \circ \phi  = q\} $$
of all automorphisms of $\hat G$ inducing the identity on $G$. 
It is the kernel of the natural homomorphism 
$$ \Phi_G \: \Aut(\hat G,N) \to \Aut(G) $$
and it obviously contains 
$\ker \Phi \cong Z^1_s(G,Z(N))$ (Proposition~A.3). 

Each automorphism of $\hat G$ preserving $N$ that induces the identity on 
$G$ can be written in the form 
$\phi = \phi_f := f \cdot \id_{\hat G}$, 
where $f \: \hat G \to N$ is a smooth 
function. We then have 
$$ \phi_{f_1} \circ \phi_{f_2} = \phi_{(f_1 \circ \phi_{f_2}) \cdot f_2}. $$
If $f\res_N \cdot \id_N \: N \to N$ is a diffeomorphism, then the map 
$\phi_f$ is a diffeomorphism of $\hat G$ with 
$\phi_f^{-1} = \phi_{f^{-1} \circ \phi_f^{-1}}.$ 

It is easy to see that $\phi_f$ is a group homomorphism if and only 
$f \in Z^1_s(\hat G,N)_{C_N}$ with respect to the conjugation action 
$C_N$ of $\hat G$ on $N$. This means that 
$$ \Gau(\hat G) \cong (Z^1_s(\hat G,N),*)^\times $$
is the unit group of the monoid $Z^1_s(\hat G,N)$, 
where the monoid structure on $Z^1_s(\hat G,N)_{C_N}$ is given by 
$\phi_{f_1 * f_2} = \phi_{f_1} \circ \phi_{f_2}$
and 
$$ \eqalign{ (f_1 * f_2)(\hat\gamma) 
&= f_1(f_2(\hat\gamma)\hat\gamma) f_2(\hat\gamma) 
= f_1(f_2(\hat\gamma)) c_{f_2(\hat\gamma)}(f_1(\hat\gamma)) f_2(\hat\gamma) 
= f_1(f_2(\hat\gamma)) f_2(\hat\gamma) f_1(\hat\gamma). \cr} $$
This formula shows in which way the pointwise product on 
the group $C^1_s(\hat G,N)$ is twisted to obtain the monoid structure on 
$Z^1_s(\hat G,N)$. 
\qed

Instead of considering the subgroup 
$\im(\Phi) = (\Aut(N) \times \Aut(G))_{[{\cal E}]},$
one may also consider the larger subgroup 
$$ \Comp(S) := (\Aut(N) \times \Aut(G))_{[S]} 
= \{ (\phi, \psi) \in \Aut(N) \times \Aut(G) \: (\phi, \psi).S \sim S\} $$
of {\it $S$-compatible pairs of automorphisms of $N$ and $G$}, 
and then identify $\im(\Phi)$ as a subgroup of $\Comp(S)$. The  group 
$\Comp(S)$ acts on $\Ext(G, N)_{[S]}$ because its action on 
$\Ext(G,N)$ preserves the subset $\Ext(G,N)_{[S]} \subeq \Ext(G,N)$. 
We also have a simply transitive action of 
the abelian group $H^2_{ss}(G,Z(N))$ on 
$\Ext(G, N)_{[S]}$ induced from the action of $Z^2_{ss}(G,Z(N))$ on 
$Z^2_{ss}(G,N)$ by 
$$ \beta.(S,\omega) := (S, \beta \cdot \omega) $$
(Corollary~II.13). 
For $(\phi, \psi) \in \Comp(S)$ we then have 
$$ (\phi,\psi).(\beta.(S,\omega)) 
=  (\phi,\psi).(S, \beta \cdot \omega) 
=  \big((\phi,\psi).S, ((\phi, \psi).\beta) \cdot (\phi, \psi).\omega\big). $$
This implies that the actions of $\Comp(S)$ and $H^2_{ss}(G,Z(N))$ on 
$\Ext(G,N)_{[S]}$ satisfy 
$$  (\phi,\psi).([\beta].[\hat G]) 
=  ((\phi,\psi).[\beta]).((\phi,\psi).[\hat G]). $$

\Proposition A.8. For $(\phi, \psi) \in \Comp(S)$ let 
$I(\phi,\psi) \in H^2_{ss}(G, Z(N))_S$ denote the unique cohomology class with 
$ (\phi, \psi).[\hat G] = I(\phi, \psi).[\hat G].$ 
Then the following assertions hold: 
\litem{(1)} The action of $\Comp(S)$ on 
$\Ext(G,N)_{[S]} = H^2_{ss}(G,Z(N))_S.[\hat G]$ has the form 
$$  (\phi,\psi).([\beta].[\hat G]) 
=  ((\phi,\psi).[\beta] + I(\phi,\psi)).[\hat G]. $$
\litem{(2)} $I \: \Comp(S) \to H^2_{ss}(G,Z(N))_S$ is a $1$-cocycle 
for the action of $\Comp(S)$ on $H^2_{ss}(G,Z(N))_S$. 
\litem{(3)} If $\hat G \cong N \times_{(S,\omega)} G$, then there exists 
an $h_0 \in C^1_s(G,N)$ with $(\phi, \psi).S = h_0.S$ and 
$$ I(\phi, \psi) = [(\phi, \psi).\omega - h_0 *_S \omega]. $$
\litem{(4)} If $I(\phi, \psi) = 0$, then there exists some 
$\gamma \in C^1_s(G,Z(N))$ with $d_S \omega = (\phi,\psi).\omega  - h_0 *_S \omega$, and for 
$h := h_0 \cdot \gamma$ the map 
$$ \nu \: N \times_{(S,\omega)} G \to N \times_{(S,\omega)} G, \quad 
(n,g) \mapsto (\phi(n) h(\psi(g)), \psi(g)) $$
is an element of $\Aut(\hat G,N)$ with $\Phi(\nu) = (\phi, \psi)$. 

\Proof. (1) is clear. 

 (2) For $(\phi, \psi), (\phi', \psi') \in \Comp(S)$ we have 
$$ \eqalign{ I((\phi, \psi)(\phi', \psi')).[\hat G]
&= ((\phi, \psi)(\phi', \psi')).[\hat G]
=  (\phi, \psi).((\phi', \psi').[\hat G]) \cr
&=  (\phi, \psi).(I(\phi', \psi').[\hat G]) 
=  ((\phi, \psi).I(\phi', \psi')).((\phi, \psi).[\hat G]) \cr
&=  ((\phi, \psi).I(\phi', \psi')).(I(\phi, \psi).[\hat G]) 
=  ((\phi, \psi).I(\phi', \psi') + I(\phi, \psi)).[\hat G]. \cr} $$

 (3) As $(\phi,\psi) \in \Comp(S)$, we have 
$S' := (\phi,\psi).S \sim S$, so that there exists an $h_0 \in C^1_s(G,N)$ 
with $S' = h_0.S = (C_N \circ h_0) \cdot S$. 
Let $\omega' := (\phi,\psi).\omega$. Then 
$$ C_N \circ \omega' = \delta_{S'} = C_N \circ (h_0 *_S \omega) $$
implies that 
$\beta := \omega' - h_0 *_S \omega \in Z^2_{ss}(G,Z(N))_S,$
and we have 
$$ (\phi,\psi).(S,\omega) = (S',\omega') 
= \beta.(S', h_0 *_S \omega)
= \beta.(h_0.S, h_0 *_S \omega)
\sim \beta.(S, \omega). $$
Therefore the cocycle $I$ is given by 
$$ I(\phi, \psi) = [\beta] = [\omega' - h_0 *_S \omega] 
= [(\phi,\psi).\omega - h_0 *_S \omega] \in H^2_{ss}(G,Z(N))_S. $$

 (4) If $I(\phi, \psi) = \1$, then there exists some $\gamma \in C^1_s(G,Z(N))$ with 
$\beta = d_S \gamma$, so that we obtain for 
$h := h_0 \cdot \gamma$ the relations $S' = h.S$ and 
$\omega' = h_0 *_S \omega + d_S \gamma = h *_S \omega,$
i.e., $(S',\omega') = h.(S,\omega)$. Then 
the map 
$$ \lambda_h \: N \times_{(S',\omega')} G \to N \times_{(S,\omega)} G, \quad 
(n,g) \mapsto (n h(g), g) $$
is an equivalence of extensions (Proposition~II.10), and composition with the isomorphism 
$$ \mu_{(\phi,\psi)} \: N \times_{(S,\omega)} G \to N \times_{(S',\omega')} G, \quad 
(n,g) \mapsto (\phi(n), \psi(g)) $$
(Lemma~A.1) leads to the automorphism $\nu$. 
\qed

We now see that  the exact sequence from above has a prolongation:

\Theorem A.9. The following sequence is exact 
$$ \1 \to Z^1_s(G,Z(N))_S \to \Aut(\hat G,N) \sssmapright{\Phi} 
\Comp(S) \sssmapright{I} H^2_{ss}(G,Z(N))_S. $$

\Proof. This follows from Lemma A.2 and Propositions~A.3 and A.8, which yields 
$$\im(\Phi) = \Comp(S)_{[\hat G]} = I^{-1}(0). 
\qeddis 

\subheadline{Automorphism groups of Lie group extensions} 

So far we dealt with pairs 
$(\phi, \psi) \in \Aut(N) \times \Aut(G)$ coming from automorphisms of 
the extension $q \: \hat G \to G$ with kernel $N$, which lead to the exact sequence 
$$ \1\to 
Z^1_s(G,Z(N))_S \into \Aut(\hat G, N) \onto (\Aut(N) \times \Aut(G))_{[\hat G]}\to\1. \leqno(A.4) $$
If we are given a group homomorphism 
$$ \psi \: H \to 
(\Aut(N) \times \Aut(G))_{[\hat G]} 
= (\Aut(N) \times \Aut(G))_{[S,\omega]} = \im(\Phi), $$
then each element $\psi(h)$ comes from an automorphism 
$\hat\psi(h) \in \Aut(\hat G,N)$, but this automorphism is not unique if 
$Z^1_s(G,Z(N))_S \not=\1$. To obtain a group homomorphism 
$\hat\psi \: H \to \Aut(\hat G,N)$ lifting $\psi$, a certain cohomology class 
in $H^2(H, Z^1_s(G, Z(N))_S)$ has to vanish, i.e., 
the abelian extension $\psi^* \Aut(\hat G,N)$ of $H$ by 
$Z^1_s(G,Z(N))_S$ has to split. That this is not always the case is shown in the 
following example, where we describe a central extension 
$Z \into \hat G \onto G$ for which (A.4) does not split. 

\Example A.10. (a) We construct an example of a central extension 
$ Z\into \hat G = Z \times_f G \onto G$ of finite-dimensional real Lie
groups, for which the sequence 
$$ Z^1_s(G,Z) \cong \Hom(G,Z) \into \Aut(\hat G,Z) \onto \big(\Aut(Z) \times
\Aut(G))_{[f]} $$ 
is not split. As this is an exact sequence of finite-dimensional real
Lie groups, it suffices to show that the corresponding sequence 
$$ Z^1(\g,\z) \cong \Hom_{\rm Lie}(\g,\z) \into \der(\hat \g,\z) \onto
\big(\der(\z) \times \der(\g))_{[\omega]} $$ 
of Lie algebras 
is not split if $\omega \in Z^2(\g,\z)$ satisfies $\hat\g\sim \z \oplus_\omega \g$. 

Let $G$ be the $3$-dimensional real Heisenberg group and $\g$
is Lie algebra. Then $\g$ has a basis of the form 
$p,q,z$ satisfying 
$$ [p,q] = z, \quad [p,z] = [q,z] = 0 $$
and we identify $G$ with $\g$, endowed with the Campbell--Hausdorff
product 
$$ x * y := x + y + {1\over 2}[x,y]. $$
Note that all left multiplications $\lambda_x(y) := x * y$ are affine maps.

We consirder the cocycle $\omega \in Z^2(\g,\R)$ with 
$$ \omega(p,z) = 1, \quad \omega(q,z) = \omega(q,z) = 0. $$
We then obtain a central extension $\hat\g := \R \oplus_\omega \g$ of $\g$ by 
$\z := \R$, and in [Ne04b, Ex.~A.8] it is shown that the exact sequence 
$$ \Hom_{\rm Lie}(\hat\g,\z) \cong \Hom_{\rm Lie}(\g,\z) 
\into \der(\hat\g,\z) \onto (\der(\g) \times \der(\z))_{[\omega]} $$
is not split. More concretely,  the
action of $\b := \Hom_{\rm Lie}(\g,\z(\g)) \subeq \der(\g)$ preserves the class
$[\omega] \in H^2(\g,\z)$, but the action of $\b$ on $\g$ does not
lift to an action of the abelian Lie algebra $\b$ on $\hat\g$. 

Let $\hat\b := \g \rtimes \b$ be the semidirect sum. Then 
$[\omega] \in H^2(\g,\z)^\b$, but there is no representation
$\hat S$ of $\b$ on $\hat\g$ lifting the representation on $\g$. 
\qed

\Lemma A.11. Suppose that $\hat G = N \times_{(S,\omega)} G$. Further let 
$\theta \: H \to C^1_s(G,N)$ be a map with 
$$ h.(S,\omega) = \theta(h)^{-1}.(S,\omega)  \quad \hbox{ for all } \quad h \in H $$
and define the automorphism 
$\hat\psi(h) \in \Aut(\hat G,N)$ by 
$$ \hat\psi(h)(n,g) := (h.n \cdot \theta(h)(h.g)^{-1}, h.g). $$
Then $\hat \psi \: H \to \Aut(\hat G,N)$ is a homomorphism if and only if 
$\theta \in  Z^1(H, C^1_s(G,N)_\psi),$
where $C^1_s(G,N)_\psi$ stands for the group $C^1_s(G,N)$, endowed with the 
$H$-action given by $(h.\chi):= \psi(h).\chi$. 

\Proof. That the maps $\hat\psi(h)$ define elements of 
$\Aut(\hat G,N)$ follows from Proposition~A.3(b). We have 
$$ \eqalign{ \hat\psi(h)\hat\psi(h')(n,g) 
&= \hat\psi(h)(h'.n \cdot \theta(h')(h'.g)^{-1}, h'.g) \cr
&= \Big(h.(h'.n \cdot  \theta(h')(h'.g)^{-1}) \cdot \theta(h)(hh'.g)^{-1}, hh'.g\Big)\cr
&= \Big(hh'.n \cdot h.(\theta(h')(h'.g)^{-1}) \cdot \theta(h)(hh'.g)^{-1}, hh'.g\Big) \cr
&= \Big(hh'.n \cdot (\theta(h) \cdot h.\theta(h'))(hh'.g)^{-1}, hh'.g\Big). \cr} $$
This implies that 
$$  \delta_{\hat\psi}(h,h') 
= \hat\psi(h) \hat\psi(h') \hat\psi(hh')^{-1} \in Z^1_s(G,Z(N)) = \ker \Phi 
\subeq \Aut(\hat G,N).$$
In view of 
$$ \hat\psi(h)^{-1}(n,g) 
= (h^{-1}.n \cdot h^{-1}.(\theta(h)(g)), h^{-1}.g) 
= (h^{-1}.n \cdot (h^{-1}.\theta(h))(h^{-1}.g), h^{-1}.g), $$
the cocycle $\delta_{\hat\psi}$ is given by 
$$  \theta(h) \cdot h.\theta(h') \cdot hh'.((hh')^{-1}.\theta(hh'))^{-1}
=  \theta(h) \cdot h.\theta(h') \cdot \theta(hh')^{-1}
= (d_H\theta)(h, h'). $$
We conclude that $\hat\psi$ is a homomorphism if and only if $\theta$ is a cocycle. 
\qed

\Remark A.12. For a given homomorphism $\psi \: H \to (\Aut(N) \times
\Aut(G))_{[(S,\omega)]}$ we may choose a map $\theta \: H \to C^1_s(G,N)$ with 
$h.(S,\omega) = \theta(h)^{-1}.(S,\omega)$ for each $h \in H$. 
We have seen in the preceding proof that 
$$ \im(d_H\theta) \subeq Z^1_s(G,Z(N)) $$
and that for the corresponding map $\hat\psi \: H \to \Aut(\hat G,N)$ 
we have 
$\delta_{\hat\psi} = d_H \theta.$
This is the cocycle of the abelian extension 
$\psi^*\Aut(\hat G,N)$ of $H$ by $Z^1_s(G,Z(N)) = \ker \Phi$ obtained by pulling back the 
abelian extension 
$$ Z^1_s(G, Z(N)) \into \Aut(\hat G,N) \to (\Aut(N) \times \Aut(G))_{[(S,\omega)]}. $$
We may thus assign to $\psi$ the cohomology class 
$$ [d_H\theta] = [\psi^* \Aut(\hat G,N)] \in H^2(H,Z^1_s(G,Z(N))_\psi). $$

If the cohomology class $[d_H\theta]$ vanishes, then there exists a
map $\eta \: H \to Z^1_s(G,Z(N))$ with $d_H \theta = d_H \eta$, which
implies that $\tilde\theta := \theta \cdot \eta^{-1}$ satisfies $d_H \tilde\theta = 0$ and, in
addition, $h.(S,\omega) = \tilde\theta(h)^{-1}.(S,\omega)$ for each
$h \in H$. Therefore $\tilde\theta \in Z^1(H,C^1_s(G,N)_\psi)$ defines 
a lift $\hat\psi$ of~$\psi$. Hence $\psi$ lifts to a
homomorphism $\hat\psi \: H \to \Aut(\hat G,N)$ if and only if the
cohomology class $[d_H \theta]$ vanishes, which means that there
exists a cocycle $\theta \: H \to C^1_s(G,N)$ with 
$h.(S,\omega) = \theta(h)^{-1}.(S,\omega)$ for each $h \in H$. 
\qed

In addition to the setting of the preceding subsection, we now assume that  
$H$ is a Lie group and ask for smoothness properties of actions on $\hat G$. 

\Proposition A.13. Suppose 
we are given smooth actions of $H$ on $G$ and $N$ by automorphisms 
and that we have an action of $H$ on $\hat G := N \times_{(S,\omega)} G$ given by 
$\theta \in Z^1(H, C^1_s(G,N))$
via 
$$ h.(n,g) = (h.n \cdot \theta(h)(h.g)^{-1}, h.g). $$
Then the following assertions hold: 
\litem{(1)} Each $h \in H$ acts by a Lie group automorphism.
\litem{(2)} The action map $\sigma \: H \times \hat G \to \hat G$ is smooth in 
an identity neighborhood of $H \times \hat G$ if and only if the map 
$$ \tilde\theta \: H \times G \to N, \quad (h,g) \mapsto \theta(h)(g) $$
is smooth in an identity neighborhood. 
\litem{(3)} If $N$ is abelian, then the 
action has smooth orbit maps if and only if all maps 
$$ \theta_g \: H \to N, \quad h \mapsto \theta(h)(h.g)^{-1} \cdot \omega(h.g, g^{-1}) $$
are smooth in an identity neighborhood of $H$. 
\litem{(4)} The action $\sigma$ of $H$ on $\hat G$ is smooth if and only if 
$\sigma$ is smooth in an identity neighborhood of $H \times\hat G$ 
and all orbit maps are smooth. If $G$ is connected, then 
the smoothness of the orbit maps follows from the local smoothness of the 
action. 

\Proof. (1) Each $h \in H$ acts by a group automorphism which is smooth in an 
identity neighborhood, hence a Lie group automorphism. 

 (2) Suppose that the action is smooth in an identity neighborhood. 
Then the function 
$(h,g) \mapsto \theta(h)(h.g)$ has this property and the map 
$$ H \times  G\to  H \times G, \quad (h,g) \mapsto (h,h^{-1}.g) $$
is a diffeomorphism fixing $(\1,\1)$. Therefore 
$\tilde\theta$  is smooth on an identity neighborhood. 
The converse is clear. 

 (3) Since $H$ acts by Lie group automorphisms, the orbit maps are smooth if and 
only if they are smooth in an identity neighborhood of $H$. 

It suffices to consider orbit maps of 
elements of the form $(\1,g)$ because $H$ acts smoothly on $N$. 
For $g \in G$ the orbit map is smooth in an identity neighborhood of $H$ if and only 
if the map 
$$ \eqalign{ h \mapsto \big(h.(\1,g)\big)(\1,g)^{-1} 
&= (\theta(h)(h.g)^{-1}, h.g)( \omega(g^{-1},g)^{-1},g^{-1}) \cr
&= \Big(\theta(h)(h.g)^{-1} S(h.g)\big(\omega(g^{-1},g)^{-1}\big) 
\cdot \omega (h.g,g^{-1}), (h.g)g^{-1}\Big) \cr}$$
(cf.\ Lemma~II.8 for the inversion formula) 
is smooth in an identity neighborhood, which is equivalent to the smoothness of 
$$  h \mapsto 
\theta(h)(h.g)^{-1} \cdot 
S(h.g)\big(\omega(g^{-1},g)^{-1}\big) \cdot \omega (h.g,g^{-1})$$
in an identity neighborhood. 

If, in addition, $N$ is abelian, then $S$ defines a smooth 
action of $G$ on $N$, 
so that  $S(h.g)\big(\omega(g^{-1},g)^{-1}\big)$ is a smooth function of $h$. 
Hence it suffices that the map 
$$  h \mapsto 
\theta(h)(h.g)^{-1} \cdot \omega (h.g,g^{-1})$$
is smooth. 

 (4) Since $H$ acts by smooth automorphisms of $\hat G$, 
the action of $H$ is smooth if and only if all orbit maps are smooth and 
it is locally smooth in the sense of (2). 

The set of all elements of $\hat G$ with smooth orbit maps is a subgroup 
containing $N$. If the action is locally smooth, then this subgroup also contains 
an identity neighborhood, hence the group $N \times_{(S_0, \omega_0)} G_0$, where 
$G_0$ is the identity component of $G$, $S_0 := S\res_{G_0}$ and 
$\omega_0 := \omega\res_{G_0 \times G_0}$. If $G$ is connected, this argument shows 
that the local smoothness of the action already implies the smoothness of 
all orbit maps. 
\qed

\section{Appendix B. Lie group structures on groups} 

\Theorem B.1. Let $G$ be a group and 
$K = K^{-1}$ a symmetric subset. 
We further assume that $K$ is a smooth manifold such that 
\litem{(L1)} there exists an open $\1$-neighborhood $V \subeq
K$ with $V^2 = V \cdot V \subeq K$ such that the group multiplication 
$\mu_V \: V \times V \to K$ is smooth, 
\litem{(L2)} the inversion map $\eta_K \: K \to K, k \mapsto k^{-1}$ is
smooth, and 
\litem{(L3)} for each $g \in G$ there exists an open $\1$-neighborhood $K_g \subeq
K$ with $c_g(K_g) \subeq K$ and such that the conjugation map 
$$ c_g \: K_g \to K, \quad x \mapsto gxg^{-1} $$
is smooth. 

Then there exists a unique structure of a
Lie group on $G$ for which there exists an open $\1$-neighborhood 
$K_1 \subeq K$ such that the inclusion map $K_1 \to G$ induces a 
diffeomorphism onto an open subset of $G$. 

If $G$ is generated by each $\1$-neighborhood $U \subeq K$, then condition {\rm(L3)} can be
omitted.

\Proof. [Ne04a, Theorem A.4 and Remark A.5]
\qed

\def\entries{


\[Baer34 Baer, R., {\it Erweiterungen von Gruppen und ihren Isomorphismen}, Math. Zeit. 
{\bf 38} (1934), 375--416 


\[Bo93 Borovoi, M. V., {\it Abelianization of the second nonabelian 
Galois cohomology}, Duke Math. Journal {\bf 72:1} (1993), 217--239 

\[Br71 Brown, L.~G., {\it Extensions of topological groups}, Pac. J. Math. {\bf 39:1} (1971), 
71--78 

\[CFNW94 Cederwall, M., Ferretti, G., Nilsson, B.E.W., and A. Westerberg, {\it 
Higher-di\-men\-sio\-nal loop algebras, non-abelian extensions and p-branes}, 
Nuclear Phys. B  {\bf 424:1}  (1994), 97--123

\[EML42 Eilenberg, S. and S. MacLane, {\it Group extensions and homology}, Ann. of Math. 
{\bf 43:4}  (1942), 757--831 

\[EML47 ---, {\it Cohomology theory in abstract groups. II. 
Group extensions with a non-Abelian kernel}, Ann. of Math. (2)  {\bf 48}  (1947), 326--341

\[Gl01 Gl\"ockner, H., {\it Infinite-dimensional Lie groups without completeness 
condition}, in ``Geometry and Analysis on Finite-
and Infinite-Dimensional Lie Groups,'' A.~Strasburger et al Eds., 
Banach Center Publications {\bf 55}, Warszawa 2002; 53--59 

\[Go69 Goto, M., {\it On an arcwise connected subgroup of a Lie group}, Procedings of the Amer. 
Math. Soc. {\bf 20} (1969), 157--162  



\[Hu81 Huebschmann, J., {\it Automorphisms of group extensions and differentials in 
the Lyndon--Hochschild--Serre spectral sequence}, J. Alg. {\bf 72:2} (1981), 296--334 


\[KM97 Kriegl, A., and P.\ Michor, ``The Convenient Setting of
Global Analysis,'' Math.\ Surveys and Monographs {\bf 53}, Amer.\
Math.\ Soc., 1997 

\[MacL63 MacLane, S., ``Homological Algebra,'' Springer-Verlag, 1963 

\[Mil58 Milnor, J., {\it On the existence of a connection with curvature zero}, Comment. Math. 
Helv. {\bf 32} (1958), 215--223 

\[Mil83 ---, {\it Remarks on infinite-dimensional Lie groups},
Proc. Summer School on Quantum Gravity, B. DeWitt ed., Les Houches, 1983

\[Ne02 Neeb, K.-H., {\it Central extensions of infinite-dimensional
Lie groups}, Annales de l'Inst. Fourier 52:5 (2002), 1365--1442 

\[Ne04a ---, {\it Abelian extensions of infinite-dimensional Lie
groups}, Travaux Math., to appear 

\[Ne04b ---, {\it Non-abelian extensions of topological Lie algebras}, 
Commun. in Algebra, to appear 

\[Ne05 ---, {\it Exact sequences for Lie group cohomology with non-abelian coefficients}, 
in preparation 

\[RSW00 Raeburn, I., Sims, A., and D. P. Williams, {\it Twisted actions and obstructions 
in group cohomology}, in ``$C^*$-algebras,'' (M\"unster, 1999), Springer, Berlin, 2000; 
 161--181 

\[Ro84  Robinson, D. J. S., {\it Automorphisms of group extensions}, in
``Algebra and its Applications,'' International Symp. on Algebra and
Its Applications in New Delhi 1981, Marcel Dekker, 1984, New York,
Lecture Notes in Pure and Applied Math. {\bf 91}; 163--167  


\[Sch26a Schreier, O., {\it \"Uber die Erweiterungen von Gruppen I},  
Monatshefte f. Math. {\bf 34} (1926), 165--180 

\[Sch26b ---, {\it \"Uber die Erweiterungen von Gruppen II},  
Abhandlungen Hamburg {\bf 4} (1926), 321--346 

\[Tu38 Turing, A.~M., {\it The extensions of a group}, Compos. Math. {\bf 5} (1938), 357--367 

\[Wei95 Weibel, C. A., ``An Introduction to Homological Algebra,''
Cambridge studies in advanced math. {\bf 38}, Cambridge Univ. Press,
1995 

\[We71 Wells, Ch., {\it Automorphisms of group extensions},
Transactions of the Amer. Math. Soc. {\bf 155:1} (1971), 189--194 

} 

\references 
\lastpage

\bye